%% file: nonstationary_graphical_models.tex
\documentclass[12pt]{article}
\usepackage{graphicx, fullpage}
\usepackage{amsmath, amssymb,theorem,verbatim}
\usepackage[round]{natbib} 
\usepackage{chapterbib}
\usepackage{xcolor}
\usepackage{setspace}
\usepackage{caption}
\usepackage{subcaption}
\usepackage{hyperref}
\usepackage{listings}

\newcommand{\shortarrow}[1][3pt]{\mathrel{%
   \vcenter{\hbox{\rule[-.5\fontdimen8\textfont3]{#1}{\fontdimen8\textfont3}}}%
   \mkern-4mu\hbox{\usefont{U}{lasy}{m}{n}\symbol{41}}}}

\usepackage{amsfonts} 
\DeclareMathSymbol{\shortminus}{\mathbin}{AMSa}{"39}

\newtheorem{theorem}{Theorem}[section]

\newtheorem{assumption}{Assumption}[section]
\newtheorem{defin}{Definition}[section]
\newtheorem{lemma}{Lemma}[section]
\newtheorem{example}{Example}[section]
\newtheorem{corollary}{Corollary}[section]
\newtheorem{proposition}{Proposition}[section]
\newtheorem{remark}{Remark}[section]

\newcommand{\cov}{\mathop{\rm {\mathbb C}ov}\nolimits}%
\newcommand{\var}{\mathop{\rm {\mathbb V}ar}\nolimits}%
\newcommand{\cor}{\mathop{\rm {\mathbb C}orr}\nolimits}

\newcommand{\mmod}{\mathrm{mod}}

\newcommand{\Ex}{\mathrm{E}}

\newcommand{\colred}[1]{{\color{red} #1}}

\title{Graphical models for nonstationary time series}
\author{{Sumanta Basu\footnote{Department of Statistics and Data Science,
Cornell University \texttt{sumbose@cornell.edu}} and Suhasini Subba
    Rao\footnote{Department of Statistics, Texas A\&M University \texttt{suhasini@stat.tamu.edu}}}}
\date{\today}

\begin{document}

\maketitle


\begin{abstract}
  We propose NonStGM, a general nonparametric graphical modeling
  framework for studying dynamic associations among the components of
  a nonstationary multivariate time series. It builds on the framework
  of Gaussian Graphical Models (GGM) and stationary time series
 Graphical models (StGM),
  and complements existing works on parametric graphical models
  based on change point vector autoregressions (VAR).
Analogous to StGM, the proposed framework 
captures conditional noncorrelations (both intertemporal and
contemporaneous) in the form of an undirected
graph.
In addition, to describe the more nuanced nonstationary relationships
among the components of the time series, we introduce the
new notion of  \textit{conditional nonstationarity/stationarity} and
incorporate it within the graph architecture.  
This allows one to distinguish between direct and indirect nonstationary relationships among system
  components, and can be used to search for small subnetworks that
  serve as the ``source'' of nonstationarity in a large system.
Together, the two concepts of conditional noncorrelation and nonstationarity/stationarity
  provide a  parsimonious description of the dependence structure of the time series. 

In GGM, the graphical model structure is encoded in the sparsity pattern of the
inverse covariance matrix. Analogously, we
explicitly connect conditional noncorrelation and stationarity 
between and within components of the multivariate time series to zero and 
Toeplitz embeddings of an infinite-dimensional inverse covariance
operator. In order to learn the graph, we move to the Fourier domain.
We show that in the Fourier domain, conditional
stationarity and noncorrelation relationships in the inverse covariance operator are
encoded with a specific sparsity
structure of its integral kernel operator. Within the local stationary
framework we show that these sparsity patterns can be recovered from
finite-length time series by node-wise regression of discrete
Fourier Transforms (DFT) across different Fourier frequencies.
We illustrate the features of our general framework under the special
case of time-varying 
Vector Autoregressive models. We demonstrate the feasibility of
learning  NonStGM structure from data using simulation studies.

\vspace{0.5em}


\end{abstract}

\input{section1}

\input{section2_Revision}

\input{section3_Revision}

\input{section4_Revision}

\input{section5_Revision}

\input{section6_Revision}

\input{conclusions_acknowledgements}

\subsection*{Acknowledgements}

SB and SSR  acknowledge the partial support of the National
Science Foundation (grants DMS-1812054 and DMS-1812128).
In addition, SB acknowledges partial support from the National Institute of
Health (grants R01GM135926 and R21NS120227).
The authors thank Gregory Berkolaiko for several useful suggestions
and Jonas Krampe for careful reading. The authors thank the Associate
Editor and two anonymous referees for their thoughtful comments
and suggestions which substantially improved the paper.

\bibliography{bib_tsreg,bib_sb}

\newpage

\appendix


\input{Appendix_Section2}

\input{Appendix_Section3}

\input{Appendix_Section4}

\input{Appendix_Section5}

\input{Appendix_example}

\input{Appendix_stationary}

\bibliographystyle{plainnat}

\end{document}

%% file: section1.tex
\section{Introduction}\label{sec:intro}

Graphical modeling of multivariate time series has received
considerable attention in the past decade as a tool to study
dynamic relationships among the components of a large system
observed over time. Key applications include, among others,
analysis of brain networks in neuroscience \cite{lurie2020questions}
and understanding linkages among firms for measuring systemic
risk buildup in financial markets \cite{diebold2014network}.

The vast majority of graphical models for time series
focuses on the stationary setting (see \cite{p:bri-96, p:dah-eic-97,
  p:dah-00b, p:eic-04, p:eich-07,
  bohm2009shrinkage, jung2015graphical, basu2015regularized, 
  davis2016sparse, p:zha-17, qiu2016joint,
  sun2018large, fiecas2019spectral, p:sac-20}, to name but a
few). While the assumption of stationarity may be realistic in many
situations, it is well known that \textit{nonstationarity} arises in many applications. In
neuroscience, for example, task based fMRI data sets are known
to exhibit considerable nonstationarity in  the network connections,
a phenomenon known as \textit{dynamic functional connectivity}, see
\cite{preti2017dynamic}. A naive application
of graphical modeling methods designed for stationary processes
can lead to spurious network edges if the actual time series
is nonstationary.

The limited body of work on graphical models for nonstationary
time series has so far focused on a restricted class of nonstationary
models, where the data generating process can be well approximated
by a finite order change point vector autoregressive (VAR)
model. Within this framework,
\cite{wang2019localizing} and \cite{safikhani2020joint} 
have proposed methods for constructing
a ``dynamically changing'' network at each of the estimated change
points. However, these methods are designed for time
  series which are piece-wise stationary and follow a finite order stationary
  VAR over each segment. For many data sets, these conditions can be too
restrictive, for example they do not allow for smoothly changing parameters. Analogous to stationary
time series where spectral methods allow for a nonparametric approach, it would be useful to define meaningful
networks for nonstationary time series.

The objective of this paper is to move away from semi-parametric models and propose a general
framework for the graphical modeling of multivariate (say, $p$-dimensional)
nonstationary time series. Our motivation comes from Gaussian graphical
models (GGM), where the edges of a conditional dependence graph 
can distinguish between
the direct and indirect nature of dependence in multivariate
Gaussian random vectors.
We argue that a general graphical model
framework for nonstationary time series should have the capability
to distinguish between  two types of nonstationarity; the source
of nonstationarity and one that inherits their nonstationarity by
way of its connection with the source. This way of dimension
reduction will be useful for modeling large systems where the
nonstationarity arises only from a small subset of the process
and then permeates through the entire system. Moreover, the
identification of sources and propagation channels of
nonstationarity may also be of scientific interest.

Analogous to GGM, in our framework, the
presence/absence of edges in the network encodes conditional
correlation/non-correlation relationships amongst the $p$ components
(nodes) of the time series. An additional attribute distinguishes
between the types of nonstationarity. A graphical model is built using 
{\it conditional} relations. In this spirit, we introduce the
concept of {\it conditional stationarity} and {\it nonstationarity}. 
To the best of our knowledge this is a new notion. A solid edge
between two nodes in the network implies that their linear relationship, conditional 
on all the other nodes, does not change over time. In contrast, a dashed edge implies
 that their conditional relationship changes over time. We formalize these notions
 in Section \ref{sec:framework}. Nodes in our network also have
{\it self-loops} to indicate whether the time series is nonstationary
on its own, or if it inherits nonstationarity from some other component
in the system. The self loops are denoted by a circle (solid or dashed)
round the node. 

The time-varying Autoregressive model is often used to
  model nonstationarity. To illustrate the above ideas, in the following example we 
connect the parameters of a time-varying Autoregressive model (tvVAR), 
which is a mixture of constant and time dependent parameters, to the concepts
introduced above.


\vspace{1mm} 

\begin{figure}[h!]
\begin{subfigure}{0.55\textwidth}
\centering  
\includegraphics[width=1\textwidth]{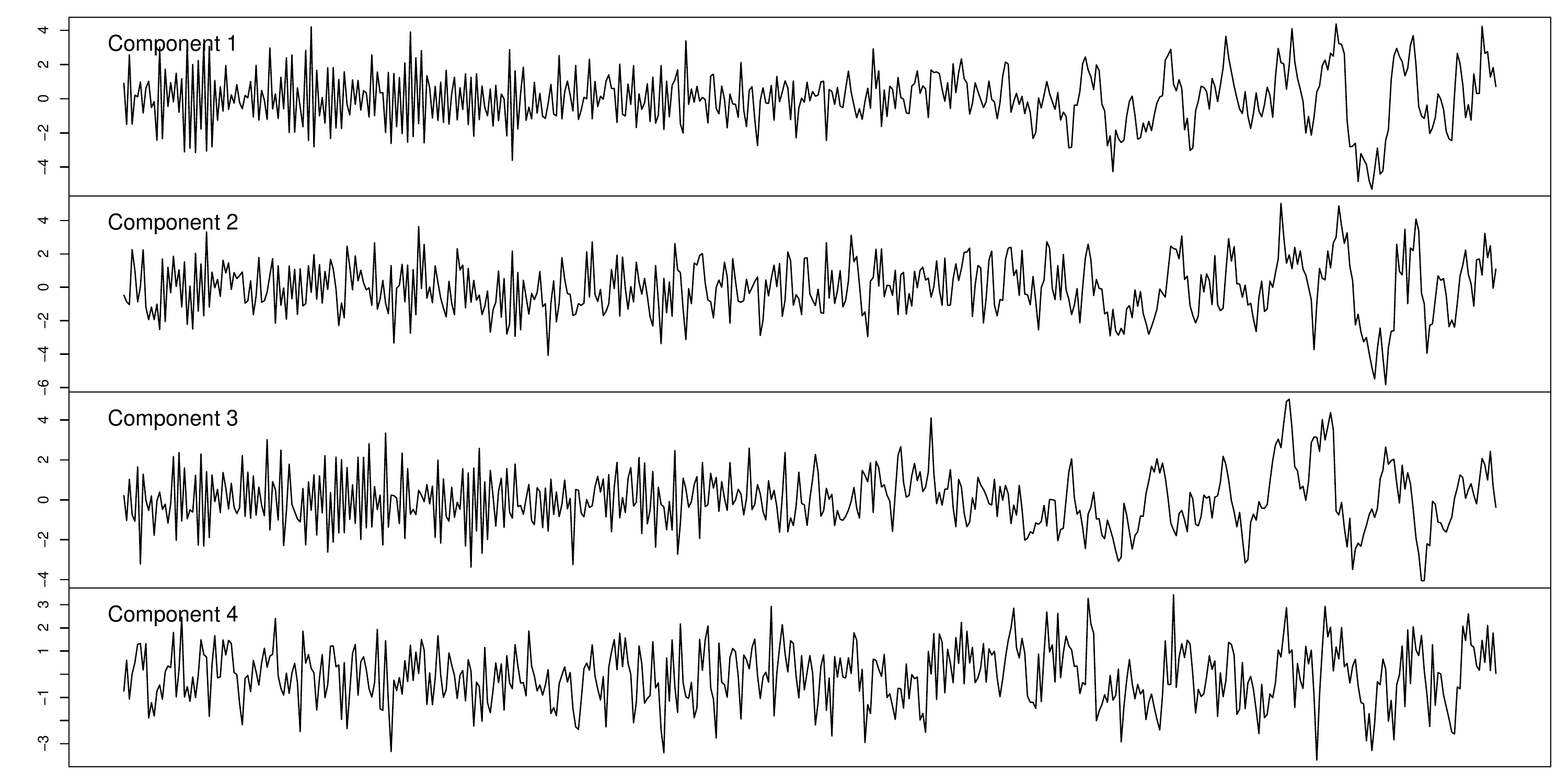}
\caption{Trajectories of a $4$-dimensional time
  series generated by a tvVAR(1) model. 
The  multivariate process is jointly nonstationary. However, components $1$ and
    $3$ are the source of nonstationarity, while the other two components
inherit the nonstationarity by means of their conditional dependence
structure. \label{fig:example_p_4a}}
\end{subfigure}
\hfill 
\begin{subfigure}{0.4\textwidth}
\centering
 \includegraphics[width = 0.7\textwidth]{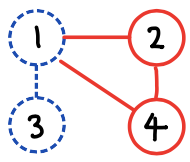}
\caption{The joint governing architecture is described in the
  graph. Dashed edges and self-loops represent conditional
  nonstationarity, while solid edges and self-loops represent
  conditional invariance and stationarity, notions of which we formalize in our new
  graphical modeling framework. \label{fig:example_p_4b}}
\end{subfigure}
\caption{Time series and conditional dependence graph of a time-varying VAR model.\label{fig:example_p_4}}
\end{figure}

\noindent 
    {\bf Toy Example }  Consider the trajectories of a $4$-dimensional time
    series given in Figure \ref{fig:example_p_4}.  The time series plots of all the components
    exhibit negative autocorrelation at the
    start of the time series that slowly changes to positive autocorrelation
    towards the end. Thus the nonstationarity of each individual time
    series, at least from a visual inspection, is  apparent. 
The data is generated from a time-varying vector autoregressive
 (tvVAR$(1)$) model (see Section \ref{sec:framework} for details),
 where components $1$ and $3$ are the \textit{sources of nonstationarity},
i.e. they are affected by their own past through a (smoothly)
time-varying parameter. In addition, component $3$ affects component $1$.
Component $2$ and $4$ affect each other in a time-invariant way. Component
$2$ is also affected by $1$ and component $1$ and $4$ are affected
through $2$.  As a result, components
$2$ and $4$  \textit{inherit the nonstationarity from the sources} $1$
and $3$.  As far as we are aware, there currently does not exist
 tools that adequately describe the 
nuanced differences in their dependencies and nonstationarity.
Our aim in this paper is to capture these relationships in the form of the
schematic diagram in Figure \ref{fig:example_p_4b}. We note that the tvVAR model is a special case of our general
 framework, which does not make any explicit assumptions on the data generating process.

It is interesting
to contrast the networks constructed using the ``dynamically
changing'' approach developed in \cite{wang2019localizing} and 
\cite{safikhani2020joint} 
for change point VAR models with our approach.  Both networks convey
different information about the nonstationary time series. The 
``dynamically changing'' network can be considered as local in the
sense that it identifies regions of stationarity and constructs a
directed graph over each 
of the stationary periods. While the graph in our approach is undirected and
yields global information about relationships between the nodes.

In order to connect the proposed framework to the current literature, 
we conclude this section by briefly reviewing the existing graphical modeling frameworks for
Gaussian random vectors (GGM) and multivariate stationary time series
(StGM). In Section \ref{sec:framework} we lay the
foundations for our nonstationary graphical models
(NonStGM) approach. In particular, we formally define the notions of conditional
noncorrelation and stationarity of nodes, edges, and subgraphs in
terms of zero and Toeplitz embeddings of an infinite dimensional inverse
covariance operator. We show that this framework offers a natural
generalizations to existing notions of
 conditional noncorrelation in GGM and StGM. 
It should be emphasized that we do not assume that the
  underlying time series is Gaussian. All the relationships that we
  describe are in terms of the partial covariance and therefore apply
  to any multivariate  time series whose covariance exists.
In Section \ref{sec:fourier}, we switch to Fourier domain and show
that the conditional noncorrelation and nonstationarity relationships
are explicitly encoded in the sparsity pattern of the integral kernel
of the inverse covariance operator. 
This connection opens the door to learning the graph structure from
finite length time series data with the discrete Fourier transforms (DFT). In Section 
\ref{sec:finitelocal} we focus on locally stationary time series.  We
show that by conducting nodewise regression of discrete Fourier transforms (DFT) of
the multivariate time series across different Fourier frequencies it
is possible to learn the network.
Section \ref{sec:tvVAR} describes how the proposed 
general framework looks in the special case of tvVAR models, where the
notions of conditional noncorrelation and nonstationarity are
transparent in the transition matrix.
Some numerical results are presented in Section \ref{sec:simulation} to 
illustrate the methodology. All the proofs for the results in this
paper can be found in the Appendix. 
 
\vspace{1mm}

\noindent {\bf Background. } We outline some relevant works in
graphical models and tests for stationarity that underpin the
technical development of NonStGM. 

\vspace{1mm}
\noindent \underline{\textit{Graphical Models.}}  A graphical
model describes the relationships among the components of
a $p$-dimensional system in the form of a graph with a set of vertices
$V = \{1, 2, \ldots, p\}$, and an edge set $E \subseteq V \times V$
containing pairs of system components which exhibit strong
association even after conditioning on the other components.


The focus of GGM is on the
conditional independence relationships in a 
$p$-dimensional (centered) Gaussian random vector $\underline{X} =
(X^{(1)}, X^{(2)}, \ldots, X^{(p)})^\top$. The non-zero
partial correlations 
$\rho^{(a,b)}$, defined as $\cor \left(X^{(a)}, X^{(b)}|X^{-\{a, b\}} \right)$ and 
also encoded in the sparsity structure of the \textit{
precision matrix} $\Theta = \left[\var(\underline{X})\right]^{-1}$, are
used to define the edge set $E$.
The task of \textit{graphical model selection}, i.e. learning the
edge set $E$ from finite sample,  is accomplished by estimating
$\Theta$ with a penalized likelihood estimator as in graphical Lasso
(\cite{friedman2008sparse}), or by nodewise regression (\cite{p:mei-bue-06})
where each component of the random vector is regressed on the other
$(p-1)$ components. 

Switching to the time series setting, consider 
$\{ \underline{X}_t = (X^{(1)}_t, \ldots, X^{(a)}_t, \ldots, X^{(p)}_t)^{\top} \}_{t \in   \mathbb{Z}}$, 
a $p$-dimensional time series
with autocovariance function $\cov(\underline{X}_t,
\underline{X}_\tau) = {\bf C}(t,\tau)$. Note that in future we usually
use $\{\underline{X}_{t}\}$ to denote the sequence $\{\underline{X}_{t}\}_{t\in \mathbb{Z}}$.
A direct adaptation of the GGM framework that estimates the \textit{contemporaneous}
precision matrix ${\bf C}^{-1}(0, 0)$ (see \cite{p:zha-17,
  qiu2016joint}) does not
provide conditional relationships between the \textit{entire time series}.
\cite{p:bri-96}  and \cite{p:dah-00b} laid the foundation of graphical models
in stationary time series, where the conditional relationships 
between the {\it entire} time series
$\{X^{(a)}_t\}$ and $\{X^{(b)}_t \}$ is captured. They show that the inverse of the
\textit{multivariate spectral density function}
$\Sigma(\omega):=
(1/2\pi) \sum_{\ell = -\infty}^\infty {\bf C}(\ell) \exp[-i\ell
  \omega]$ $\omega \in [0, \pi]$
explicitly encodes the conditional uncorrelated relationships. To be precise,
 $[\Sigma^{-1}(\omega)]_{a,b} = 0$ for all $\omega \in [0,\pi]$
if and only if $\{X^{(a)}_t\}$ and $\{X^{(b)}_t \}$
are conditionally uncorrelated, given \textit{all  the other time
  series}. The graphical model selection problem reduces to finding
  all pairs $(a,b)$ where $[\Sigma^{-1}(\omega)]_{a,b} \neq 0$ for some
  $\omega \in [0,\pi]$. For Gaussian time series the graph is a
  conditional independence graph while for non-Gaussian time series
  the graph encodes partial correlation information.
  For brevity, we refer to this approach as StGM (stationary time
  series graphical models).
 Estimation of $\Sigma^{-1}(\omega)$ is typically done
using the Discrete Fourier transform of the time series (see \cite{p:eic-08}). 
More recently, for relatively ``large'' $p$, penalized methods such as
GLASSO \citep{jung2015graphical} and CLIME \citep{fiecas2019spectral}
have been used to estimate $\Sigma^{-1}(\omega)$. 
This framework crucially relies on stationarity, in particular, the
Toeplitz property of the autocovariance function ${\bf C}_{t, \tau} =
{\bf C}(t - \tau)$, and is not immediately generalizable to the
nonstationary case.

\noindent \underline{\textit{Testing for stationarity.}}  There is a
rich literature on testing for nonstationarity of a time
series. 
Most methods are based on testing for invariance of the
spectral density function or autocovariance function over time (see \cite{p:pri-sub-69},
 \cite{p:pap-09}, \cite{p:nas-13} to name but a few). 
An alternative approach is based on the fact that the 
Discrete Fourier Transform at certain frequencies is close to uncorrelated for stationary time
 series.  \cite{p:eph-01}, \cite{p:dwi-sub-11}, \cite{p:jen-sub-15}, 
\cite{p:aue-20} use this property to test for nonzero correlation between
DFTs of different frequencies to detect for departures from stationarity.
The above mentioned tests focus on the ``marginal''  
notion of nonstationarity  instead of the conditional notion defined in this paper. 
Tests for marginal nonstationarity are not equipped to delineate between
direct and indirect nature of conditionally nonstationary relationships
among the components of a multivariate time series. However, in this
paper, we show that analogous to marginal tests, it is possible to
utilize the Fourier domain to detect for different types of conditional (non)stationarity.

%% file: section2_Revision.tex
\section{Graphical models and conditional
  stationarity}\label{sec:framework}

For a $p$-dimensional nonstationary time series $\{\underline{X}_t\}$, all the \textit{pairwise covariance}
information are contained in the infinite set of $p \times p$ autocovariance
matrices ${\bf C}_{t, \tau} = \cov[\underline{X}_t, \underline{X}_\tau]$, for $t, \tau \in \mathbb{Z}$.
We aggregate this information into an operator $C$, and show that its inverse
operator $D$ captures meaningful \textit{conditional (partial) covariance} relationships (Section \ref{sec:covariance}).
Leveraging this connection, we first define a graphical model, and conditional stationarity of
its nodes, edges and subgraphs, in terms of the operator $D$ (Section \ref{sec:network}).
Then we show that these notions can be viewed as natural generalizations of the GGM and StGM
frameworks (Sections \ref{sec:covD} and \ref{sec:covtimeD}). We start by introducing some notation that will be
used to formally define these structures (this can be skipped on first
reading). Let ${\bf A}=(A_{a,b}:1\leq a\leq d_1,1\leq b \leq d_2)$ denote a
$d_1\times d_2$-dimensional matrix, then we define $\|{\bf A}\|_{2}^{2} =
\sum_{a,b}^{}|A_{a,b}|^{2}$, $\|{\bf A}\|_{1} =
\sum_{a,b}^{}|A_{a,b}|$ and $\|{\bf A}\|_{\infty} = \sup_{ a,b}|A_{a,b}|$.

\subsection{Definitions and notation}

We use $\ell_2$ and $\ell_{2,p}$ to denote the sequence space 
$\{u=(\ldots,u_{-1},u_{0},u_{1},\ldots)'
; u_{j}\in \mathbb{C}\textrm{ and }\sum_{j}|u_{j}|^{2}<\infty\}$
 and the (column) sequence space 
$\{w={\text{vec}[u^{(1)},\ldots,u^{(p)}]};u^{(s)}\in \ell_{2}
\textrm{ for all }1\leq s \leq p\}$ respectively ($\text{vec}$ denotes the
vectorisation of a matrix).
On the spaces $\ell_{2}$ and $\ell_{2,p}$ we define the two inner products
$\langle u,v \rangle = \sum_{j\in \mathbb{Z}}u_{j}v_{j}^{*}$ (where
$*$ denotes the complex conjugate), for  $u,v\in \ell_{2}$ and $\langle x,y \rangle = \sum_{s=1}^{p}\langle u^{(s)},v^{(s)} \rangle$ for $x=(u^{(1)},\ldots,u^{(p)})^{\prime},y=(v^{(1)},\ldots,v^{(p)})\in \ell_{2,p}$,
such that $\ell_{2}$ and $\ell_{2,p}$ are two Hilbert spaces. 
For $x\in \ell_{2,p}$ let $\|x\|_{2} = \langle x,x\rangle$. 
For $s_1, s_2 \in \mathbb{Z}$, we use
$A_{s_1,s_2}$ to denote the $(s_1,s_2)$ entry in the matrix $A$, which can be
infinite dimensional and involve negative indices.

We consider the $p$-dimensional real-valued time series 
$\{\underline{X}_{t}\}_{t \in \mathbb{Z}}$, 
$\underline{X}_{t} =
(X_{t}^{(1)},\ldots,X_{t}^{(p)})^{\prime}$, where the univariate random variables $X_t^{(a)}, a = 1, \ldots, p$, are defined on the probability space
$(\Omega,\mathcal{F},P)$.
 We assume for all $t$ that
$\Ex[\underline{X}_{t}] = 0$ this condition is not
  necessary in Sections \ref{sec:framework} and \ref{sec:fourier}, but
it simplifies the exposition.
Let $L^{2}(\Omega,\mathcal{F},P)$
denote all univariate random variables $X$
where $\var[X]<\infty$, and for any $X,Y\in
L^{2}(\Omega,\mathcal{F},P)$ we define the inner product $\langle
X,Y\rangle  = \cov[X,Y]$. For every $t, \tau \in \mathbb{Z}$, we define the $p \times p$ covariance  
${\bf C}_{t,\tau} =\cov[\underline{X}_{t},\underline{X}_{\tau}]$ 
and assume $\sup_{t\in  \mathbb{Z}}\|{\bf
  C}_{t,t}\|_{\infty}<\infty$.  Under this
assumption, for all $t\in \mathbb{Z}$ and  $1\leq c\leq p$
$X_{t}^{(c)}\in L^{2}(\Omega,\mathcal{F},P)$. Let 
$\mathcal{H} = \overline{\textrm{sp}}(X_{t}^{(c)};t\in \mathbb{Z},1\leq
c\leq p)\subset L^{2}(\Omega,\mathcal{F},P)$ 
be the closure of the space spanned by $(X_{t}^{(c)};t\in \mathbb{Z},1\leq
c\leq p)$.  Since
$L^{2}(\Omega,\mathcal{F},P)$ defines a Hilbert space,
$\mathcal{H}$ is also a Hilbert space. Therefore, by the projection theorem,
for any closed subspace $\mathcal{M}$  of $\mathcal{H}$, there is a unique projection of $Y\in \mathcal{H}$
onto $\mathcal{M}$ which minimises  $\Ex(Y -
X)^{2}$ over all ${X\in \mathcal{M}}$ (see Theorem 2.3.1, \cite{b:bro-dav-06}).
We will use $P_{\mathcal{M}}(Y)$ to denote this projection. 
In this paper, we will primarily use the following subspaces 
\begin{eqnarray*}
\mathcal{H}  - X_{t}^{(a)} &=& \overline{\textrm{sp}}[X_{s}^{(c)};s\in \mathbb{Z},1\leq
c\leq p, (s,c)\neq (t,a)],\\
\mathcal{H}  - (X_{t}^{(a)},X_{\tau}^{(b)}) &=& \overline{\textrm{sp}}[X_{s}^{(c)};s\in \mathbb{Z},1\leq
c\leq p,(s,c) \notin \{(t,a), (\tau,b) \}], \\
\mathcal{H} - (X^{(c)};c\in \mathcal{S}) &=&
                                             \overline{\textrm{sp}}[X_{s}^{(c)};s\in \mathbb{Z},c\in \mathcal{S}^{\prime}],
\end{eqnarray*}
where $\mathcal{S}^{\prime}$ denotes the complement of $\mathcal{S}$.

Using the covariance ${\bf C}_{t,\tau}$ we define the infinite dimensional matrix operator $C$ 
as $C = (C_{a,b};a,b\in
\{1,\ldots,p\})$ where $C_{a,b}$ denotes an infinite dimensional submatrix  with entries 
$[C_{a,b}]_{t,\tau} =  [{\bf C}_{t,\tau}]_{a,b}$ for all $t,\tau\in
\mathbb{Z}$. For any $u\in
\ell_{2}$,
we define the (column) sequence
$C_{a,b}u =
\{[C_{a,b}u]_{t};t\in \mathbb{Z}\}$ where
 $[C_{a,b}u]_{t} = \sum_{\tau\in
   \mathbb{Z}}[C_{a,b}]_{t,\tau}u_{\tau}$. For any   
$v = \text{vec}[u^{(1)},\ldots,u^{(p)}] \in \ell_{2,p}$ we define the (column) sequence $Cv$ as 
\begin{eqnarray}
\label{eq:CcovarianceOp}
Cv = 
\left(
\begin{array}{cccc}
C_{1,1} & C_{1,2} &\ldots & C_{1,p} \\
C_{2,1} & C_{2,2} &\ldots & C_{2,p} \\
\vdots & \vdots &\ddots & \vdots \\
C_{p,1} & C_{p,2} &\ldots & C_{p,p} \\
\end{array}
\right)
\left(
\begin{array}{c}
u^{(1)} \\
u^{(2)} \\
\vdots \\
u^{(p)} \\ 
\end{array}
\right)
=
\left(
\begin{array}{c}
\sum_{s=1}^{p}C_{1,s}u^{(s)} \\
\sum_{s=1}^{p}C_{2,s}u^{(s)} \\
\vdots \\
\sum_{s=1}^{p}C_{p,s}u^{(s)} \\ 
\end{array}
\right).
\end{eqnarray}
An infinite dimensional matrix operator, $B$, is said to be zero, if
all its entries are zero. An infinite dimensional matrix operator
$A$ is said to be Toeplitz
if its entries satisfy $A_{t,\tau} = a_{t-\tau}$ for all
$t,\tau\in \mathbb{Z}$ and for some sequence $\{a_{r};r\in
\mathbb{Z}\}$.

\subsection{Covariance and inverse covariance operators}\label{sec:covariance}

Within the nonstationary framework we
require the following assumptions on $C$ to show that $C$ is a mapping from
$\ell_{2,p}$ to $\ell_{2,p}$ (and later that $C^{-1}$ is a mapping from
$\ell_{2,p}$ to $\ell_{2,p}$). For stationary time series
analogous assumptions are often made on the spectral density function
(see Remark \ref{remark:assum}). 

\begin{assumption}
\label{assum:lambda}
Define 
$\lambda_{\sup} =
  \sup_{v \in \ell_{2,p}, \|v\|_{2}=1}\langle v, Cv \rangle$,   
$\lambda_{\inf} =
  \inf_{v \in \ell_{2,p}, \|v\|_{2}=1}\langle v, Cv \rangle$.
Then 
\begin{eqnarray}
\label{assum:eigenvalues}
0 <\lambda_{\inf} \leq \lambda_{\sup} <\infty.
\end{eqnarray}
\end{assumption}
Assumption \ref{assum:lambda} implies that 
  $\sup_{t}\sup_{a}\sum_{\tau\in \mathbb{Z}}\sum_{b=1}^{p}[{\bf
    C}_{t,\tau}]_{a,b}^{2}<\infty$ and also the coefficients of
  the inverse are square summable. It can be shown that if 
$\sup_{t\in \mathbb{Z}}\sum_{\tau\in \mathbb{Z}}\|{\bf
  C}_{t,\tau}\|_{\infty}<\infty$, then 
    $\sup_{v \in \ell_{2,p}, \|v\|_{2}=1}\langle v, Cv
\rangle<\infty$. This is analogous to a short memory condition for
stationary time series. It is worth keeping in mind
  that cointegrated time series do not
  satisfy this condition. The theory developed in Sections  \ref{sec:framework} and
\ref{sec:fourier} only require Assumption
\ref{assum:lambda}. However, to estimate the network stronger
conditions on ${\bf D}_{t,\tau}$ are required and these are stated in Section \ref{sec:finitelocal}.

Under the above assumption $C:\ell_{2,p}\rightarrow \ell_{2,p}$,
and since ${\bf C}_{t,\tau} = {\bf C}_{\tau,t}^{\prime}$, $\langle v, Cu
\rangle = \langle Cv,u \rangle$, thus $C$ is a self-adjoint, bounded
operator with  $\|C\| = \lambda_{\sup}$, where $\|\cdot\|$ denotes the
operator norm: $\|A\| = \sup_{u\in \ell_{2,p}, \|u\|_{2} = 
  1}\|Au\|_{2}$. 

\begin{remark}\label{remark:assum}
In the case of stationary time series sufficient conditions for
Assumption \ref{assum:lambda} to hold is that
the eigenvalues of the spectral
density matrix $\Sigma(\omega)$ are uniformly bounded away from zero
and away from $\infty$ overall $\omega \in [0,\pi]$ (see, for example, 
\cite{b:bro-dav-06}, Proposition 4.5.3).
\end{remark}

The core theme of GGM is to learn conditional (partial) covariances
between two variables after conditioning on a set of other variables.
These conditional
relationships can be derived from the inverse covariance matrix. Now we will define a suitable
inverse covariance operator $D = C^{-1}$ and show how its entries capture the
conditional relationships. We will define these conditional 
relations in terms of \textit{projections} with respect to the
$\ell_2$-norm, this is equivalent to the least squares regression coefficients at the population level.
We consider the projection of $X_{t}^{(a)}$ onto
$\mathcal{H}-X_{t}^{(a)}$, given by 
\begin{eqnarray}
\label{eq:projectXat}
P_{\mathcal{H}-X_{t}^{(a)}}(X_{t}^{(a)}) = \sum_{\tau\in \mathbb{Z}}\sum_{b=1}^{p}\beta_{(\tau,b)\shortarrow (t,a)}X_{\tau}^{(b)},
\end{eqnarray}
with $\beta_{(t,a)\shortarrow (t,a)}=0$ (note the coefficients
$\{\beta_{(\tau,b)\shortarrow (t,a)}\}$ are unique since $C$ is non-singular). 
Let $\sigma_{a,t}^{2} =
\Ex[X_{t}^{(a)} - P_{\mathcal{H}-X_{t}^{(a)}}(X_{t}^{(a)})]^{2}$, it
can be shown that $\sigma_{a,t}^{2}\geq \lambda_{\min}$ (see Appendix \ref{sec:covarianceproofs}).  
Analogous to finite dimensional
covariance matrices, to obtain the entries of the
inverse we use the coefficients of the projections of  $X_{t}^{(a)}$ onto
$\mathcal{H}-X_{t}^{(a)}$. For all $1\leq t,\tau\leq p$, we define the $p\times p$-dimensional
matrices ${\bf D}_{t,\tau}$ as follows
\begin{eqnarray}
\label{eq:Dttau}
[{\bf D}_{t,\tau}]_{a,b} = 
\left\{
\begin{array}{cc}
\frac{1}{\sigma_{a,t}^{2}} & a = b \textrm{ and } t = \tau \\
 -\frac{1}{\sigma_{a,t}^{2}}\beta_{(\tau,b)\shortarrow (t,a)} &
                     \textrm{ otherwise }.
\end{array}
\right.
\end{eqnarray}
Using ${\bf D}_{t,\tau}$ we define the infinite dimensional matrix 
\begin{eqnarray}
\label{eq:Dab}
D_{a,b}  = \left\{[D_{a,b}]_{t,\tau}=[{\bf D}_{t,\tau}]_{a,b};t,\tau \in \mathbb{Z}\right\}.
\end{eqnarray}
Analogous to the definition of $C$, we define 
$D = (D_{a,b};a,b\in \{1,\ldots,p\})$.

Our next lemma shows that the operator $D$ is indeed the inverse of the
  covariance operator $C$. We also state some upper bounds on its entries which will be
useful in our technical analysis.

\begin{lemma}\label{lemma:inverse}
Suppose Assumption \ref{assum:lambda} holds. 
Let $D$ be defined as in (\ref{eq:Dttau}). Then $C^{-1}=D$ 
and $\|D\|=\lambda_{\inf}^{-1}$. Further, for all $a,b\in \{1,\ldots,p\}$, 
$\|D_{a,b}\|\leq \lambda_{\inf}^{-1}$,
$\|D_{a,a}^{-1}\|\leq \lambda_{\sup}$ and $\sup_{t}\sum_{\tau\in \mathbb{Z}}\|{\bf D}_{t,\tau}\|_{2}^{2}\leq
p\lambda_{\inf}^{-2}$.
\end{lemma}
\noindent {\bf PROOF} In Appendix \ref{sec:covarianceproofs}. \hfill $\Box$

\subsection{Nonstationary graphical models (NonStGM)}\label{sec:network}

The operators $C$ and $D$ provide us with the objects needed to formally
define the edges in our network, and connect them to the notions of
conditional uncorrelatedness and conditional stationarity.

At this point, we note an important distinction between edge construction in GGM
and StGM, an issue that is crucial for generalizing graphical models to the nonstationarity
 case. In GGM, conditional uncorrelatedness between two random variables is
defined after conditioning on all the other random variables in the
system. On the other hand, 
in StGM, the conditional uncorrelatedness between two time series is
defined after conditioning on all the other time series. This leads to
two, potentially, different generalizations in the nonstationary
setup. A direct generalization of the
GGM framework would use the \textit{partial covariances} $\cov(X_{t}^{(a)}, X_{\tau}^{(b)} | \mathcal{S}_{1}^{\prime})$, where
$\mathcal{S}^{\prime}_{1} = \{X_{s}^{(c)}: (s,c) \notin \{(t,a), (\tau, b) \}
\}$. While a generalization of the StGM framework, would suggest using \textit{time series partial
  covariances} $\cov(X_{t}^{(a)}, X_{\tau}^{(b)} | \mathcal{S}_2^{\prime})$,
where $\mathcal{S}_2^{\prime} = \{X_{s}^{(c)}: s \in \mathbb{Z}, ~c
\notin \{a, b\} \}$.

To address this issue, we start by using the inverse covariance
operator $D$ to define edges that encode conditional
uncorrelatedness and (non)stationarity. We show that, as expected,
these notions are a direct generalization of the GGM framework.
Then we present a surprising result (Theorem \ref{lemma:conditional}),
that the encoding of the partial
covariances in terms of the operator $D$ remains unchanged
even if we adopt the StGM notion of partial covariance, i.e. the
conditionally uncorrelated and
conditionally (non)stationary nodes, edges, subgraphs are preserved under the
two frameworks.

We  now define the 
network corresponding to the multivariate time series. 
Each edge in our network $(V,E)$ will have 
an indicator to denote conditional invariance and  conditional time-varying, a new notion we now
introduce. The edge set $E$ will contain all pairs $(a,b)$ where 
$\{X_{t}^{(a)}\}$ and $\{X_{t}^{(b)}\}$ are conditionally correlated. 
The edge set $E$ will also contain self-loops, that convey important information about the network.  
We start by formally defining the notions of conditional
noncorrelation and (non)stationarity. This is stated in terms
of the submatrices $\{D_{a,b}\}$ of $D$.

\begin{defin}[Nonstationary network]\label{def:network} 
Conditional covariance and (non)stationarity of the components
of a $p$-dimensional nonstationary time series are represented using a
graph $G = (V, E)$, where $V = \{ 1, 2, \ldots, p \}$ is the set of
nodes, and $E \subseteq V \times V$ is a set of undirected edges 
($(a,b) \equiv (b,a)$), and includes \textit{self-loops} of the form $(a,a)$. 

\begin{itemize}
\item {\bf  Conditional Noncorrelation} 
The two time series $\{X_{t}^{(a)}\}$ and $\{X_{t}^{(b)}\}$ are conditionally uncorrelated
if $D_{a,b} =0$. As in GGM and StGM, this is represented by the absence
of an edge between nodes $a$ and $b$ in the network, i.e. $(a,b) \notin E$.
\item {\bf Conditionally Stationary Node} 
The time series $\{X_{t}^{(a)}\}$ is conditionally
stationary if $D_{a,a}$ is Toeplitz operator. We denote this using a
solid self-loop $(a,a)$ around the node $a$. 
\item  {\bf Conditionally Time-invariant Edge} 
If $a\neq b$ and $D_{a,b}$ is a Toeplitz operator, then $(a,b)$ is a
time invariant edge. 
We represent a conditionally time-invariant edge $(a,b)$ in our
network with a solid edge.   

\item {\bf Conditionally Stationary Subgraph} 
A subnetwork of nodes $\mathcal{S}\subset \{1,\ldots,p\}$ is a called
a conditionally stationary subgraph if for all $a,b\in
  \mathcal{S}$, $D_{a, b}$ are Toeplitz operators
  i.e. $D_{\mathcal{S},\mathcal{S}}$ is a block Toeplitz operator. 

As a special case of the above, we call a conditionally stationary subgraph of
order two (consisting of the nodes $\{a,b\}$) a conditionally
stationary pair if $D_{a,a},D_{a,b}$ and $D_{b,b}$ are Toeplitz.
\item {\bf Conditionally Nonstationary Node/Time-varying Edge:} (i) If $D_{a,a}$ is not
Toeplitz then $\{X_{t}^{(a)}\}$ is conditionally nonstationary. 
(ii) For $a\neq b$, if $D_{a,b}$ is not Toeplitz 
then $(a,b)$ has a conditionally time-varying edge. 

 We represent conditional nonstationary nodes using a \textit{dashed}
 self-loop and a conditionally time-varying edge with a dashed edge. 

\end{itemize}
\end{defin}
In Section \ref{sec:tvVAR} we show how the parameters of a general tvVAR
model are related to the operator $D$, and can be used to identify
the network structure in NonStGM. As a concrete example, 
below we describe the network corresponding to the tvVAR$(1)$
considered in the introduction.

\begin{example}\label{example:running}
  Consider the following tvVAR(1) model for a $4$-dimensional time
  series
\begin{eqnarray*}
\label{eq:example1}
\left(
\begin{array}{c}
X_{t}^{(1)} \\
X_{t}^{(2)} \\
X_{t}^{(3)} \\
X_{t}^{(4)} \\
\end{array}
\right) = 
\left(
\begin{array}{cccc}
\alpha(t) & 0 & \alpha_3 & 0 \\
\beta_{1} & \beta_{2} & 0 & \beta_{4} \\
0 & 0 & \gamma_{}(t) & 0 \\
0 & \nu_2 & 0 & \nu_4 \\
\end{array}
\right) \left(
\begin{array}{c}
X_{t-1}^{(1)} \\
X_{t-1}^{(2)} \\
X_{t-1}^{(3)} \\
X_{t-1}^{(4)} \\
\end{array}
\right) + \underline{\varepsilon}_{t}
= A(t)\underline{X}_{t-1}+\underline{\varepsilon}_{t},
\end{eqnarray*}
where $\{\underline{\varepsilon}_{t}\}$ are independent random variables (i.i.d)
with $\underline{\varepsilon}_{t}\sim N(0,I_{4})$, and $\alpha(t)$, $\gamma(t)$
are smoothly varying functions of $t$. The four 
time series are marginally nonstationary, in the sense that for each $1\leq
a\leq 4$, the time series $\{X_{t}^{(a)}\}_{}$ is second order nonstationary. 

The inverse operator and network corresponding
to $\{\underline{X}_{t}\}_{}$ is given below and is deduced from the
transition matrix $A(t)$
(the explicit connection between $D$ and $\{A(t)\}$ is given in Section \ref{sec:tvVAR}). Note that 
red and blue denote Toeplitz and non-Toeplitz matrix operators
respectively. 

\begin{minipage}{0.5\textwidth}
\begin{eqnarray*}
  D = \left(
\begin{array}{cccc}
{\color{blue}D_{1,1}} & {\color{red}D_{1,2}} & {\color{blue}D_{1,3}} & {\color{red}D_{1,4}} \\
{\color{red}D_{2,1}} & {\color{red}D_{2,2}} & 0 & {\color{red}D_{2,4}} \\
{\color{blue}D_{3,1}} & 0 & {\color{blue}D_{3,3}} & 0 \\
{\color{red}D_{4,1}} & {\color{red}D_{4,2}} & 0 & {\color{red}D_{4,4}} \\ 
\end{array}
\right)
\end{eqnarray*}
\end{minipage}
\begin{minipage}{0.5\textwidth}
\includegraphics[width = 0.4\textwidth]{plotspaper/corrected_graph.png}
\end{minipage}

\noindent \underline{Connecting the transition matrix to the network} The connections between the nodes
is because node $1$ is connected to node $3$ (if
$\alpha_3\alpha(t)\neq 0$ for some $t$), node $4$ (if $\beta_{1}\beta_{2}\neq 0$)
and node 2 (if $\beta_{1}\beta_{4}\neq 0$). By a similar argument, 
nodes $2$ and $4$ are connected (if $\beta_{1}\beta_{4}\neq 0$ or
$\nu_{2}\nu_{4}\neq 0$).

The nonstationarity of the multivariate time series is due to the
time-varying parameters $\alpha(t)$ and $\gamma(t)$. Specifically, the
parameter $\alpha(t)$ is the reason that node $1$ is nonstationary,
and by a similar argument the time-varying parameter $\gamma(t)$ is the reason
node $3$ is nonstationary. Since the coefficients on the second and fourth columns are not
time-varying, nodes $2$ and $4$ have
``inherited'' their nonstationarity from nodes
$1$ and $3$. Thus nodes $1$ and $3$ are conditionally stationary
whereas nodes $2$ and $4$ are conditionally stationary.
The connections between nodes $1$ to $2$ and $1$ to $4$ are
time-invariant because $\beta_{1}\beta_{2}$ and $\beta_{1}\beta_{4}$
are time-invariant respectively. 

\end{example}

\begin{remark}[Connection to GGM]
Let $X^{(a)} = (X_{t}^{(a)};t\in \mathbb{Z})$. It is clear that the density
of the infinite dimensional vector 
$(X^{(a)};1\leq a \leq p)$ is not well defined. However, we can
informally view the joint density (at least in the Gaussian case) as ``being proportional to''
\begin{eqnarray*}
\exp\left(-\frac{1}{2}\sum_{a=1}^{p}\langle X^{(a)}, D_{a,a}X^{(a)}  \rangle    -  
\frac{1}{2}\sum_{(a,b)\in E,a\neq b}\langle X^{(a)},D_{a,b}X^{(b)}
  \rangle\right)
\end{eqnarray*}
This is analogous to the representation of multivariate Gaussian
vector in terms of its inverse covariance. Using the
  above representation we conjecture that the above notions of
conditional correlation/stationarity/nonstationarity can be
generalized to time series which is
not necessarily continuous valued, for example binary valued time series.
\end{remark}

\subsection{NonStGM as a generalization of GGM}\label{sec:covD}

We start by defining partial covariances in the spirit of the definition
used in GGM but for infinite dimensional
random variables. This is defined by removing {\it two 
random variables} from the spanning set of $\mathcal{H}$
\begin{eqnarray}
\label{eq:partial1}
\rho_{t,\tau}^{(a,b)} = \cov\left[X_{t}^{(a)} - P_{\mathcal{H}  -
  (X_{t}^{(a)},X_{\tau}^{(b)})}(X_{t}^{(a)}), X_{\tau}^{(b)} - P_{\mathcal{H}  -
  (X_{t}^{(a)},X_{\tau}^{(b)})}(X_{\tau}^{(b)})\right].
\end{eqnarray}
Note that for the case $t=\tau$ and $a=b$ the above reduces to 
\begin{eqnarray}
\label{eq:partial1K}
\rho_{t,t}^{(a,a)} = \var\left[X_{t}^{(a)} - P_{\mathcal{H}  -
  X_{t}^{(a)}}(X_{t}^{(a)})\right] = \sigma_{a,t}^{2}.
\end{eqnarray}
In  the discussion below we refer to the infinite dimensional
  conditional covariance matrices
$\rho^{(a,b)} = (\rho_{t,\tau}^{(a,b)};t,\tau\in \mathbb{Z})$ and
$\rho^{(a,a)} = (\rho_{t,\tau}^{(a,a)};t,\tau\in \mathbb{Z})$.
In GGM the partial covariances are encoded in the precision
matrix. In a similar spirit, we show that
$\rho^{(a,b)}_{t,\tau}$ is encoded in the inverse covariance operator  $D$.

\begin{lemma}\label{lemma:D}
Suppose Assumption \ref{assum:lambda} holds. 
Let $D_{a,b}$  be defined as in (\ref{eq:Dab}). Then the entries of $D_{a,b}$ satisfy the identities
\begin{eqnarray}
\label{eq:corpartial}
\cor\left[X_{t}^{(a)} - P_{\mathcal{H}  -
  (X_{t}^{(a)},X_{\tau}^{(b)})}(X_{t}^{(a)}), X_{\tau}^{(b)} - P_{\mathcal{H}  -
  (X_{t}^{(a)},X_{\tau}^{(b)})}(X_{\tau}^{(b)})\right]
= -\frac{[D_{a,b}]_{t,\tau}}{\sqrt{[D_{a,a}]_{t,t} [D_{b,b}]_{\tau,\tau}}}
\end{eqnarray}
and 
\begin{eqnarray}
\label{eq:covpartial}
\var\left[
\left(
\begin{array}{c}
X_{t}^{(a)} - P_{\mathcal{H}  -
  (X_{t}^{(a)},X_{\tau}^{(b)})}(X_{t}^{(a)}) \\
X_{\tau}^{(b)} - P_{\mathcal{H}  -
  (X_{t}^{(a)},X_{\tau}^{(b)})}(X_{\tau}^{(b)}) \\
\end{array}
\right)
\right] = 
\left(
\begin{array}{cc}
[D_{a,a}]_{t,t} & [D_{a,b}]_{t,\tau} \\
{} [D_{b,a}]_{\tau,t} & [D_{b,b}]_{\tau,\tau} \\ 
\end{array}
\right)^{-1}.
\end{eqnarray}
\end{lemma}
\noindent {\bf PROOF} See Appendix \ref{sec:covDproofs}. \hfill $\Box$

\vspace{2mm}

An immediate consequence of Lemma \ref{lemma:D} is that the notions of conditional noncorrelation
  and conditional stationarity can be equivalently defined in terms of
  the properties of the partial covariances $\rho^{(a, b)}$.
  In particular, conditional noncorrelation between the two series $a$ and $b$ translates to zero $\rho^{(a,b)}$,
  while conditional stationarity of the pair $(a,b)$ translates to
  Toeplitz structures on $\rho^{(a,a)}$, $\rho^{(b,b)}$ and $\rho^{(a,b)}$. It is
  worth noting that the Toeplitz structure of $\rho^{(a,a)}$ (the  partial covariance of $a$)
captured in our framework is an important property, viz., the conditional
(non)stationarity of a node. A similar role
 on the diagonal entries of the precision or spectral precision matrices
($\Theta_{a,a}$ or $[\Sigma^{-1}(\omega)]_{a,a}$) is 
  absent in both the classical GGM and  StGM frameworks.

\begin{proposition}[NonStGM in terms of $\rho^{(a,b)}_{t,\tau}$]\label{lemma:networkpartial}
Suppose Assumption \ref{assum:lambda} holds. Let $\rho^{(a,b)}_{t,\tau}$ be defined as
in (\ref{eq:partial1}).  Then
\begin{itemize}
\item {\bf  Conditional Noncorrelation} 
$\rho_{t,\tau}^{(a,b)}=0$ for all $t$ and $\tau$ (i.e. $\rho^{(a,b)}=0$) iff
$D_{a,b} =0$
\item {\bf Conditionally Stationary Node} $D_{a,a}$ is Toeplitz
 iff for all $t$ and $\tau$
\begin{eqnarray*}
\rho_{t,\tau}^{(a,a)} = \rho_{0,t-\tau}^{(a,a)},
\end{eqnarray*}
i.e. $\rho^{(a,a)}$ is Toeplitz.
\item  {\bf Conditionally Stationary Pair} $D_{a,a}$, $D_{b,b}$ and
  $D_{a,b}$ are Toeplitz iff for all $t$ and $\tau$
\begin{eqnarray*}
\var\left[
\left(
\begin{array}{c}
X_{t}^{(a)} - P_{\mathcal{H}  -
  (X_{t}^{(a)},X_{\tau}^{(b)})}(X_{t}^{(a)}) \\
X_{\tau}^{(b)} - P_{\mathcal{H}  -
  (X_{t}^{(a)},X_{\tau}^{(b)})}(X_{\tau}^{(b)}) \\
\end{array}
\right)
\right] = 
\left(
\begin{array}{cc}
\rho_{0,t-\tau}^{(a,a)} & \rho_{0,t-\tau}^{(a,b)} \\
\rho_{0,t-\tau}^{(a,b)} & \rho_{0,t-\tau}^{(b,b)} \\
\end{array}
\right).
\end{eqnarray*}
i.e. $\rho^{(a,a)}, \rho^{(b,b)}$ and $\rho^{(a,b)}$ are Toeplitz.
\end{itemize}
\end{proposition}
\noindent{\bf PROOF} See Appendix \ref{sec:covDproofs}. \hfill $\Box$

\subsection{NonStGM as a generalization of StGM}\label{sec:covtimeD}

Now we define the time series partial covariance analogous to that 
used in StGM. We recall that the classical time series definition of partial
covariance in a multivariate time series evaluates the covariance between two random variables
$X^{(a)}_{t}$ and $X^{(b)}_\tau$, after conditioning on all random
variables in the $(p-2)$ component series $V\backslash \{a,b\}$. In
other words, we exclude the {\it entire} time series $a$ and $b$ from
the conditioning set. 

Formally, for any $\mathcal{S} \subseteq V$, we define the residual of
$X_{t}^{(a)}$ after projecting on 
$\overline{\textrm{sp}}
  (X_{s}^{(c)};s\in \mathbb{Z},c\notin \mathcal{S}) = \mathcal{H} -
  (X^{(c)};c\in \mathcal{S})$  as 
\begin{eqnarray*}
X_{t}^{(a)|\shortminus\mathcal{S}} &:=& X_{t}^{(a)} - P_{\mathcal{H} -
  (X^{(c)};c\in \mathcal{S})}(X_{t}^{(a)})  \textrm{ for }t\in \mathbb{Z}. 
\end{eqnarray*}
In the definitions below we focus on the two sets
$\mathcal{S}=\{a,b\}$ and $\mathcal{S}=\{a\}$. We mention that the set $\mathcal{S}=\{a\}$ is not
considered in StGMM but plays an important role in NonStGM.
Using the above, we define the edge partial covariance
\begin{eqnarray}
\label{eq:partialTS1}
\left(
\begin{array}{cc}
\rho_{t,\tau}^{(a,a)|\shortminus\{a,b\}} & \rho_{t,\tau}^{(a,b)|\shortminus\{a,b\}} \\
\rho_{t,\tau}^{(b,a)|\shortminus\{a,b\}} &\rho_{t,\tau}^{(b,b)|\shortminus\{a,b\}}\\
\end{array}
\right) :=
\cov\left[ \left(
\begin{array}{c}
X_{t}^{(a)|\shortminus\{a,b\}} \\
X_{t}^{(b)|\shortminus\{a,b\}} \\
\end{array}
 \right),
\left(
\begin{array}{c}
X_{\tau}^{(a)|\shortminus\{a,b\}} \\
X_{\tau}^{(b)|\shortminus\{a,b\}} \\
\end{array}
 \right)
\right]
\end{eqnarray}
and node partial covariance
\begin{eqnarray}
\label{eq:partialTS2}
\rho_{t,\tau}^{(a,a)|\shortminus\{a\}} = \cov[X_{t}^{(a)|\shortminus\{a\}}, X_{\tau}^{(a)|\shortminus\{a\}}].
\end{eqnarray}

We will show that the  partial covariance in (\ref{eq:partialTS1}) and (\ref{eq:partialTS2}) are closely
related to the partial covariance in (\ref{eq:partial1}).
In Lemma \ref{lemma:D} we have shown that the partial correlations $\rho_{t,\tau}^{(a,b)}$
define the entries of the operator $D$. We now connect the time series
definition of a partial covariance to the operator $D=(D_{a,b};a,b\in
\{1,\ldots,p\})$. Before we present the equivalent definitions of our nonstationary networks in terms of the
  time series partial covariances
  $\rho_{t,\tau}^{(a,b)|\shortminus\mathcal{S}}$, we  show that
  $\rho_{t,\tau}^{(a,b)|\shortminus\mathcal{S}}$ can be expressed 
in terms of the inverse covariance operator $D$.

\begin{theorem}\label{lemma:GGMstGGM}
Suppose Assumption \ref{assum:lambda} holds. Let 
$\rho_{t,\tau}^{(a,a)|\shortminus\{a,b\}}$,$\rho_{t,\tau}^{(a,b)|\shortminus\{a,b\}}$
and $\rho_{t,\tau}^{(a,a)|\shortminus\{a\}}$ be defined as in (\ref{eq:partialTS1})
and (\ref{eq:partialTS2}) respectively.
Then 
\begin{itemize}
\item[(i)] $\rho_{t,\tau}^{(a,a)|\shortminus\{a\}}= [D_{a,a}^{-1}]_{t,\tau}$
\item[(ii)] If $a\neq b$, then
\begin{eqnarray}
\label{eq:Dabinv}
\var\left[X_{t}^{(c)|\shortminus\{a,b\}};t\in \mathbb{Z}, c\in \{a,b\}\right]
 &=& 
\left(
\begin{array}{cc}
D_{a,a} & D_{a,b} \\
D_{b,a} & D_{b,b} \\
\end{array}
\right)^{-1} 
\end{eqnarray}
with
\begin{eqnarray*}
\rho_{t,\tau}^{(a,a)|\shortminus\{a,b\}}
&=&  [-(D_{a,a} - D_{a,b}D_{b,b}^{-1}D_{b,a})^{-1}D_{a,b}D_{b,b}^{-1}
                                             ]_{t,\tau} \\
\rho_{t,\tau}^{(a,b)|\shortminus\{a,b\}} &=&[(D_{a,a} -
                                             D_{a,b}D_{b,b}^{-1}D_{b,a})^{-1}]_{t,\tau} \\
\rho_{t,\tau}^{(b,b)|\shortminus\{a,b\}}&=&  [(D_{b,b} -
                                             D_{b,a}D_{a,a}^{-1}D_{a,b})^{-1}]_{t,\tau}.
\end{eqnarray*}
\end{itemize}
\end{theorem}
\noindent {\bf PROOF} See Appendix \ref{sec:covtimeDproofs}. \hfill $\Box$

\vspace{2mm}

\noindent 
A careful examination of the expressions for the GGM
  covariance   $\rho_{t,\tau}^{(a,b)}$ given in Lemma \ref{lemma:D}
  with the StGM covariance given in $\rho_{t,\tau}^{(a,b)|-\{a,b\}}$
  shows they are very different quantities. Therefore, it is suprising that despite
  these stark differences they preserve the same structures. More precisely, 
in Proposition  \ref{lemma:networkpartial} we showed that 
Definition \ref{def:network} had a clear interpretation in terms of 
$\rho_{t,\tau}^{(a,b)}$. We show below that  the
network definition given in Definition \ref{def:network} can be
interpreted in terms of the 
conditional dependence (or residuals) of the time series. The fact
that two very different conditional covariance definitions lead to the
same conditional graph is due to the property that infinite dimensional Toeplitz
operators remain Toeplitz even after  inversion and multiplication
with other Toeplitz operators. 

\begin{theorem}\label{lemma:conditional}[NonStGM in terms of $\rho^{(a,b) | \shortminus \{a,b \}}$]
Suppose Assumption \ref{assum:lambda} holds.
Let $\rho_{t,\tau}^{(a,a)|\shortminus\{a,b\}}$, $\rho_{t,\tau}^{(a,b)|\shortminus\{a,b\}}$
and $\rho_{t,\tau}^{(a,a)|\shortminus\{a\}}$ be defined as in (\ref{eq:partialTS1})
and (\ref{eq:partialTS2}) respectively.
Then
 \begin{itemize}
\item[(i)]   {\bf  Conditional noncorrelation} 
$D_{a,b}=0$ iff 
$\rho_{t,\tau}^{(a,a)|\shortminus\{a\}}=0$ for all $t$ and $\tau$.
\item[(ii)] {\bf Conditionally stationary node}  $D_{a,a}$ is a Toeplitz
  operator iff for all $t$ and $\tau$, $\rho_{t,\tau}^{(a,a)|\shortminus\{a\}}= \rho_{0,t-\tau}^{(a,a)|\shortminus\{a\}}$.
\item[(iii)] {\bf Conditionally stationary pair}  $D_{a,a}$,
  $D_{b,b}$ and $D_{a,b}$ are Toeplitz iff  for all $t$ and $\tau$,
$\rho_{t,\tau}^{(a,a)|\shortminus\{a,b\}}= \rho_{0,t-\tau}^{(a,a)|\shortminus\{a,b\}}$,
$\rho_{t,\tau}^{(b,b)|\shortminus\{a,b\}}=
\rho_{0,t-\tau}^{(b,b)|\shortminus\{a,b\}}$ and $\rho_{t,\tau}^{(a,b)|\shortminus\{a,b\}}= \rho_{0,t-\tau}^{(a,a)|\shortminus\{a,b\}}$.
\end{itemize}
\end{theorem}
{\bf PROOF} See Appendix \ref{sec:covtimeDproofs}. \hfill $\Box$

\vspace{2mm}

\noindent We show in the following result that
the \emph{time series partial covariances}
  can be used to define conditional stationarity of a subgraph 
containing three or more nodes.

\begin{corollary}[Conditionally stationary subgraph]\label{cor:nodes}
Let $\mathcal{S}= \{\alpha_{1},\ldots,\alpha_{r}\}$ be a subset 
of $\{1,\ldots,p\}$ and $\mathcal{S}^{\prime}$
denote the complement of $\mathcal{S}$.
Suppose for all $a,b\in \mathcal{S}$, $D_{a,b}$ are Toeplitz
(including the case $a=b$). Then 
$\{X_{t}^{(a)};t\in \mathbb{Z},a\in \mathcal{S}\}$ is a
conditionally stationary subgraph where 
\begin{eqnarray*}
\var\left[X_{t}^{(a)} -  
P_{\mathcal{H} - (X^{c};c\in
  \mathcal{S}^{\prime})}(X_{t}^{(a)});t\in \mathbb{Z},a\in
  \mathcal{S}\right]  =  P
\end{eqnarray*}
with
\begin{eqnarray*}
P^{-1} = 
\left(
\begin{array}{cccc}
D_{\alpha_{1},\alpha_{1}} & D_{\alpha_{1},\alpha_{2}} & \ldots &
                                                                 D_{\alpha_{1},\alpha_{r}} 
  \\
 D_{\alpha_{2},\alpha_{1}} & D_{\alpha_{2},\alpha_{2}} & \ldots &
                                                                 D_{\alpha_{2},\alpha_{r}} \\
\vdots & \vdots & \ddots & \vdots \\
D_{\alpha_{r},\alpha_{1}} & D_{\alpha_{r},\alpha_{2}} & \ldots &
                                                                 D_{\alpha_{r},\alpha_{r}} \\
\end{array}
\right).
\end{eqnarray*}
\end{corollary}

\noindent {\bf PROOF} See Appendix \ref{sec:covtimeDproofs}. \hfill $\Box$


%% file: section3_Revision.tex
\section{Sparse characterisations within the Fourier domain}\label{sec:fourier}

For general nonstationary processes 
it is infeasible to estimate the operator $D$ and learn its network
within the time domain. The problem is akin to StGM, where it is
difficult to learn the graph structure in 
the time domain by studying all the autocovariance
matrices. Estimation is typically carried out in the
Fourier domain by detecting conditional independence from the zeros of 
$\Sigma^{-1}(\omega)$. Following the same route, we will switch to the Fourier domain and construct a 
quantity that can be used to ``detect zeros and non-zeros''.  In
addition, within the Fourier domain we will define 
  meaningful notions of  weights/strengths of conditionally stationary
nodes and pairs that are analogous to well-known partial spectral coherence measures used in StGM.

 \textit{Notation} We first summarize some of the notation we will use
in this section. 
We define the function space of square integrable functions
$L_{2}[0,2\pi)$ as all complex functions where $g\in L_{2}[0,2\pi)$ if
$\int_{0}^{2\pi}|g_{}(\omega)|^{2}d\omega<\infty$. 
We define the function space of all square summable vector complex functions $L_{2}[0,2\pi)^{p}$, where
$\underline{g}(\omega)^{\prime} =
(g_{1}(\omega),\ldots,g_{p}(\omega))\in L_{2}[0,2\pi)^{p}$ if for all $1\leq j \leq
p$ $g_{j}\in L_{2}[0,2\pi)$. For all
$\underline{g},\underline{h}\in L_{2}[0,2\pi)^{p}$ we define the
inner-product $\langle \underline{g}, \underline{h}  \rangle =
\sum_{j=1}^{p} \langle g_{j},h_{j}\rangle$, where 
$\langle g_{j},h_{j}\rangle =
\int_{0}^{2\pi}g_{j}(\omega)h_{j}(\omega)^{*}d\omega$. Note that 
$L_{2}[0,2\pi)^{p}$ is a Hilbert space. We use $\delta_{\omega,\lambda}$ to denote
the Dirac delta function and set $i=\sqrt{-1}$.

\subsection{Transformation to the Fourier domain}\label{sec:fourierback}

In this section we summarize results which are pivotal to the
development in the subsequent sections. This section can be skipped on
first reading. 

To connect the time and Fourier domain we define a transformation
between the sequence and function space. We define the functions 
$F:L_{2}[0,2\pi)\rightarrow \ell_{2}$ and $F^{*}:
\ell_{2}\rightarrow L_{2}[0,2\pi)$ 
\begin{eqnarray}
\label{eq:FourierDef}
[F(g)]_{j}  =
  \frac{1}{2\pi}\int_{0}^{2\pi}\underline{g}(\lambda)\exp(i
  j\lambda)d\lambda \quad\textrm{and}\quad
F^{*}(v)(\omega) = \sum_{j\in \mathbb{Z}}v_{j}\exp(-ij\omega).
\end{eqnarray}
It is well known that $F$ and $F^{*}$ are isomorphisms between
$\ell_{2}$ and $L_{2}[0,2\pi)$ (see, for example, \cite{b:bro-dav-06}, Section 2.9). For $d>1$ the transformations 
$F(\underline{g}) = (F(g_{1}),\ldots, F(g_{d}))$ and 
$F^{*}v =
(F^{*}v^{(1)},\ldots,F^{*}v^{(d)})$ where $v =
(v^{(1)},\ldots,v^{(d)})$ are isomorphisms between $\ell_{2,d}$ and
$L_{2}[0,2\pi)^{d}$. Often we use that $d=p$.
These two isomorphims will provide a link between the infinite dimensional
matrix operators $D$ defined in the time domain to an equivalent
operator in the Fourier domain. 

Let 
$A = (A_{a,b};a,b\in \{1,\ldots,d\})$, if $A:
\ell_{2,d}\rightarrow \ell_{2,d}$ is a bounded operator, then standard
results show that $F^{*}AF: L_{2}[0,2\pi)^{d}\rightarrow
L_{2}[0,2\pi)^{d}$ is a bounded operator (see \cite{b:con-90}, Chapter
II). 
$F^{*}AF$ is an integral operator, such that for all $g\in L_{2}[0,2\pi)^{d}$
\begin{eqnarray}
\label{eq:integraloperator}
F^{*}AF(g)[\omega] =  \frac{1}{2\pi}\int_{0}^{2\pi}{\bf A}(\omega,\lambda)\underline{g}(\lambda)d\lambda,
\end{eqnarray}
and  ${\bf A}$ is the $d\times d$-dimensional matrix integral kernel
where 
\begin{eqnarray*}
{\bf A}(\omega,\lambda) &=&  \left( \sum_{t\in \mathbb{Z}}\sum_{\tau \in
\mathbb{Z}}[A_{a,b}]_{t,\tau}\exp(it\omega - i\tau \lambda);a,b\in \{1,\ldots,d\} \right).
\end{eqnarray*} 
To understand how $A$ and ${\bf A}(\omega,\lambda)$ are related we
focus on the case $d=1$ and note that the $(t,\tau)$ entry of the
infinite dimensional matrix $A$ is 
\begin{eqnarray*}
A_{t,\tau} = \frac{1}{(2\pi)^{2}}\int_{0}^{2\pi}\int_{0}^{2\pi}
  {\bf A}(\omega,\lambda)\exp(-it\omega+i\tau\lambda)d\omega d\lambda
  \quad \textrm{for all }t,\tau\in \mathbb{Z}.
\end{eqnarray*}

\begin{remark}[Connection with covariances and stationary time series]
We note if $C$ were a covariance operator of a univariate time series
$\{X_{t}\}$ with integral kernel $G$ then 
\begin{eqnarray}
\label{eq:dual}
\cov[X_{t},X_{\tau}] =  C_{t,\tau} = \frac{1}{(2\pi)^{2}}\int_{0}^{2\pi}\int_{0}^{2\pi}
  G(\omega,\lambda)\exp(-it\omega+i\tau\lambda)d\omega d\lambda,
\end{eqnarray}
where $G(\omega,\lambda)$ is the Lo\`{e}ve dual frequency spectrum. 
The Lo\`{e}ve dual
frequency spectrum is used to describe nonstationary features in a
time series and has been extensively studied in
\cite{p:gla-63},
\cite{p:lun-95}, \cite{p:lii-02}, \cite{p:jen-col-07}, \cite{p:hin-olh-10}, \cite{p:olh-11}, \cite{p:olh-omb-13},
\cite{p:gor-19}, \cite{p:ast-19}. 

If $\{X_{t}\}$ were a second order stationary time series, then (\ref{eq:dual})
reduces to Bochner's Theorem
\begin{eqnarray*}
\cov[X_{t},X_{\tau}] =  C_{0,t-\tau} = \frac{1}{(2\pi)}\int_{0}^{2\pi}
  f(\omega)\exp(-i(t-\tau)\omega)d\omega.
\end{eqnarray*}
The relationship between the  spectral density function
$f(\omega)$ and the Lo\`{e}ve dual frequency spectrum
$G(\omega,\lambda)$ 
is made apparent in Lemma \ref{lemma:diagonal} below.
\end{remark}

${\bf A}(\omega,\lambda)$ is a formal representation and typically
it will not be a well defined function over $[0,2\pi)^{2}$, as it is likely to have
singularities. 
Despite this, it has a very specific sparsity
structure when the operator $A$ is Toeplitz. For the
identification of nodes and edges in the nonstationary networks it is
the location of zeros in
${\bf A}(\omega,\lambda)$ that we will exploit. 
This will become apparent in the following lemma due to
\cite{p:toe-11} (we state the result for the case $d=1$).

\begin{lemma}\label{lemma:diagonal}
Suppose $A$ is an infinite dimensional bounded matrix operator
$A:\ell_{2}\rightarrow \ell_{2}$. The matrix operator $A$ is Toeplitz 
iff the integral kernel associated with $F^{*}AF$ has the form
\begin{eqnarray*}
A(\omega,\lambda) = \delta_{\omega,\lambda}A(\omega)
\end{eqnarray*}
where $A(\omega) \in L_{2}[0,2\pi)$ and $\delta_{\omega,\lambda}$ is 
the Dirac delta function.
\end{lemma}
\noindent {\bf PROOF} See Appendix \ref{sec:fourierproofgen} for details. \hfill $\Box$ 

\vspace{1mm}
\noindent The crucial observation in the above lemma is that 
$A(\omega,\lambda)=0$ for $\lambda \neq \omega$ iff $A$ is a Toeplitz matrix.
Below we  generalize the above to the case that $A$ (and its
inverse) is a block Toeplitz matrix operator.

\begin{lemma}\label{lemma:inversegeneral}
Suppose that $A$ is an infinite dimensional, symmetric, block 
matrix operator $A:\ell_{2,d}\rightarrow \ell_{2,d}$ where $0 < \inf_{\|v\|_{2}=1}\langle
v,Av\rangle \leq \sup_{\|v\|_{2}=1}\langle
v,Av\rangle <\infty$ with $A = (A_{a,b};a,b\in \{1,\ldots,d\})$ and
$A_{a,b}$ is Toeplitz. Then the integral kernel associated with
$F^{*}AF$ is ${\bf A}(\omega,\lambda) = {\bf A}(\omega)\delta_{\omega,\lambda}$ 
where ${\bf A}(\omega)$ is a $d\times d$ matrix with entries
$[{\bf A}(\omega)]_{a,b} = \sum_{r\in
  \mathbb{Z}}[A_{a,b}]_{0,r}\exp(ir\omega)$. 
Further the integral kernel associated with
$F^{*}A^{-1}F$ is ${\bf A}(\omega)^{-1}\delta_{\omega,\lambda}$. 
\end{lemma}
\noindent {\bf PROOF} In Appendix \ref{sec:fourierproofgen}. \hfill $\Box$

\vspace{2mm}

\noindent From now on we say that the kernel ${\bf A}(\omega,\lambda)$
is diagonal if it can be represented as $\delta_{\omega,\lambda}{\bf A}(\omega)$.

We use the operators $F:L_{2}[0,2\pi)^{p}\rightarrow \ell_{2,p}$ and
$F^{*}:\ell_{2,p}\rightarrow L_{2}[0,2\pi)^{p}$ to
recast the covariance and inverse covariance operators of a multivariate time series within the
Fourier domain.  We recall that $C$ is the covariance operator of the
time series $\{\underline{X}_{t}\}$ and by using (\ref{eq:integraloperator})
$F^{*}CF$ is an integral operator with matrix kernel ${\bf
  C}(\omega,\lambda) = (C_{a,b}(\omega,\lambda);a,b\in \{1,\ldots,p\})$ where 
$C_{a,b}(\omega,\lambda) = \sum_{t\in \mathbb{Z}}\sum_{\tau \in
\mathbb{Z}}[C_{a,b}]_{t,\tau}\exp(it\omega - i\tau \lambda)$.

In the case that $\{\underline{X}_{t}\}_{}$ is second order
stationary, then
$C_{a,b}(\omega,\lambda)=[\Sigma(\omega)]_{a,b}\delta_{\omega,\lambda}$
where $\Sigma(\cdot)$ is the spectral density matrix of
$\{\underline{X}_{t}\}_{}$. However, if 
$\{\underline{X}_{t}\}_{}$ is second order nonstationary, then by 
Lemma \ref{lemma:diagonal} at least one of the kernels
$C_{a,b}(\omega,\lambda)$ will be non-diagonal. The dichotomy that the
mass of ${\bf C}(\omega,\lambda)$ lies on
the diagonal $\omega = \lambda$ if and only if 
the  underlying process is multivariate second order stationary is
used in  \citep{p:eph-01, p:dwi-sub-11, p:jen-sub-15} to test for
second order  stationarity. 


\subsection{The nonstationary inverse covariance in the Fourier domain}\label{sec:precisionfourier}

The covariance operator $C$ and corresponding integral kernel ${\bf C}(\omega,\lambda)$ does not distinguish between
\textit{direct} and \textit{indirect} nonstationary relationships. 
We have shown in Section \ref{sec:framework} that conditional relationships are encoded
in the inverse covariance $D$. Therefore in this section we study the
properties of the integral kernel corresponding to $F^{*}DF$. 
Under Assumption \ref{assum:lambda}, $D = C^{-1}$ is a bounded operator,
thus $F^{*}DF$ is a bounded operator defined by the matrix kernel
 ${\bf K}(\omega,\lambda) = (K_{a,b}(\omega,\lambda);a,b\in \{1,\ldots,p\})$ 
where
\begin{eqnarray}
\label{eq:integralKnonstatab}
K_{a,b}(\omega,\lambda) = \sum_{t\in \mathbb{Z}}\sum_{\tau \in
\mathbb{Z}}[D_{a,b}]_{t,\tau}\exp(it\omega - i\tau \lambda) = 
\sum_{t\in \mathbb{Z}}\Gamma_{t}^{(a,b)}(\lambda)\exp(it(\omega-\lambda))
\end{eqnarray}
and
\begin{eqnarray}
\label{eq:spectrumpartial}
\Gamma_{t}^{(a,b)}(\lambda) = \sum_{r\in \mathbb{Z}}[D_{a,b}]_{t,t+r}\exp(ir\lambda).
\end{eqnarray}
Note that under Assumption \ref{assum:lambda}
  $\|D\|<\infty$, this implies for all $a,b\in \{1,\ldots,p\}$ that
the sequence $\{[D_{a,b}]_{t,t+r}\}_{r}\in
  \ell_{2}$, thus $\Gamma_{t}^{(a,b)}(\cdot)\in L_{2}[0,2\pi]$.
As far as we are aware,  neither $K_{a,b}(\omega,\lambda)$ nor
$\Gamma_{t}^{(a,b)}(\lambda)$ haven been studied previously. But 
$\Gamma_{t}^{(a,b)}(\lambda)$ can be viewed as the inverse covariance version of the
time-varying spectrum that is commonly used to analyze
nonstationary covariances (see \cite{p:pri-65}, \cite{p:mar-fla-85},
\cite{p:dah-96}, \cite{p:bir-18}). We observe that 
$D_{a,b}$ is Toeplitz if and only if $\Gamma_{t}^{(a,b)}(\lambda)$ does not depend
on $t$. 


In the following theorem we show that ${\bf K}(\omega,\lambda)$ defines a very clear sparsity
pattern depending on the  conditional properties of $\{\underline{X}_{t}\}$. This will allow
us to discriminate between different
types of edges in a network. In particular, zero matrices $D_{a,b}$
map to zero kernels and Toeplitz matrices $D_{a,b}$ map to diagonal kernels.


\begin{theorem}\label{lemma:spectral}
Suppose Assumption \ref{assum:lambda} holds. Then
\begin{itemize}
\item[(i)] {\bf Conditionally noncorrelated}
$\{X_{t}^{(a)},X_{t}^{(b)}\}_{t}$ are conditionally
  noncorrelated  iff \\ $K_{a,b}(\omega,\lambda)\equiv 0 $ for all $\omega,\lambda\in
[0,2\pi]$.
\item[(ii)] {\bf Conditionally stationary node}
$\{X_{t}^{(a)}\}_{t}$ is conditionally stationary iff the integral
kernel 
$K_{a,a}(\omega,\lambda)$ is diagonal.
 \item[(iii)] {\bf Conditionally time-invariant edge}
The edge $(a,b)$ is conditionally time-invariant iff the integral kernel
$K_{a,b}(\omega,\lambda)$ is diagonal.
\end{itemize}
\end{theorem}
{\bf PROOF} In Appendix \ref{appendix:preccoh}. \hfill $\Box$

\vspace{2mm}
\noindent These equivalences show that conditional noncorrelatedness
and stationarity relationships 
in the graphical model, as defined by the $D$ operator, are encoded
in the object $K(.,.)$. 
This provides the foundation for an alternate route to learning the graph structure in the frequency domain.

\begin{example}
We return to tvAR$(1)$ model described in Example
\ref{example:running}. In Figure \ref{fig:plotfourier} we give a schematic illustration of the matrix
$D$ in the frequency domain 
\begin{figure}
\begin{center}
\includegraphics[width = 0.8\textwidth]{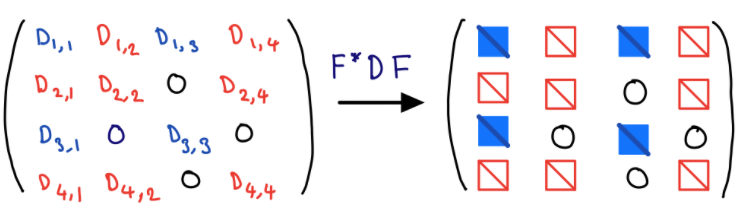}
\caption{Illustration of the mapping of matrix $D$ to the integral
  kernel corresponding to $F^{*}DF$. 
The diagonal red box indicates the mass of
$F^{*}D_{a.b}F$ lies only on the diagonal (it
corresponds to a Toeplitz matrix). The blue
  filled box indicates that the mass of
  $F^{*}D_{a.b}F$  lies both on the diagonal and elsewhere
  (it corresponds to a non-Toeplitz matrix). 
\label{fig:plotfourier}}
\end{center}
\end{figure}
\end{example}

\subsection{Partial spectrum for conditionally stationary
  time series}\label{sec:coherence}

So far we have considered the construction of an undirected, unweighted
  network which encodes the conditional uncorrelation and nonstationarity properties of
  time series components. In practice, we would be interested in assigning weights
  to network edges that represent the strength or magnitude of these conditional relationships.
  This will also be useful for learning the graph structure from finite samples.
  In GGM, partial correlation values are used to define edge weights,
  In StGM the partial
  spectral coherence (the frequency domain analogue of partial correlation) is used to
define suitable edge weights. 

We now define the notion of partial spectral coherence 
for conditionally stationary  time series. 
We start by interpreting $\Gamma^{(a,a)}_{t}(\omega)$ and
$\Gamma^{(a,b)}_{t}(\omega)$, defined in (\ref{eq:spectrumpartial}),
in the case that the node or edge is conditionally stationary. In the
following proposition we relate 
these quantities to the partial covariance $\rho_{t,\tau}^{(a,b)}$. 
Analogous to the definition of $\rho_{a,b}^{(t,\tau)}$ we define the
partial correlation
\begin{eqnarray}
\label{eq:partialVar} 
\phi_{t,\tau}^{(a,b)} &=& 
\cor[X_{t}^{(a)} -  P_{\mathcal{H}  -
  (X_{t}^{(a)},X_{\tau}^{(b)})}(X_{t}^{(a)}),X_{\tau}^{(b)} -  P_{\mathcal{H}  -
  (X_{t}^{(a)},X_{\tau}^{(b)})}(X_{\tau}^{(b)})].
\end{eqnarray}
Using the above we obtain an expression for $\Gamma_{t}^{(a,a)}(\omega)$
and $\Gamma_{t}^{(a,b)}(\omega)$ in the case that an edge or a node is
conditionally stationary.

\begin{theorem}\label{lemma:spectral2}
Suppose Assumption \ref{assum:lambda} 
holds. Let $\rho_{t,\tau}^{(a,b)}$ and $\phi_{t,\tau}^{(a,b)}$ be
defined as in (\ref{eq:partial1}) and (\ref{eq:partialVar}).
 \begin{itemize}
\item[(i)] If the node $a$ is conditionally stationary, then  
$\Gamma_{t}^{(a,a)}(\omega) = \Gamma^{(a,a)}(\omega)$ for all $t$, where 
\begin{eqnarray*}
\Gamma^{(a.a)}(\omega) = 
\sum_{r=-\infty}^{\infty}[D_{a,a}]_{(0,r)}\exp(ir\omega) = 
\frac{1}{\rho_{0,0}^{(a,a)}}\left[ 1-
  \sum_{r\in \mathbb{Z}\backslash \{0\}}\phi_{0,r}^{(a,a)}\exp(ir\omega)\right]
\end{eqnarray*}
 \item[(ii)] If $(a,b)$ is a conditionally stationary pair, 
   then expressions for $\Gamma_{}^{(a,a)}(\omega)$ and
$\Gamma_{}^{(b,b)}(\omega)$ are given in (i) and 
$\Gamma_{t}^{(a,b)}(\omega) = \Gamma^{(a,b)}(\omega)$ for all $t$, where 
\begin{eqnarray*}
\Gamma^{(a,b)}(\omega) = 
\sum_{r=-\infty}^{\infty}[D_{a,b}]_{(0,r)}\exp(ir\omega) =
-\frac{1}{(\rho_{0,0}^{(a,a)}\rho_{0,0}^{(b,b)})^{1/2}}
  \sum_{r\in \mathbb{Z}}\phi_{0,r}^{(a,b)}\exp(ir\omega).
\end{eqnarray*}
\end{itemize}
\end{theorem}
{\bf PROOF} In Appendix \ref{appendix:preccoh}. \hfill $\Box$

\vspace{2mm}

\noindent For StGM, the partial spectral
coherence is typically defined in terms of the Fourier transform of
the partial time series covariances (see \cite{b:pri-81}, Section 9.3,
and \cite{p:dah-00b}). We now show that an analogous result
holds in the case of conditional stationarity.

\begin{theorem}\label{lemma:spectral3}
Suppose Assumption \ref{assum:lambda} holds.
\begin{itemize}
\item[(i)] If the node $a$ is conditionally stationary, then  
\begin{eqnarray*}
\sum_{r\in \mathbb{Z}}\rho_{0,r}^{(a,a)|\shortminus\{a\}}\exp(ir\omega)
  = \Gamma^{(a,a)}(\omega)^{-1}.
\end{eqnarray*}
\item[(ii)] If $(a,b)$ is a conditionally stationary pair, then 
\begin{eqnarray*}
\sum_{r\in \mathbb{Z}}
\left(
\begin{array}{cc}
\rho^{(a,a)|\shortminus\{a,b\}}_{0,r} & \rho^{(a,b)|\shortminus\{a,b\}}_{0,r}\\
\rho^{(b,a)|\shortminus\{a,b\}}_{0,r} &  \rho^{(b,b)|\shortminus\{a,b\}}_{0,r} \\
\end{array}
\right)
\exp(ir\omega) 
&=& 
\left(
\begin{array}{cc}
\Gamma^{(a,a)}(\omega) & \Gamma^{(a,b)}(\omega) \\
\Gamma^{(a,b)}(\omega)^{*} & \Gamma^{(b,b)}(\omega) \\
\end{array}
\right)^{-1}.
\end{eqnarray*}
\end{itemize}
\end{theorem}
{\bf PROOF} In Appendix \ref{appendix:preccoh}. \hfill $\Box$

The above allows us to define the notion of spectral partial  coherence
in the case that underlying time series is nonstationary. We recall
that the spectral partial coherence between $\{X_{t}^{(a)}\}_{t}$
and $\{X_{t}^{(b)}\}_{t}$ for stationary time series is the
standardized spectral  
conditional covariance (see \cite{p:dah-00b}). Analogously, by using
Theorem \ref{lemma:spectral3}(ii) the 
spectral partial coherence between the conditionally
stationary pair $(a,b)$ is
\begin{eqnarray}
\label{eq:Rabomega}
R_{a,b}(\omega) = -\frac{\Gamma^{(a,b)}(\omega)}{\sqrt{\Gamma^{(a,a)}(\omega) \Gamma^{(b,b)}(\omega)}}.
\end{eqnarray}
In Appendix \ref{sec:stat} we show how this expression is related to
the spectral partial coherence for stationary time series.

\subsection{Connection to node-wise regression}\label{sec:noderegression}

In Lemma \ref{lemma:inverse} we connected the coefficients of $D$ to
the coefficients in a linear regression.
The regressors are in the spanning set of $\mathcal{H}-X_{t}^{(a)}$.
In contrast, in node-wise regression each node is regressed on all of the
\emph{other} nodes (the coefficients in this regression can also be
connected to the precision matrix). 
We now derive an analogous
result for multivariate time series. In particular, we regress the time series at node $a$ ($\{X_{t}^{(a)}\}_{t}$)
onto all the other time series (excluding node $a$ i.e. the spanning
set of $\mathcal{H} - (X^{(a)})$) and connect these
to the matrix $D$. 
These results can be used to encode conditions for a
conditionally stationary edge in terms of the regression
coefficients. Furthermore, they allow us to deduce the time
series at node $a$ conditioned on all the other nodes (if the time
series is Gaussian).

The best linear predictor of $X_{t}^{(a)}$ given the ``other'' time series
$\{X_{s}^{(b)};s\in \mathbb{Z}, b\neq a\}$ is
\begin{eqnarray}
\label{eq:expectation1}
P_{\mathcal{H} - (X^{(a)})}(X_{t}^{(a)}) = \sum_{b\neq
  a}\sum_{\tau\in \mathbb{Z}}\alpha_{(\tau,b)\shortarrow (t,a)}X_{\tau}^{(b)}.
\end{eqnarray}
We group the coefficients according to time
series and define the infinite dimensional matrix $B_{b\shortarrow a}$
with entries 
\begin{eqnarray}
\label{eq:Bac}
[B_{b\shortarrow a}]_{t,\tau} = \alpha_{(\tau,b)\shortarrow
  (t,a)} \textrm{ for all } t,\tau\in \mathbb{Z}. 
\end{eqnarray}
In the lemma below we
connect the coefficients in the infinite dimensional matrix  $B_{b\shortarrow a}$
to $D_{a,b}$
\begin{proposition}\label{lemma:expectation}
Suppose Assumption \ref{assum:lambda} holds. Let $(D_{a,b};1\leq
a,b\leq p)$ be defined as in
(\ref{eq:Dab}). Then for all $b\neq a$ we have 
\begin{eqnarray}
\label{eq:expectation2}
D_{a,b}  = -D_{a,a}B_{b\shortarrow a}.
\end{eqnarray} 
\end{proposition}
\noindent {\bf PROOF} See Appendix \ref{appendix:regression}. \hfill $\Box$

\vspace{2mm}

In the following theorem we rewrite the conditions for conditional
noncorrelation and conditional time-invariant 
  edge in terms of node-regression coefficients.
\begin{theorem}\label{lemma:regressiontoeplitz}
Suppose Assumption \ref{assum:lambda} holds. Let $B_{b\rightarrow a}$ be defined as
in (\ref{eq:Bac}). Then
\begin{itemize}
\item[(i)] $B_{b\shortarrow a} = 0$ iff $D_{a,b}=0$. 
\item[(ii)] If $D_{a,a}$ and  $B_{b\shortarrow a}$ are Toeplitz, then $D_{a,b} =
-D_{a,a}B_{b\shortarrow a}$ is Toeplitz. 
\end{itemize}
\end{theorem}
{\bf PROOF} See Appendix \ref{appendix:regression}. \hfill $\Box$

\vspace{2mm}
\noindent Below we show that the integral kernel associated with
$B_{b\shortarrow a}$ has a clear sparsity structure. 
\begin{corollary}\label{cor:regfourier}
Suppose Assumption \ref{assum:lambda}  holds. Let $B_{b\shortarrow a}$ be defined as in
(\ref{eq:expectation2}). Let $K_{b\shortarrow a}(\omega,\lambda)$
denote the integral kernel associated with $B_{b\shortarrow a}$. Then
\begin{itemize}
\item[(i)] $B_{b\shortarrow a}$ is a bounded operator.
\item[(ii)] {\bf Conditionally noncorrelated}
$\{X_{t}^{(a)},X_{t}^{(b)}\}_{t}$ are conditionally
  noncorrelated  iff \\ $K_{a\shortarrow b}(\omega,\lambda)\equiv 0 $.
\item[(iii)] {\bf Conditionally stationary pair} 
$\{X_{t}^{(a)},X_{t}^{(b)}\}_{t}$ are conditionally jointly 
  stationary iff the kernels $K_{a,a}(\omega,\lambda)$,
  $K_{b,b}(\omega,\lambda)$ and $K_{b\shortarrow a}(\omega,\lambda)$
  are diagonal.
\end{itemize}
\end{corollary}
\noindent {\bf PROOF}  In Appendix \ref{appendix:regression}. \hfill $\Box$

\vspace{2mm}
We use the results above to deduce the conditional distribution of
$X^{(a)}$ under the assumption that the time series $\{\underline{X}_{t}\}_{}$ is jointly
Gaussian. The conditional distribution of $X^{(a)}$ given
$\mathcal{H} - (X^{(a)})$ is Gaussian where 
\begin{eqnarray*}
X^{(a)}|\mathcal{H} - (X^{(a)})\sim N\left( \sum_{b=1,b\neq
  a}^{p}B_{b\shortarrow a}X^{(b)}, D_{aa}^{-1}\right)
\end{eqnarray*}
with $\Ex[X^{(a)}|\mathcal{H} - (X^{(a)})] = \sum_{b=1,b\neq
  a}^{p}B_{b\shortarrow a}X^{(b)}$ and $\var[X^{(a)}|\mathcal{H} - (X^{(a)})]=D_{aa}^{-1}$.
Some interesting simplications can be made if the nodes and
corresponding edges are
conditionally stationary and time-invariant. If $X^{(a)}$ has a conditionally stationary
node, then by Theorem \ref{lemma:conditional}(ii)
the conditional variance
will be stationary (Toeplitz). If, in addition, the conditionally stationary node $a$
is connected to the set of nodes
$\mathcal{S}_{a}$ and all the edge connections are
conditionally time-invariant then by Theorem \ref{lemma:regressiontoeplitz}
the coefficients in the conditional
expectation are shift invariant where
\begin{eqnarray*}
\Ex[X_{t}^{(a)}|\mathcal{H} - (X^{(a)})] = \sum_{b\in \mathcal{S}_{a}}\sum_{j\in \mathbb{Z}}\alpha_{j}^{(b\shortarrow a)}X_{t-j}^{(b)}.
\end{eqnarray*}
Therefore, if the node $a$ is conditionally stationary and all its
connecting edges are conditionally time-invariant then the conditional
distribution $X^{(a)}|(\mathcal{H}-(X^{(a)}))$ is stationary.

%% file: section4_Revision.tex
\section{Learning the network  from finite length time series}\label{sec:finitelocal}

The network structure of $\{\underline{X}_{t}\}_{t}$ is succinctly described in terms of 
${\bf K}(\omega,\lambda)$. However, for the purpose of estimation,
there are three problems. The first is that ${\bf K}(\omega,\lambda)$ is
a singular kernel  making direct estimation impossible. The second
is that for conditional nonstationary time series
the structure of
$[{\bf K}(\omega,\lambda)]_{a,b}$ is not well defined. 
Finally, in practice we only observe a finite length sample
$\{\underline{X}_{t}\}_{t=1}^{n}$.  Thus our object of interest
changes from ${\bf K}(\omega,\lambda)$ to its finite dimensional
counterpart (which we define below). For the purpose of network identification, we show that
the finite dimensional version of ${\bf K}(\omega,\lambda)$
 inherits the sparse properties of  
${\bf K}(\omega,\lambda)$. Moreover, in a useful twist, whereas ${\bf K}(\omega,\lambda)$  
is a singular kernel its finite dimensional counterpart is a well defined
matrix, making estimation possible.

\subsection{Finite dimensional approximation}\label{sec:finite}

To obtain the finite dimensional version of ${\bf K}(\omega,\lambda)$,
we recall that the Discrete Fourier transform (DFT) can be viewed as the
analogous version of the Fourier operator $F$ (defined in  (\ref{eq:FourierDef}))  in finite dimensions.
Let $F_{n}$ denote the $(np\times np)$-dimension DFT transformation
matrix. It comprises of $p^{2}$ identical $(n\times n)$-dimension DFT
matrices, which we denote as $\mathcal{F}_{n}$. Define the concatenated 
$np$-dimension vector 
${\bf X}_{n}^{\prime} =
((\underline{X}^{(1)})^{\prime},\ldots,(\underline{X}^{(p)})^{\prime})$, where
$\underline{X}^{(a)} = (X_{1}^{(a)},\ldots,X_{n}^{(a)})^{\prime}$ for
$a\in \{1,\ldots,p\}$. 
Then $F_{n}^{*}{\bf X}_{n}$ is a $np$-dimension vector where
$(F_{n}^{*}{\bf X}_{n})^{\prime} =
((\mathcal{F}_{n}^{*}\underline{X}^{(1)})^{\prime},\ldots, 
(\mathcal{F}_{n}^{*}\underline{X}^{(p)})^{\prime})$
with 
\begin{eqnarray}
\label{eq:Jdef}
J_{k}^{(a)} = [\mathcal{F}_{n}^{*}\underline{X}^{(a)}]_{k} = 
  \frac{1}{\sqrt{n}}\sum_{t=1}^{n}X_{t}^{(a)}\exp(it\omega_{k}) 
  \quad k = 1,\ldots,n \textrm{ and }\omega_{k} = \frac{2\pi k}{n}.
\end{eqnarray}
Let $\var[{\bf X}_{n}] = C_{n}$, then 
$\var[F_{n}^{*}{\bf X}_{n}] = F_{n}^{*}C_{n}F_{n}$. Our focus will be
on the $(np\times np)$-dimensional inverse matrix
\begin{eqnarray*}
{\bf K}_{n} = [\var[F_{n}^*{\bf X}_{n}]]^{-1}=
  [F_{n}^{*}C_{n}F_{n}]^{-1} = [F_{n}]^{-1}C_{n}^{-1}[F_{n}^{*}]^{-1} 
= F_{n}^{*}\widetilde{D}_{n}F_{n},
\end{eqnarray*}
where $\widetilde{D}_{n}  =
C_{n}^{-1}$ and the above follows from the identity $F_{n}^{-1}=F_{n}^{*}$. Let  ${\bf K}_{n} = ( [{\bf K}_{n}]_{a,b};a,b\in
\{1,\ldots,p\})$ where $[{\bf K}_{n}]_{a,b}$ denotes the $(n\times n)$-dimensional sub-matrix of
${\bf K}_{n}$ and $[{\bf K}_{n}(\omega_{k_1},\omega_{k_2})]_{a,b}$ denotes the
$(k_{1},k_{2})$th entry in the submatrix matrix $[{\bf
  K}_{n}]_{a,b}$. For future reference we define the $(p\times
p)$-dimensional matrix ${\bf K}_{n}(\omega_{k_1},\omega_{k_2}) = 
([{\bf K}_{n}(\omega_{k_1},\omega_{k_2})]_{a,b};1\leq a,b\leq p)$.
We show below that $[{\bf K}_{n}(\omega_{k_1},\omega_{k_2})]_{a,b}$
can be viewed as the finite dimensional version of $K_{a,b}(\omega,\lambda)$.

\vspace{1mm}
The covariance matrix $C_{n}=\var[{\bf X}_{n}]$ is a submatrix of the infinite dimensional $C$.
Unfortunately its inverse $\widetilde{D}_{n} = C_{n}^{-1}$ is not a submatrix
of $D$. As our aim is to show that the properties of the inverse
covariance map to those in finite dimensions we will show that under
suitable conditions $\widetilde{D}_{n}$ can be approximated by a
finite dimensional submatrix of $D$. To do this we  represent
$\widetilde{D}_{n}$ as $p\times p$
submatrices each of dimension $n\times n$ 
\begin{eqnarray}
\label{eq:defDtilde}
\widetilde{D}_{n} = \left(\widetilde{D}_{a,b;n};a,b\in \{1,\ldots,p\} \right).
\end{eqnarray}
Analogously, we define $p\times p$ submatrices of $D$ each of dimension
$n\times n$
\begin{eqnarray}
\label{eq:defD}
D_{n} = \left(D_{a,b;n};a,b\in \{1,\ldots,p\} \right)
\end{eqnarray}
where $D_{a,b;n} = \{[D_{a,b}]_{t,\tau};t,\tau \in \{1,\ldots,n\}\}$.
Below we show that under suitable conditions,  $\widetilde{D}_{n} =
C_{n}^{-1}$ can be approximated well by $D_{n}$. This result
requires the following conditions on the rate of decay of the inverse covariances ${\bf
  D}_{t,\tau}$ which is stronger than the conditions in Assumption
\ref{assum:lambda}. 

\begin{assumption}
\label{assum:invcovarianceK}
The inverse covariance ${\bf D}_{t,\tau}$ defined in (\ref{eq:Dttau})
satisfy the condition 
\\*
$\sup_{t}\sum_{j\neq  0}|j|^{K}\|{\bf D}_{t,t+j}\|_{\infty}<\infty$ (for some
$K\geq 3/2$). 
\end{assumption}

The conditions in Assumption \ref{assum:invcovarianceK} are analogous
to those used in the analysis of stationary
time series, where certain conditions on the rate of decay of the
autocovariances coefficients are
often used. \cite{p:kra-sub-22} obtain an equivalence between the
  rate of decay on ${\bf D}_{t,\tau}$ and ${\bf C}_{t,\tau}$. In
  particular, \cite{p:kra-sub-22}
Theorem 2.1, show that under Assumption \ref{assum:lambda} and if
for some $K>7/2$ and all $|r|\neq 0 $ we have that $\sup_{t}\|{\bf C}_{t,t+r}\|_{}<K|r|^{-K}$
(where $\|\cdot\|$ denotes the spectral norm), then 
$\sup_{t}\|{\bf  D}_{t,t+r}\|_{}<K((1+\log|r|)/|r|)^{K-1}$. Thus 
Assumption \ref{assum:invcovarianceK} holds.

In the lemma below we obtain a bound between the rows of 
$\widetilde{D}_{n}$ and $D_{n}$. 

\begin{theorem}\label{lemma:meyer3a}
Suppose Assumptions \ref{assum:lambda}  and
\ref{assum:invcovarianceK} hold. 
Let $\widetilde{D}_{n}$ and $D_{n}$ be defined as in
(\ref{eq:defDtilde}) and (\ref{eq:defD}). 
Then for all $1\leq t \leq n$ we have 
\begin{eqnarray*}
\sup_{1\leq a\leq p}\left\| [\widetilde{D}_{n}]_{(a-1)n+t,\cdot} -  [D_{n}]_{(a-1)n+t,\cdot}\right\|_{1} 
&=&  O\left( \frac{(np)^{1/2}}{\min(|n+1-t|,|t|)^{K}}\right),
\end{eqnarray*}
where $A_{(a-1)n+t,\cdot}$ denotes the $((a-1)n+t)$th row of the
matrix $A$, or, equivalently the $t$th row along the $a$th block of $A$.
\end{theorem}
{\bf PROOF} See Appendix \ref{sec:prooffinite}. \hfill $\Box$

\vspace{2mm}
\noindent 
The theorem above shows that the further  $t$ lies from the two end
boundaries of the sequence $\{1,2,\ldots,n\}$ the better the
approximation between $[\widetilde{D}_{n}]_{(a-1)n+t,\cdot}$ and
$[D_{n}]_{(a-1)n+t,\cdot}$. For example when $t=n/2$ (recall that $p$ is
fixed)
$\| [\widetilde{D}_{n}]_{(a-1)n+t,\cdot} -  [D_{n}]_{(a-1)n+t,\cdot}\|_{1} =O(1/n^{K-1/2})$.
Using  Theorem \ref{lemma:meyer3a} we replace $F_{n}^{*}\widetilde{D}_{n}F_{n}$ with
$F_{n}^{*}D_{n}F_{n}$ to obtain the following approximation.

\begin{proposition}\label{lemma:conditionalstat}
Suppose Assumptions \ref{assum:lambda}  and
\ref{assum:invcovarianceK} hold. Let $\Gamma_{t}^{(a,b)}(\omega)$ be
defined as in (\ref{eq:spectrumpartial}).
Then 
\begin{eqnarray}
\label{eq:Gammatbound}
[{\bf K}_{n}(\omega_{k_{1}},\omega_{k_{2}})]_{a,b} &=& 
  \frac{1}{n}\sum_{t=1}^{n}\Gamma_{t}^{(a,b)}(\omega_{k_{2}})\exp(-it(\omega_{k_{1}}-\omega_{k_{2}}))
                                                       +
                                                       O\left(\frac{1}{n}
                                                       \right) \nonumber\\
&=&  \left[\frac{1}{n}\sum_{t=1}^{n}\Gamma_{t}^{(b,a)}(\omega_{k_{1}})\exp(-it(\omega_{k_{2}}-\omega_{k_{1}}))\right]^{*}
                                                       +                                                       O\left(\frac{1}{n}
                                                       \right) 
\end{eqnarray}
Further, if 
$\{X_{t}^{(a)}\}_{t}$ and $\{X_{t}^{(b)}\}_{t}$ are conditionally
stationary, then
\begin{eqnarray}
\label{eq:Gammatboundstat}
[{\bf K}_{n}(\omega_{k_{1}},\omega_{k_{2}})]_{a,b} = 
\left\{
\begin{array}{cc}
\Gamma^{(a,b)}(\omega_k) + O(n^{-1})   &  k_{1}=k_{2}(=k) \\
 O(n^{-1})  & k_{1} \neq k_{2}
\end{array}
\right.
\end{eqnarray}
where $\Gamma^{(a,b)}(\omega) = 
\sum_{r=-\infty}^{\infty} [D_{a,b}]_{(0,r)}\exp(ir\omega_{})$.
\end{proposition}
\noindent {\bf PROOF} See Appendix \ref{sec:prooffinite}. \hfill $\Box$

\subsection{Locally stationary time series}\label{sec:LS}

We showed in Proposition \ref{lemma:conditionalstat} that in the case the
node or edge $(a,b)$ is conditional stationary or conditionally time-variant   $[{\bf
  K}_{n}(\omega_{k_1},\omega_{k_2})]_{a,b}$ has a well defined structure; the
diagonal dominates the off-diagonal terms (which are of order $O(n^{-1})$). However,
in the case of conditional nonstationary node/time-varying edge the precise structure of $[{\bf
  K}_{n}(\omega_{k_1},\omega_{k_2})]_{a,b}$ is not apparent, this makes
detection of conditional nonstationarity difficult. In this section we
impose some structure on the form of the nonstationarity. We will work
under the canopy of local stationarity. It formalizes the notion
that the ``nonstationarity'' in a time series evolves ``slowly'' through time. It is
arguably one of the most popular methods for describing nonstationary
behaviour and describes a wide class of nonstationarity; various
applications are discussed in \cite{p:pri-65}, \cite{p:dah-gir-98}, \cite{p:zhou-09}, \cite{p:nas-10},
\cite{p:kle-17}, \cite{p:dah-19}, \cite{p:sun-pou-18},
\cite{p:din-zho-20}, \cite{p:omb-21}, to name
but a few. We show below that  for
locally stationary time series $[{\bf
  K}_{n}(\omega_{k_1},\omega_{k_2})]_{a,b}$ has a distinct structure
that can be  detected. 

The locally stationary process were formally proposed in
\cite{p:dah-96}. In the locally stationary framework the asymptotics hinge on the
rescaling device $n$, which is linked to the sample size. It measures
how close the nonstationary time series is to an
auxillary (latent) process $\{\underline{X}_{t}(u)\}_{t}$ which for a fixed $u$ is stationary
over $t$. More precisely, a time series $\{\underline{X}_{t,n}\}_{t}$ is said to be
locally stationary if there exists a stationary time series
$\{\underline{X}_{t}(u)\}_{t}$ where 
\begin{eqnarray}
\label{eq:localstationary}
\|\underline{X}_{t,n} - \underline{X}_{t}(u) \|_2 = O_{p}\left(\frac{1}{n} + \left|\frac{t}{n} - u\right|\right).
\end{eqnarray}
Thus for every $t$, $\underline{X}_{t,n}=(X_{t,n}^{(1)},\ldots,X_{t,n}^{(p)})^{\prime}$ can  be closely approximated by an 
auxillary variable $\underline{X}_{t}(u)$ (where $u=t/n$); see \citep{p:dah-sub-06,
p:sub-06, p:dah-12,p:dah-19}. 
However, as the difference between $t/n$ and $u$ grows, the similarity
between $X_{t,n}$ and the auxillary stationary process $X_{t}(u)$ decreases.
This asymptotic device
allows one to obtain well defined limits for nonstationary time
series which otherwise would not be possible within classical real time
asymptotics.  Though the formulation in (\ref{eq:localstationary}) is
a useful start for analysing nonstationary time series, analogous to
\cite{p:dah-pol-06}, we require additional local stationarity
conditions on the moment structure. \cite{p:dah-00a} and \cite{p:dah-pol-06}
state the conditions in terms of bounds between
$\cov[\underline{X}_{t,n},\underline{X}_{\tau,n}]$ and $\cov[\underline{X}_{0}(u),
\underline{X}_{t-\tau}(u)]$. Below we state similar conditions in
terms of the inverse covariances ${\bf D}_{t,\tau}$ and its stationary
approximation counterpart. 

\begin{assumption}
\label{assum:LS}
There exists a sequence $\{\ell(j)\}_{j}$ such that $\sum_{j\in
  \mathbb{Z}}j\ell(j)^{-1}<\infty$ and matrix function ${\bf D}_{t-\tau}:\mathbb{R}^{}\rightarrow
\mathbb{R}^{p\times p}$ where 
\begin{eqnarray}
\label{eq:qqq}
{\bf D}_{t,\tau} = {\bf D}_{t-\tau}\left(\frac{t+\tau}{2n}\right) + O\left(\frac{1}{n\ell(t-\tau)}\right)\quad t,\tau \in \mathbb{Z}.
\end{eqnarray}
Further, the matrix function ${\bf D}_{j}(\cdot)$ is such that (i)
$\sup_{u}\sum_{j\in \mathbb{Z}}\|j{\bf D}_{j}(u)\|_{1}<\infty$, \\ 
(ii) $\sup_{u}|\frac{d[{\bf D}_{j}(u)]_{a,b}}{du}|\leq \ell(j)^{-1}$, (iii) for
all $u, v\in \mathbb{R}$ $\|{\bf
  D}_{j}(u) - {\bf D}_{j}(v)\|_{1}\leq |u-v|\ell(j)^{-1}$  and (iv)
$\sup_{u}|\frac{d[{\bf D}_{j}(u)]_{a,b}}{du}|\leq \ell(j)^{-1}$. 

Standard within the locally stationary paradigm  ${\bf
  D}_{t,\tau}$ should be indexed by $n$ (but to simplify
notation we have dropped the $n$).
\end{assumption}
 Theorem 3.3 in \cite{p:kra-sub-22}  shows that 
Assumption \ref{assum:LS} is fulfilled by a large class of locally stationary time series under certain
smoothness conditions on their covariance.  

The above assumptions require that the entry wise derivative of matrix functions ${\bf
  D}_{t-\tau}(\cdot)$ exists. This technical condition can be relaxed to include matrix functions
${\bf D}_{j}(\cdot)$ of bounded variation (which would allow for change
point models as a special case) similar to
\cite{p:dah-pol-06}. 


The above assumptions allow for two important
types of behaviour (i) conditionally stationary nodes and
time-invariant edges where
$[{\bf D}_{j}(u)]_{a,b} = [{\bf D}_{j}]_{(a,b)}$  and  (ii) conditional nonstationarity where the partial covariance
between $X_{t}^{(a)}$ and $X_{t+j}^{(b)}$ (for fixed lag $j$)
evolves ``nearly'' smoothly over $t$.


\subsection{Properties of ${\bf K}_{n}(\omega_{k_1},\omega_{k_2})$
  under local stationarity}\label{sec:KLS}

Typically, the second order analysis of locally stationary time series
is conducted through its time varying spectral density matrix.
This is the spectral density matrix corresponding to the locally stationary
approximation $\{\underline{X}_{t}(u)\}_{t}$, which we denote as
${\boldsymbol \Sigma}(u;\omega)$. 
The time-varying spectral density matrix corresponding to
$\{\underline{X}_{t,n}\}_{t}$ is $\{{\boldsymbol
  \Sigma}(t/n;\omega)\}_{t}$. In contrast, in this section
our focus will be on the inverse ${\boldsymbol \Gamma}(u;\omega) =
{\boldsymbol \Sigma}(u;\omega)^{-1}$,  where by 
 Lemma \ref{lemma:inversegeneral},  
 ${\boldsymbol \Gamma}(u;\omega)$
is the Fourier transform of ${\bf D}_{j}(u)$ over the lags $j$ i.e.
\begin{eqnarray}
\label{eq:localprecision}
{\boldsymbol \Gamma}(u;\omega) = \sum_{j\in \mathbb{Z}}{\bf D}_{j}(u)\exp(ij\omega).
\end{eqnarray}
We note that ${\boldsymbol \Gamma}(u;\omega) =
(\Gamma^{(a,b)}(u;\omega);1\leq a,b\leq p)$. We use 
Assumption \ref{assum:LS} to relate $\Gamma^{(a,b)}(u;\omega)$ to 
$\Gamma^{(a,b)}_{t}(\omega)$ (defined in
(\ref{eq:spectrumpartial})). In particular, 
$\Gamma^{(a,b)}(t/n;\omega)$ is an
approximation of $\Gamma^{(a,b)}_{t}(\omega)$ and 
\begin{eqnarray}
\label{eq:GLS}
\left|\Gamma^{(a,b)}_{t}(\omega)  -  \Gamma^{(a,b)}(u;\omega)\right| \leq
  C\left(\left|\frac{t}{n}- u\right| + \frac{1}{n} \right).
\end{eqnarray}
Thus the time-varying spectral  \emph{precision} matrix corresponding to
$\{\underline{X}_{t,n}\}_{t}$ is $\{{\boldsymbol
  \Gamma}(t/n;\omega)\}_{t}$.

Our aim is to relate ${\boldsymbol \Gamma}(u;\omega)$ to ${\bf
  K}_{n}(\omega_{k_1},\omega_{k_{2}})$. First we notice that
${\boldsymbol \Gamma}(u;\omega)$ is ``local'' in the sense that it is a time
local approximation to the precision spectral density at time point $t
= \lfloor un \rfloor $. On the other
hand, ${\bf  K}_{n}(\omega_{k_1},\omega_{k_{2}})$ is ``global'' in the
sense that it is based on the entire observed time series. However, we show
below that ${\bf  K}_{n}(\omega_{k_1},\omega_{k_{2}})$ is connected to 
${\boldsymbol \Gamma}(u;\omega)$, as it measures how 
${\boldsymbol \Gamma}(u;\omega)$ evolves over time. These
insights allow us to deduce the network structure from ${\bf  K}_{n}(\omega_{k_1},\omega_{k_{2}})$.

In the following lemma we show that the entries of the matrix 
${\bf K}_{n}(\omega_{k_1},\omega_{k_{2}})$ can be approximated by the
Fourier coefficients of $\Gamma^{(a,b)}(\cdot;\omega)$, where 
\begin{eqnarray}
\label{eq:Hrab}
K_{r}^{(a,b)}(\omega) &=&\int_{0}^{1}\exp(-2\pi
    iru)\Gamma^{(a,b)}(u;\omega)du.
\end{eqnarray}
The Fourier coefficients $K_{r}^{(a,b)}(\omega)$ fully determine the
function $\Gamma^{(a,b)}(u;\omega)$. In particular (i) if all the Fourier
coefficients are zero then $\Gamma^{(a,b)}(u;\omega)=0$ (ii) if all
the Fourier coefficients are zero except 
$r=0$, then $\Gamma^{(a,b)}(u;\omega)$ does not depend on $u$. Using
this, it is clear the coefficients $K_{r}^{(a,b)}(\omega)$ hold
information on the network. We summarize these properties in the following proposition.

\begin{proposition}\label{lemma:NetworkK} 
Suppose Assumptions \ref{assum:lambda}, \ref{assum:invcovarianceK} and \ref{assum:LS}
hold. Let $K_{r}^{(a,b)}(\cdot)$ be defined as in (\ref{eq:Hrab}). Then
\begin{itemize}
\item[(i)] $\{X_{t,n}^{(a)}\}_{t=1}^{n}$ and
  $\{X_{t,n}^{(b)}\}_{t=1}^{n}$ is a (asymptotically) conditionally
  noncorrelated edge iff  $K_{r}^{(a,b)}(\omega) \equiv 0$ for all
  $r\in \mathbb{Z}$ and $\omega\in [0,2\pi]$.
\item[(ii)] $\{X_{t,n}^{(a)}\}_{t=1}^{n}$ is a (asymptotically) conditionally stationary node iff
$K_{r}^{(a,a)}(\omega) \equiv 0$ for all $r\neq 0$ and $\omega\in
[0,2\pi]$.
\item[(iii)] The edge $(a,b)$ is conditionally time-invariant iff asymptotically
$K_{r}^{(a,b)}(\omega) \equiv 0$ for all $r\neq 0$ and $\omega\in
[0,2\pi]$.
\end{itemize}
\end{proposition}
\noindent {\bf PROOF} in Appendix  \ref{sec:LSproof}. \hfill $\Box$

\vspace{2mm}

Note that the above result is asymptotic in rescaled time
($n\rightarrow \infty$). We make this precise in the following proposition
where  we show that $[{\bf
  K}_{n}(\omega_{k_1},\omega_{k_2})]_{a,b}$ closely approximates the 
Fourier coefficients $K_{k_1-k_2}^{(a,b)}(\omega_{k_2})$. 

\begin{proposition}\label{lemma:Hrab}
Suppose Assumptions \ref{assum:lambda}, \ref{assum:invcovarianceK} and \ref{assum:LS}
hold. Let $K_{r}^{(a,b)}(\cdot)$ be defined as in (\ref{eq:Hrab}). Then
\begin{eqnarray}
\label{eq:Kdef}
[{\bf K}_{n}(\omega_{k_1},\omega_{k_2})]_{a,b} &=& 
\frac{1}{n}\sum_{t=1}^{n}\exp\left(-\frac{i2\pi (k_1-k_2)t}{n}\right)\Gamma^{(a,b)}\left(\frac{t}{n};\omega_{k_2}\right)+
O\left(\frac{1}{n}\right).
\end{eqnarray}
Further, 
\begin{eqnarray}
\label{eq:Kdef2}
[{\bf K}_{n}(\omega_{k_1},\omega_{k_2})]_{a,b} &=& 
\left\{
\begin{array}{cc}
K_{k_{1}-k_{2}}^{(a,b)}(\omega_{k_2})+
O\left(\frac{1}{n}\right) & \textrm{ if }|k_1-k_2|\leq n/2 \\
K_{k_{1}-k_{2}-n}^{(a,b)}(\omega_{k_2})+
O\left(\frac{1}{n}\right) & \textrm{ if } n/2 < (k_1-k_2) < n \\
K_{k_{1}-k_{2}+n}^{(a,b)}(\omega_{k_2})+
O\left(\frac{1}{n}\right) & \textrm{ if } -n<(k_1-k_2)< -n/2 \\
\end{array}
\right.
\end{eqnarray}
where the $O(n^{-1})$ bound is uniform over all $1\leq r \leq n$ (and
$n$ is in rescaled time). 

Since $[{\bf K}_{n}(\omega_{k_1},\omega_{k_2})]_{a,b} = [{\bf
  K}_{n}(\omega_{k_2},\omega_{k_1})]_{b,a}^{*}$, then (\ref{eq:Kdef2})
can be replaced with $K_{k_{2}-k_{1}}^{(b,a)}(\omega_{k_1})^{*}$, 
$K_{k_{2}-k_{1}+n}^{(b,a)}(\omega_{k_1})^{*}$ and
$K_{k_{2}-k_{1}-n}^{(b,a)}(\omega_{k_1})^{*}$ respectively. 
\end{proposition}
\noindent {\bf PROOF} in Appendix  \ref{sec:LSproof}. \hfill $\Box$

\vspace{2mm}
\noindent Note we split $[{\bf
  K}_{n}(\omega_{k_1},\omega_{k_2})]_{a,b}$ into three separate cases
due to the circular wrapping of the DFT, which is most pronounced
when $\omega_{k_1}$ lies at the boundaries of the interval $[0,2\pi]$.

In \cite{p:dwi-sub-11} and
\cite{p:jen-sub-15} we showed that the Fourier transform
of the time-varying spectral density matrix ${\boldsymbol
  G}_{r}(\omega) = \int_{0}^{1}e^{-2\pi i r u}{\boldsymbol
  \Sigma}(u;\omega)du$ decayed to zero as $|r|\rightarrow \infty$ and was smooth
over $\omega$. In the following lemma we show that a similar result holds
for the Fourier transform of the inverse spectral density matrix.

\begin{proposition}[Properties of $K_{r}^{(a,b)}(\omega)$]\label{lemma:smoothness}
Suppose Assumption \ref{assum:LS} holds. 
Then for all $1\leq a,b\leq p$ we have 
\begin{eqnarray}
\label{eq:Hrdecay}
\sup_{\omega}|K_{r}^{(a,b)}(\omega)|\rightarrow 0 \textrm{ as } r\rightarrow \infty
\end{eqnarray}
and  $\sup_{\omega}|K_{r}^{(a,b)}(\omega)|\sim |r|^{-1}$. 
Furthermore,  for all
$\omega_{1},\omega_{2}\in [0,\pi]$ and $r\in \mathbb{Z}$
\begin{eqnarray}
\label{eq:Hrsmooth}
\left|K_{r}^{(a,b)}(\omega_{1}) -  K_{r}^{(a,b)}(\omega_{2}) \right|  \leq 
\left\{
\begin{array}{cc}
C|\omega_{1}-\omega_{2}| & r = 0 \\
C|r|^{-1}|\omega_{1}-\omega_{2}| & r \neq 0 \\
\end{array}
\right.
\end{eqnarray}
where $C$ is a finite constant that does not depend on $r$ or $\omega$.
\end{proposition}
\noindent {\bf PROOF} in Appendix  \ref{sec:LSproof}. \hfill $\Box$

\vspace{2mm}
\noindent The above results describe two important features in ${\bf
  K}_{a,b}$:
\begin{enumerate}
\item For a given subdiagonal $r$, $[{\bf K}]_{a,b}^{(r)}$
  changes smoothly along the subdiagonal, where  ${\bf K}_{a,b}^{(r)}$ denotes the $r$th subdiagonal
($-(n-1)\leq  r \leq (n-1)$).
Analogous to locally smoothing the periodogram, 
to estimate the entries of $[{\bf    K}]_{a,b}$ from the DFTs 
we use the smoothness property and frequencies in a local
neighbourhood to obtain multiple ``near replicates''.
\item For a given row $k$,  $[{\bf K}(\omega_{k},\omega_{k+r})]_{a,b}$
 is large when $r \mmod(n)$ is close to zero and decays the further it
 is from zero. 
\end{enumerate}
These observations motivate the regression method that we describe
below for learning the nonstationary network structure.

\subsection{Node-wise regression of the DFTs}\label{sec:nodeDFT}

In this section, we propose a method for estimating the entries of
$F_n^{*}D_{n}F_{n}$. 
The problem of learning the network structure from finite sample time
series is akin to the graphical model selection problem in GGM,
addressed by \cite{dempster1972covariance} for the low-dimensional and
\cite{p:mei-bue-06} for the high-dimensional setting. In particular, the neighborhood selection approach of \cite{p:mei-bue-06} regresses one component of a multivariate random vector on the other components with Lasso \citep{tibshirani1996regression}, and uses non-zero regression coefficients to select its neighborhood, i.e. the nodes which are conditionally noncorrelated with the given component. 

Assuming the multivariate time series
is locally stationary and satisfies Assumption \ref{assum:LS}, we show that the nonstationary network learning problem can be formulated 
in terms of a regression of DFTs at a specific
Fourier frequency on neighboring DFTs. 
Let $J^{(a)}_k$ denote the DFT of the time series $\{X^{(a)}_t\}_t$
at Fourier frequency $\omega_k$, as defined in \eqref{eq:Jdef}. 
We denote the $p$-dimensional
vector of DFTs at $\omega_k$ by $\underline{J}_k$, and use
$\underline{J}_k^{-(a)}$ to denote the $(p-1)$-dimensional vector
consisting of all the coordinates of $\underline{J}_k$ except $J^{(a)}_k$.

We define the space $\mathcal{G}_{n} = \overline{\textrm{sp}}(J_{k}^{(b)};1\leq k \leq n,1\leq b
\leq p)$  (note that the coefficients in this space can be
complex). Then 
\begin{eqnarray}
\label{eq:Jregress}
P_{\mathcal{G}_{n} - J_{k}^{(a)}}(J_{k}^{(a)}) =
  \sum_{b=1}^{p}\sum_{s=1}^{n}B_{(b,s)\shortarrow (a,k)}J_{s}^{(b)},
\end{eqnarray}
where we set $B_{(a,k)\shortarrow (a,k)} = 0$. Let 
\begin{eqnarray}
\label{eq:Jregress2}
\Delta_{k}^{(a)} = \var\left(J_{k}^{(a)}-P_{\mathcal{G}_{n} - J_{k}^{(a)}}(J_{k}^{(a)})\right).
\end{eqnarray}
The above allows us to rewrite the entries of $[{\bf
  K}_{n}(\omega_{k_1},\omega_{k_2})]_{a,b}$ in terms of regression
coefficients. In particular,
\begin{eqnarray}\label{eqn:Kn2B}
[{\bf K}_{n}(\omega_{k_1},\omega_{k_2})]_{a,b} = 
\left\{
\begin{array}{cc}
\frac{1}{\Delta_{k_1}^{(a)}} & k_{1}=k_{2}\textrm{ and }a=b \\
-\frac{1}{\Delta_{k_{1}}^{(a)}}B_{(b,k_{2})\shortarrow (a,k_{1})} &
                                                                    \textrm{
                                                                    otherwise
                                                                    }
\end{array}
\right..
\end{eqnarray}
Comparing the above with Proposition \ref{lemma:Hrab} 
 for $(a,k_1)\neq (b,k_{2})$ we have
\begin{eqnarray*}
B_{(b,k_{2})\shortarrow (a,k_{1})}  = B_{k_2-k_1,n}^{(b\shortarrow a)}(\omega_{k_1}) + O(n^{-1})
\textrm{ and } 
\Delta_{k}^{(a)} &=& [K^{(a,a)}_{0}(\omega_{k})]^{-1} + O(n^{-1}),
\end{eqnarray*}
where 
\begin{eqnarray}\label{eqn:B2H}
B_{r,n}^{(b\shortarrow a)}(\omega_{k}) = 
\left\{
\begin{array}{cc}
-K^{(a,a)}_{0}(\omega_{k})^{-1} K_{r}^{(b,a)}(\omega_{k})^{*}
& \textrm{ if }|r|\leq n/2, r\neq 0\\
-K^{(a,a)}_{0}(\omega_{k})^{-1}K_{r-n}^{(b,a)}(\omega_{k})^{*}
& \textrm{ if } n/2 < r < n \\
-K^{(a,a)}_{0}(\omega_{k})^{-1} K_{r+n}^{(b,a)}(\omega_{k})^{*}
 & \textrm{ if } -n<r< -n/2 \\
\end{array}
\right..
\end{eqnarray}
Thus by using Proposition \ref{lemma:smoothness} we have 
 \begin{eqnarray}
\label{eq:Bk1k2}
\left| B_{(b,k_{1}+r)\shortarrow (a,k_{1})} -B_{(b,k_{2}+r)\shortarrow
   (a,k_{2})}\right|   \leq
   A_{}|\omega_{k_{1}}-\omega_{k_{2}}| + O(n^{-1}),
\end{eqnarray}
where $A_{}$ is a finite constant. The benefit of these results is
in the estimation of the coefficients $B_{(b,k_{}+r)\shortarrow
  (a,k_{})}$. We recall (\ref{eq:Jregress}) can be expressed as 
\begin{eqnarray*}
P_{\mathcal{G}_{n} - J_{k_{}}^{(a)}}(J_{k_{}}^{(a)}) &=&
  \sum_{b=1}^{p}\sum_{r=-k_{}+1}^{n-k_{}}B_{(b,k_{}+r)\shortarrow
                                                       (a,k_{})}J_{k_{}+r}^{(b)},
\end{eqnarray*}
where the above is due to the periodic nature of $J^{(a)}_k$, which allows us to extend the definition
to frequencies outside $[0, 2\pi]$.
By using the near Lipschitz condition in (\ref{eq:Bk1k2}) if $k_{1}$ and $k_{2}$ are ``close'' then
the coefficients of the projections $P_{\mathcal{G}_{n} - J_{k_{1}}^{(a)}}(J_{k_{1}}^{(a)})$
and $P_{\mathcal{G}_{n} - J_{k_{2}}^{(a)}}(J_{k_{2}}^{(a)})$ will be
similar. This observation will allow us to estimate
$B_{(b,k+r)\shortarrow(a,k)}$ using the DFTs whose frequencies all
lie in the $M$-neighbourhood of $k$ (analogous to smoothing the
periodogram of stationary time series). We note that with these quasi
replicates the estimation would involve $(2M+1)$ (where $M<<n$) response variables and 
$pn-1$ regressors. Even with the aid of sparse estimation methods this is a large
number of regressors. However,  Proposition \ref{lemma:smoothness} allows us
to reduce the number of regressors in the regression. Since $|B_{(b,k+r)\shortarrow (a,k)}|\sim |r|^{-1}$ we
can truncate the projection to  a small number ($2\nu+1$) of regressors
about $\underline{J}_{k}$ to obtain the approximation
\begin{eqnarray*}
P_{\mathcal{G}_{n} - J_{k}^{(a)}}(J_{k}^{(a)}) &\approx&
  \sum_{b=1}^{p}\sum_{r=-\nu}^{\nu}B_{(b,k+r)\shortarrow
                                                       (a,k)}J_{k+r}^{(b)}.      
\end{eqnarray*}
Thus smoothness together with near sparsity of the coefficients make estimation of the
entries in the high-dimensional precision matrix $F_{n}^{*}D_{n}F_{n}$ feasible.

For a given choice of $M$ and $\nu$, and every value of $a, k$, we
define the $(2M+1)$-dimensional complex response vector
$\mathcal{Y}^{(a)}_k = (J^{(a)}_{k-M}, J^{(a)}_{k-M+1}, \ldots,
J^{(a)}_{k}, J^{(a)}_{k+1}, \ldots, J^{(a)}_{k+M})'$, and the $(2M+1)
\times ((2\nu+1)p - 1)$ 
dimensional complex design matrix
\begin{eqnarray*}
  \mathcal{X}^{(a)}_k =
  \left[
    \begin{array}{ccccccc}
      \underline{J}_{k-M-\nu}' & \ldots & \underline{J}_{k-M-1}' &
                                                                   (\underline{J}_{k-M}^{-(a)})' & \underline{J}_{k-M+1}' & \ldots & 
\underline{J}_{k-M+\nu}' \\
      \vdots & \vdots & \vdots & \vdots & \vdots & \vdots & \vdots \\
      \underline{J}_{k-\nu}' & \ldots & \underline  {J}_{k-1}' &
                                                                 (\underline{J}_{k}^{-(a)})' & \underline{J}_{k+1}' & \ldots & 
\underline{J}_{k+\nu}' \\
            \vdots & \vdots & \vdots & \vdots & \vdots & \vdots & \vdots \\
      \underline{J}_{k+M-\nu}' & \ldots & \underline{J}_{k+M-1}' &
                                                                   (\underline{J}_{k+M}^{-(a)})' &
 \underline{J}_{k+M+1}' & \ldots & \underline{J}_{k+M+\nu}'
    \end{array}
  \right].
\end{eqnarray*}
Then the estimator 
\begin{eqnarray*}
\hat{B}_{(.,.) \shortarrow (a,k)} &=&
\left(\hat{B}_{(1,k-\nu) \shortarrow (a,k)}, \ldots, \hat{B}_{(p,k-v)
    \shortarrow (a,k)}, \ldots, 
\hat{B}_{(1,k) \shortarrow (a,k)}, \ldots, \hat{B}_{(k-1,k) \shortarrow (a,k)}, \right.\\ 
&& \left. \hat{B}_{(k+1,k) \shortarrow (a,k)}, \hat{B}_{(p,k) \shortarrow
  (a,k)}, \ldots, \hat{B}_{(1,k+v) \shortarrow (a,k)}, \ldots,
\hat{B}_{(p,k+\nu) \shortarrow (a,k)} \right)'
\end{eqnarray*}
of $\{B_{(b,k+r)\shortarrow(a,k)};1\leq b \leq p,-\nu \leq r \leq \nu\}$
is obtained by solving the complex lasso optimization problem
\begin{eqnarray*}
\min_{\beta \in \mathbb{C}^{(2\nu+1)p-1}}\left[
 \frac{1}{2M+1} \left\| \mathcal{Y}^{(a)}_k - \mathcal{X}^{(a)}_k \beta \right\|^2_2 + \lambda \left\| \beta \right\|_1\right],
\end{eqnarray*}
where $\left\| \beta \right\|_1:= \sum_{j} |\beta_{j}|$, the sum of moduli
of all the (complex) coordinates, and $\lambda$ is a (real positive)
tuning parameter controlling the degree of regularization. 
It is well-known \citep{maleki2013asymptotic} that the above
optimization 
problem can be equivalently expressed as a group lasso optimization
over real variables,
 and can be solved using existing software. We use this property to
 compute the estimators in our numerical experiments.

From Proposition \ref{lemma:NetworkK} we 
observe that the problem of graphical model selection reduces to learning the
locations of large entries of ${\bf K}_n(\omega_{k_1}, \omega_{k_2})$ for different Fourier
frequencies $\omega_{k_1}, \omega_{k_2}$. Furthermore, from equation \eqref{eqn:Kn2B} and
\eqref{eqn:B2H} it is possible to learn the sparsity structure of 
${\bf K}_n(\omega_{k_1}, \omega_{k_2})$ from the regression
coefficients $B_{(b,k_1) \shortarrow (a,k_2)}$ (up to order $O(n^{-1})$). In particular, there is an edge $(a,b) \subseteq E$,
i.e. the components $a$ and $b$ are conditionally correlated, if
$B_{(b,k_1) \shortarrow (a,k_2)}$ is non-zero for some $k_1, k_2$ 
(within the locally stationary framework). Similarly, an edge between $a$ and $b$ is
conditionally time-varying if $B_{(b,k_1) \shortarrow (a,k_2)}$ is
non-zero for some $k_1\neq k_2$. In the above we have
  ignored the $O(n^{-1})$ terms.

In view of these connections, we define two quantities involving the estimated regression coefficients whose sparsity patterns encode information on the graph structure. In particular, we aggregate the estimated regression coefficients across different Fourier frequencies into two $p \times p$ weight matrices
\begin{eqnarray}\label{eqn:w_self_other}
\hat{W}_{self} &=& \left( \left( \sum_{k} |\hat{B}_{(b,k) \shortarrow (a,k)}|^2 \right) \right)_{1 \le a,b \le p} \\
\hat{W}_{other} &=& \left( \left( \sum_{k_1 \neq k_2} |\hat{B}_{(b,k_1) \shortarrow (a,k_2)}|^2 \right) \right)_{1 \le a,b \le p}
\end{eqnarray}
for graphical model selection in NonStGM. Two components $a$ and $b$ are
deemed conditionally noncorrelated if both the $(a,b)^{th}$ and the $(b,a)^{th}$ off-diagonal
elements of $\hat{W}_{self}+\hat{W}_{other}$ are small. In contrast, a node
$a$ is deemed conditionally stationary
if the $(a,a)^{th}$ element of $\hat{W}_{other}$ is small. Similarly,
an edge between $a$ and $b$ is deemed 
conditionally time-invariant if both the $(a,b)^{th}$ and the $(b,a)^{th}$
elements of $\hat{W}_{other}$ is small. Note that our node-wise
regression approach does not ensure that the estimated weight matrices
$\hat{W}$ are symmetric. 
However, following \citep{p:mei-bue-06}, one can formulate suitable
``and'' (or ``or'') rule to construct an undirected graph, 
where an edge $(a, b)$ is present if the $(a, b)^{th}$ and $(b,a)^{th}$ entries are both large (or at least one of them is large).

%% file: section5_Revision.tex
\section{Time-varying Vector Autoregressive Models}\label{sec:tvVAR}

In this section we link the structure of the
coefficients of the time-varying Vector Autoregressive (tvVAR) process
with the notion of conditional noncorrelation and conditional
stationarity. This gives a rigourous understanding of certain features
in a tvVAR model.

The time-varying VAR (tvVAR) model is often used to model nonstationarity (see
\cite{p:sub-70}, \cite{p:dah-00a}, \cite{p:dah-pol-06},
\cite{p:zha-20}, \cite{safikhani2020joint}). A time
series is said to have a  time-varying VAR$(\infty)$ representation if it
can be expressed as 
\begin{eqnarray}
\label{eq:tvVAR}
\underline{X}_{t} = \sum_{j=1}^{\infty}{\bf
  A}_{j}(t)\underline{X}_{t-j} + \underline{\varepsilon}_{t} \quad
  t\in \mathbb{Z}
\end{eqnarray}
where $\{\underline{\varepsilon}_{t}\}_{t}$ are i.i.d random vectors
with $\var[\underline{\varepsilon}_{t}] = {\boldsymbol \Sigma}$ and
$\Ex[\underline{\varepsilon}_{t}]=0$. For simplicity, we have centered
the time series as the focus is on the second order structure of the
time series. We assume that
(\ref{eq:tvVAR}) has a well defined time-varying moving average
representation as its solution (we show below that this allows the inverse covariance to be
expressed in terms of $\{{\bf A}_j(t)\}$). 
We show below that the inverse covariance matrix operator corresponding to
(\ref{eq:tvVAR}) has a simple form that can easily be deduced from the
VAR parameters. 

\subsection{The tvVAR model and the nonstationary network}\label{sec:tvD}

In this section we obtain an expression for $D$ in terms of the tvVAR
parameters. 

Let ${\bf C}_{t,\tau} =
\cov[\underline{X}_{t},\underline{X}_{\tau}]$ and $C$ denote the
corresponding covariance operator as defined in
(\ref{eq:CcovarianceOp}). Let  ${\bf H}$ denote the Cholesky
decomposition of ${\boldsymbol \Sigma}^{-1}$ such that
${\boldsymbol \Sigma}^{-1} = {\bf H}^{\prime}{\bf H}$ (where ${\bf H}^{\prime}$ denotes the transpose
of ${\bf H}$). To obtain $D$ we use the Gram-Schmidt
orthogonalisation.
We define the following matrices; 
\begin{eqnarray*}
\widetilde{{\bf A}}_{\ell}(t) = 
\left\{
\begin{array}{cc}
I_{p} & \ell=0 \\
-{\bf A}_{\ell}(t) & \ell >0 \\
0 & \ell<0
\end{array}
\right..
\end{eqnarray*}
Using $\{\widetilde{{\bf A}}_{\ell}(t)\}$
we define the  infinite dimensional, block, lower
triangular matrix $L$ where the $(t,\tau)$th block of $L$ is defined as $L_{t,\tau} =
{\bf H}\widetilde{{\bf A}}_{t-\tau}(t)$ for all $t,\tau \in \mathbb{Z}$. Define $X = (\ldots,
\underline{X}_{-1},\underline{X}_{0}, \underline{X}_{1},\ldots)^{}$, then
$L X$ is defined as 
\begin{eqnarray*}
 (LX)_{t} = {\bf H}\sum_{\ell=0}^{\infty}\widetilde{{\bf
  A}}_{\ell}(t)\underline{X}_{t-\ell} = {\bf H}\left(\underline{X}_{t}-
  \sum_{\ell=1}^{\infty}{\bf A}_{\ell}(t)\underline{X}_{t-\ell}\right)
  = {\bf H}\underline{\varepsilon}_{t} \quad t\in \mathbb{Z}.
\end{eqnarray*}
By definition of (\ref{eq:tvVAR}) it can be seen that $\{(LX)_{t} \}_{t}$ are uncorrelated random vectors
with $\var[(LX)_{t}] = I_{p}$. From this, it is clear that $L^{\prime}L$
is the inverse of a rearranged version of $C$. We use this to deduce the 
inverse  $D = C^{-1}$. We define ${\bf D}_{t,\tau}$  as
\begin{eqnarray}
\label{eq:DARt}
{\bf D}_{t,\tau}=\sum_{\ell=-\infty}^{\infty}\widetilde{{\bf
  A}}_{\ell}(t+\ell)^{\prime}\Sigma^{-1}\widetilde{{\bf A}}_{(\tau-t)+\ell}(t+\ell).
\end{eqnarray}
The inverse of $C$ is $D = (D_{a,b};1\leq a,b\leq p)$, where $D_{a,b}$ is defined
by substituting (\ref{eq:DARt}) into (\ref{eq:Dab}).

We now focus on the case ${\boldsymbol \Sigma}=I_{p}$ and derive conditions for
conditional noncorrelation and stationarity.  
In this case, the suboperators $D_{a,b}$ have the entries
\begin{eqnarray}
\label{eq:Dabttau3}
[D_{a,b}]_{t,t+r} = 
\left\{
\begin{array}{cc}
\sum_{\ell=1}^{\infty}\langle [{\bf
  A}_{\ell}(t+\ell)]_{\cdot,a}
[{\bf A}_{\ell+r}(t+\ell)]_{\cdot,b}\rangle
-\langle [I_{p}]_{\cdot,a},
[{\bf A}_{r}(t+\ell)]_{\cdot,b}\rangle,
& r \geq 0 \\
\sum_{\ell=1}^{\infty}\langle [{\bf A}_{\ell}(t+\ell)]_{\cdot,b},
[{\bf  A}_{\ell-r}(t+\ell)]_{\cdot,a}\rangle 
-\langle [I_{p}]_{\cdot,b},
[{\bf A}_{-r}(t+\ell)]_{\cdot,a}\rangle,
& r < 0 \\
\end{array}
\right.,
\end{eqnarray}
where ${\bf A}_{\cdot,a}$ denotes the $a^{th}$ column of the matrix ${\bf A}$ and
$\langle\cdot,\cdot \rangle$ the standard dot product on $\mathbb{R}^{p}$.
Using the above expression for $D_{a,b}$, the parameters of the tvVAR
model  can be connected to conditional
noncorrelation and conditional stationarity:
\begin{itemize}
\item[(i)] {\bf Conditional noncorrelation} If for all $\ell\in
  \mathbb{R}$ the non-zero entries in the columns $[\widetilde{{\bf
  A}}_{\ell}(t)]_{\cdot,a}$ and $[\widetilde{{\bf
  A}}_{\ell}(t)]_{\cdot,b}$ do not coincide, then
$\{X_{t}^{(a)}\}$ and $\{X_{t}^{(b)}\}$ are conditionally noncorrelated. 
\item[(ii)] {\bf Conditionally stationary node} If for all $\ell\in
  \mathbb{Z}$, the columns
$[\widetilde{{\bf A}}_{\ell}(t)]_{\cdot,a}$ do not depend on $t$ then the node $a$ is
conditionally stationary and the submatrix $D_{a,a}$ simplifies to
\begin{eqnarray*}
[D_{a,a}]_{t,t+r} = 
\sum_{\ell=1}^{\infty}\langle [{\bf
  A}_{\ell}(0)]_{\cdot,a}
[{\bf A}_{\ell+|r|}(0)]_{\cdot,a}\rangle
-\langle [I_{p}]_{\cdot,a},
[{\bf A}_{|r|}(0)]_{\cdot,a}\rangle \quad \textrm{ for all }r,t\in \mathbb{Z}.
\end{eqnarray*}
\item[(iii)] {\bf Conditionally time-invariant edge} If for all $\ell$
  and $r$ the dot products $\langle [{\bf A}_{\ell}(t)]_{\cdot,b},
[{\bf  A}_{\ell-r}(t)]_{\cdot,a}\rangle$ do not depend on $t$ and $[{\bf
  A}_{r}(t)]_{a,b}$ and $[{\bf A}_{r}(t)]_{b,a}$ does not depend
on $t$ then $D_{a,b}$ is Toeplitz where 
\begin{eqnarray*}
[D_{a,b}]_{t,t+r} = 
\left\{
\begin{array}{cc}
\sum_{\ell=1}^{\infty}\langle [{\bf
  A}_{\ell}(0)]_{\cdot,a}
[{\bf A}_{\ell+r}(0)]_{\cdot,b}\rangle
-\langle [I_{p}]_{\cdot,a},
[{\bf A}_{r}(0)]_{\cdot,b}\rangle,
& r \geq 0 \\
\sum_{\ell=1}^{\infty}\langle [{\bf A}_{\ell}(0)]_{\cdot,b},
[{\bf  A}_{\ell-r}(0)]_{\cdot,a}\rangle 
-\langle [I_{p}]_{\cdot,b},
[{\bf A}_{-r}(0)]_{\cdot,a}\rangle,
& r < 0 \\
\end{array}
\right. 
\end{eqnarray*}
for all $t\in \mathbb{Z}$.
\end{itemize}
There can arise situations where some $[\widetilde{{\bf
  A}}_{\ell}(t)]_{\cdot,a}$ and $[\widetilde{{\bf
  A}}_{\ell}(t)]_{\cdot,b}$ depend on $t$, but the corresponding node
or edge is conditionally stationary or time-invariant. This happens when there is a cancellation in the
entries of ${\bf A}_{\ell}(t)$. However, these cases are 
quite exceptional.  

In Appendix \ref{sec:proofAR} we state conditions
on the tvVAR process such that Assumptions
  \ref{assum:lambda}, \ref{assum:invcovarianceK} and \ref{assum:LS}
  are satisfied.

\begin{remark}[The time-varying AR approximation of locally stationary
  time series]
In \cite{p:kra-sub-22}, Theorem 3.3 it is shown that if a 
multivariate nonstationary time series satisfies certain
second order locally stationary conditions, then the time series
has a  tvAR$(\infty)$ representation with nearly smooth VAR
parameters i.e.
\begin{eqnarray*}
X_{t} - \sum_{j=1}^{\infty}\Phi_{j}(t/n)X_{t-j} \approx H(t/n)\varepsilon_{t},
\end{eqnarray*}
where $H(\cdot)$ is a lower triangular matrix, $H(\cdot)$ and
$\Phi_{j}(\cdot)$ are Lipschitz continuous 
and $\{\varepsilon_{t}\}_{t}$ are uncorrelated random variables with
$\var[\varepsilon_{t}]=I_p$. 
Using $\{\Phi_{j}(\cdot)\}_{j}$ and $H(\cdot)$ it would be possible to
determine the approximate network of a nonstationary time series based on  the conditions
(i,ii,iii) stated above. 
\end{remark}

%% file: section6_Revision.tex
\section{Numerical Experiments}\label{sec:simulation}

We demonstrate the applicability of node-wise regression in selecting NonStGM on two systems of multivariate time series, a small ($p=4$) dimensional tvVAR(1) process described in Example \ref{example:running}, and a large ($p=10$) dimensional tvVAR(1) process.

\subsection{Small System}

\begin{figure}[!t]
\begin{center}
  \includegraphics[width = 0.25\textwidth]{plotspaper/corrected_graph}
  \includegraphics[width=0.6\textwidth, trim = 0 0.5in 0 0.8in, clip]{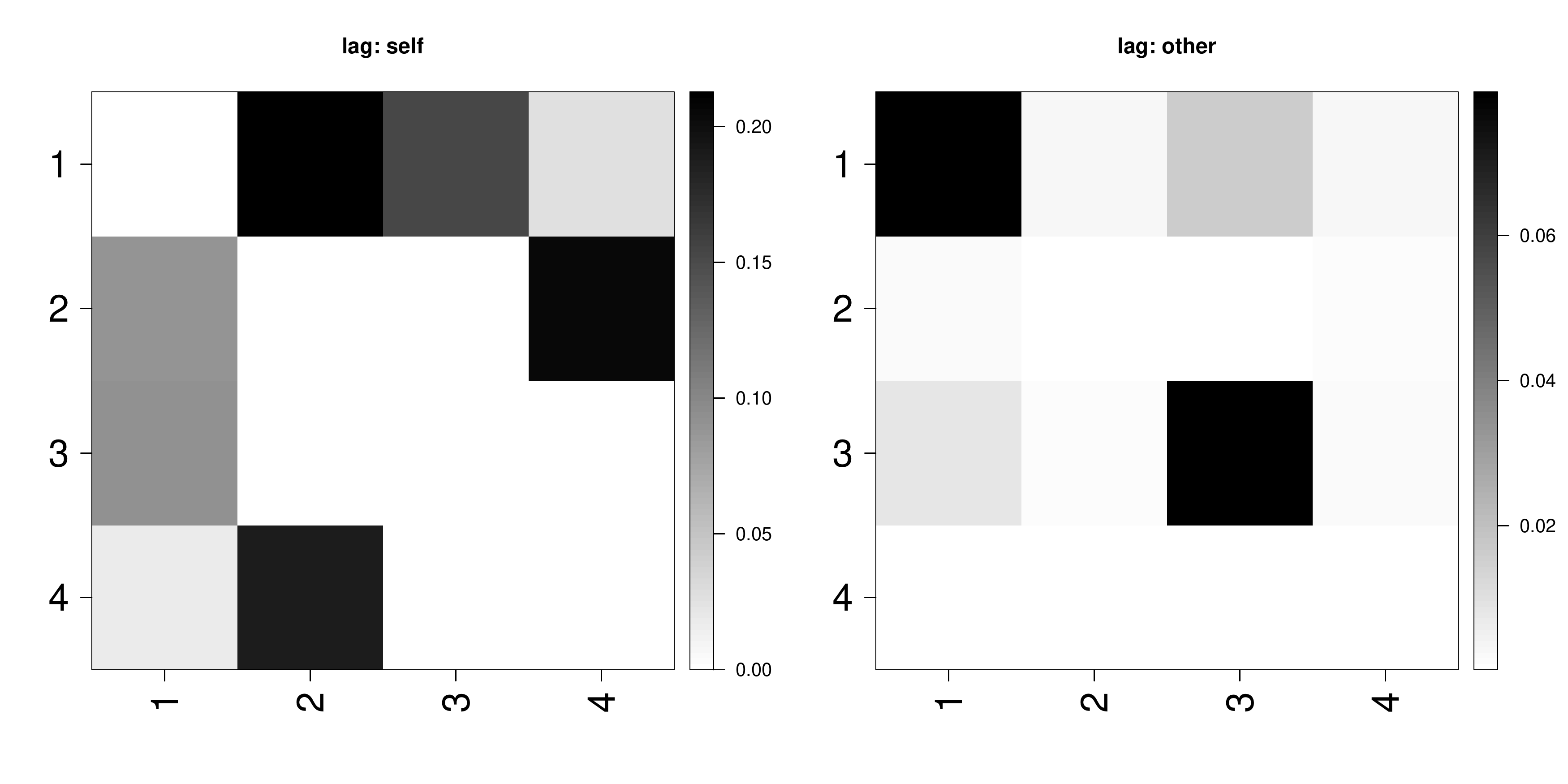}
\caption{NonStGM selection with node-wise regression for a $p=4$ dimensional system. [Left]: True graph structure. [Middle]: Heat map of $\hat{W}_{self}$ showing conditional noncorrelation between components $(1,2)$, $(1,3)$, $(1,4)$ and $(2,4)$. [Right]: Heat map of $\hat{W}_{other}$ showing conditional nonstationarity of nodes $1$ and $3$, and the conditionally time-varying edge $(1, 3)$. Results are aggregated over $20$ replicates.}
\label{fig:heatmap_small}
\end{center}
\end{figure}

We simulate the $p=4$ dimensional tvVAR(1) system described in 
Example \ref{example:running}, where all the time-invariant parameters set to $0.4$ and with $n=5000$ observations.
The two time-varying parameters $\alpha(t)$ and $\gamma(t)$ change from $-0.8$ to $0.8$ as $t$ varies from $1$ to $n$ according to the function $f(t) = -0.8 + 1.6\times e^{-5+10(t-1)/(n-1)}/(1 + e^{-5+10(t-1)/(n-1)})$.  
Using the results from Section \ref{sec:tvD}, 
 nodes $1, 3$ are conditionally nonstationay  and the edge $(1,3)$ is
 conditionally time-varying. On the other hand, the
nodes $2, 4$ are conditionally stationary
and the edge $(2,4)$, $(1,2)$ and (1,4) are conditionally
time-invariant. 
As Figure \ref{fig:example_p_4a} shows, these nuanced relationships are not prominent from the four time series trajectories. 

We perform node-wise regression of DFTs with $M = \lceil \sqrt{n} \rceil$ and $\nu = 1$. The tuning parameters in the individual group lasso regressions were selected using cross-validation. The estimated regression coefficients $\hat{B}$ were used to construct the weight matrices $\hat{W}_{self}$ and $\hat{W}_{other}$. The heat maps of these weight matrices, aggregated over $20$ replicates, are displayed in Figure \ref{fig:heatmap_small}.

The true graph structure (left) has two conditionally nonstationary nodes $1, 3$, and two stationary nodes $2, 4$. 
A heat map of $\hat{W}_{self}$ (middle) clearly shows the edges $(1,3), (1,2)$, $(2,4)$ and $(1,4)$
capturing conditional noncorrelation in the true graph structure. The
heat map of $\hat{W}_{other}$ (right) shows the conditionally
nonstationary nodes $1$ and $3$ on the diagonal. The conditionally
time-varying edge $(1,3)$ is also clearly visible on this heat map.

\subsection{Large System}

\begin{figure}[!t]
\begin{center}
\includegraphics[width = 0.3\textwidth]{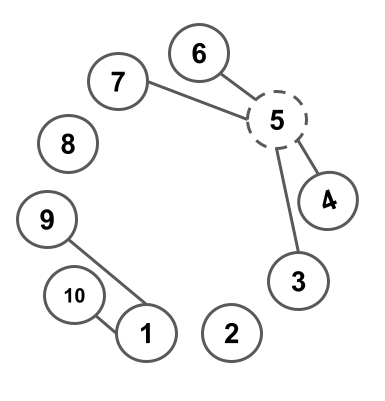}
\includegraphics[width=0.6\textwidth, trim = 0 0.5in 0 0.8in, clip]{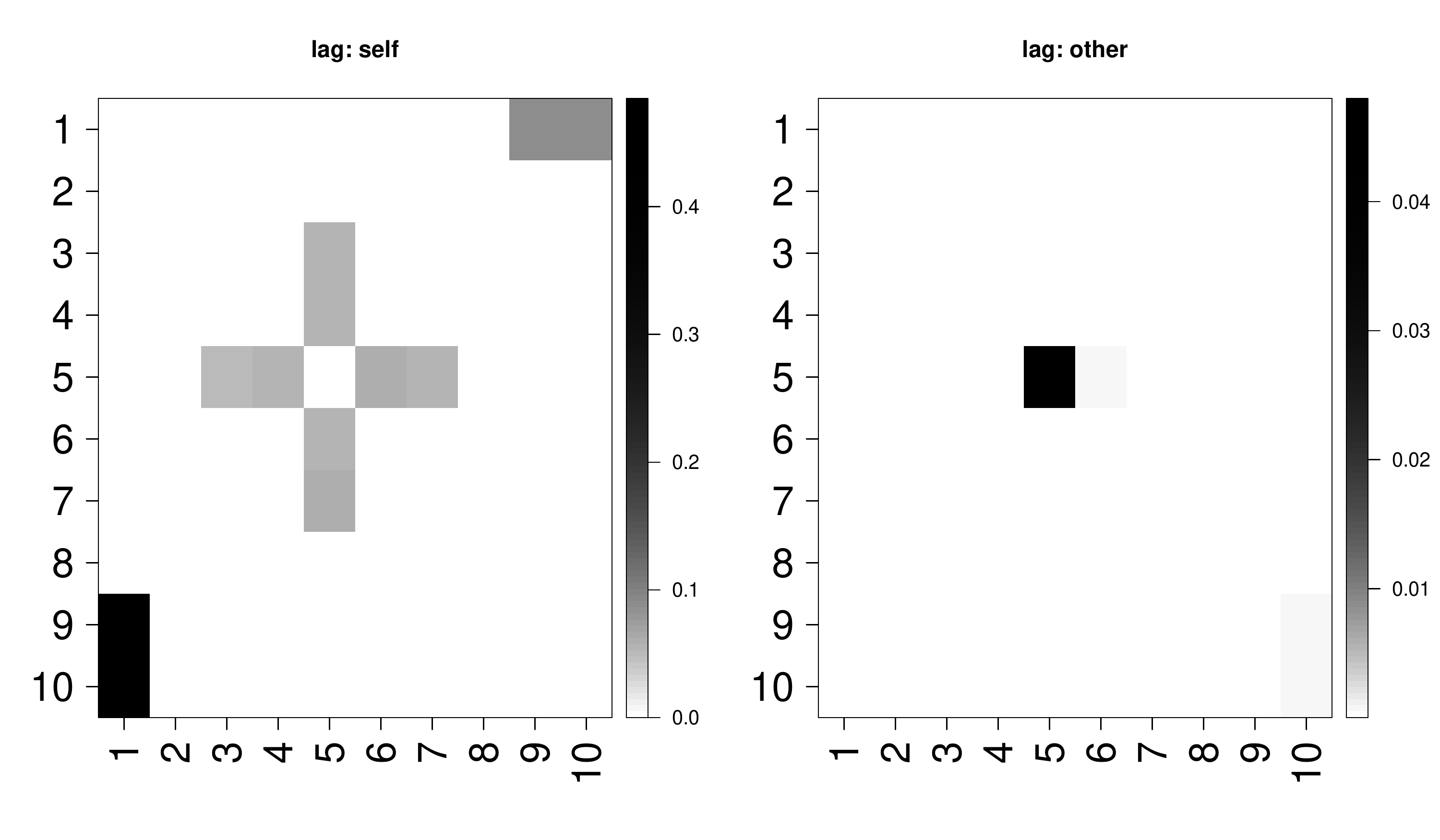}
\caption{NonStGM selection with node-wise regression for a $p=10$ dimensional system. [Left]: True graph structure. [Middle]: Heat map of $\hat{W}_{self}$ showing conditional noncorrelation captured by the edges. [Right]: Heat map of $\hat{W}_{other}$ showing conditional nonstationarity of node $5$.  Results are aggregated over $20$ replicates.}
\label{fig:heatmap_large}
\end{center}
\end{figure}

We now consider a larger system of $p=10$. The data generating process is tvVAR(1) $X_t = A(t) X_{t-1}+\varepsilon_t$. Here $\varepsilon_t \stackrel{i.i.d.}{\sim} N(0, I_{10})$. Non-zero time-invariant entries of the transition matrix $A(t)$ are constant functions as follows: $A_{j,j}(t) = 0.5$ for all $j \neq 5$, $A_{9,1}(t) = A_{10,1}(t) = A_{3,5}(t) = A_{4,5}(t) = A_{6,5(t)} = A_{7,5}(t) = 0.3$. The only time-varying entry is $A_{5,5}(t) = \alpha(t)$
, 
where $\alpha(t)$ decays exponentially from $0.7$ to $-0.7$ as $t$ varies from $1$ to $n$ according to the function $f(t) = 0.7  - 1.4\times e^{-5+10(t-1)/(n-1)}/(1 + e^{-5+10(t-1)/(n-1)})$. As we can see from the structure of $A(t)$ 
(and the true graph structure in the left panel of Figure \ref{fig:heatmap_large}), 
this network has two connected components and two isolated nodes ($2$ and $8$). These two nodes are independent of the other nodes, and are treated
as the ``control''.
The component consisting of $(1,9,10)$ is stationary (due to time invariant AR parameters). 
On the other hand, the component $(3,4,5,6,7)$ is nonstationary. However, the source of nonstationarity is node $5$ which permeates through to nodes $3,4,6$ and $7$. Thus the four nodes $3,4,6$ and $7$ are conditionally stationary (due to time-invariant parameters). 

We simulate $n = 15000$ observations from this system, and perform node-wise regression of DFTs with $M = \lceil \sqrt{n} \rceil$ and $\nu = 1$. The tuning parameters in the individual group lasso regressions were selected using cross-validation. The estimated regression coefficients $\hat{B}$ were used to construct the weight matrices $\hat{W}_{self}$ and $\hat{W}_{other}$. The heat maps of these weight matrices, aggregated over $20$ replicates, are displayed in Figure \ref{fig:heatmap_large}.

We observe that the edges for both components $(1,9,10)$ and $(3,4,5,6,7)$ 
are visible in the heat map of $\hat{W}_{self}$ (middle). As expected the isolated nodes do not show up.
The heat map of $\hat{W}_{other}$ (right) correctly identifies node $5$ as conditionally nonstationary. 

%% file: conclusions_acknowledgements.tex
\subsection*{Conclusion}

We introduced a general graphical modeling framework for describing
conditional relationships among the components of a multivariate
nonstationary time series using an undirected network. In this
network, absence of an edge corresponds to conditional noncorrelation
relationships, as is common in GGM and StGM. An additional node or edge
attribute (dashed or solid) further describes a newly introduced notion of \textit{conditional
  nonstationarity}, which can be used to provide a parsimonious
description of nonstationarity inherent in the overall system.
We showed that this framework is a natural generalization of
the existing GGM and StGM network. Under the locally stationary
framework, we proposed methods to
learn the nonstationary graph structure from finite-length time series
in the Fourier domain. Numerical experiments on simulated data
demonstrate the feasibility of our proposed method.

For stationary time series, there is well-established asymptotic theory 
for spectral density  matrix 
  estimators (see, e.g. \citep{woodroofe1967maximum, b:bri-01, wu2018asymptotic,
    rosuel2021asymptotic}). To estimate the inverse of moderate to high-dimensional spectral
  density matrices, penalized estimation methods for detecting non-zero off-diagonal
    entries \citep{fiecas2019spectral} have shown promise. 
These methods are based on learning the conditional correlation structure of the DFTs at 
different nodes at the same
frequency. Using the results in Section \ref{sec:nodeDFT}
we conjecture that the nonstationary
network can be estimated by learning 
 the non-zero coefficients of node-wise DFT regression across 
 different frequencies.  In future work, we hope to develop 
a complete statistical theory for graphical model estimation and
inference.

%% file: Appendix_Section2.tex
\section{Proofs for Section \ref{sec:framework}}\label{appendix:A}

\subsection{Proofs of results in Section \ref{sec:covariance}}\label{sec:covarianceproofs}

We first show that $\sigma^{2}_{a,t}\geq \lambda_{\inf}$. This ensures
that $\sigma^{2}_{a,t}>0$ and the operator $D$ is well defined (see (\ref{eq:Dttau})).  We recall that 
\begin{eqnarray*}
P_{\mathcal{H}-X_{t}^{(a)}}(X_{t}^{(a)}) = \sum_{\tau\in \mathbb{Z}}\sum_{b=1}^{p}\beta_{(\tau,b)\shortarrow (t,a)}X_{\tau}^{(b)}
\end{eqnarray*}
where $\beta_{(t,a)\shortarrow (t,a)}=0$. For all $(\tau,b)$ except
$(t,a)$ let $v_{(\tau,b)\shortarrow (t,a)} =
-\beta_{(\tau,b)\shortarrow (t,a)}$ and let
$v_{(t,a)\shortarrow (t,a)}=1$. For every $b \in \{1, \ldots, p\}$ define
$v^{(b)} = (v_{(\tau,b)\shortarrow (t,a)};\tau \in
\mathbb{Z})$ and $v =
\text{vec}[v^{(1)},\ldots,v^{(p)}]$. It is easily seen that 
\begin{eqnarray}
\label{eq:sigmaat}
\sigma_{a,t}^{2} = 
\Ex\left[X_{t}^{(a)} - P_{\mathcal{H}  -  X_{t}^{(a)}}(X_{t}^{(a)})
  \right]^{2} = \langle v, Cv \rangle.
\end{eqnarray}
Since $\|v\|_{2}\geq 1$, by Assumption \ref{assum:lambda}  we have 
\begin{eqnarray}
\label{eq:sigmaatbound}
\sigma_{a,t}^{2} = \langle v, Cv \rangle  \geq \lambda_{\inf}. 
\end{eqnarray}
We use this result and the notation above to prove Lemma \ref{lemma:inverse}.

\vspace{2mm}

\noindent {\bf PROOF of Lemma \ref{lemma:inverse}} 
For $1\leq b \leq  p$ we define the column vectors $X^{(b)} = (\ldots,X_{-1}^{(b)},X_{0}^{(b)},X_{1}^{(b)},\ldots)^{\prime}$
and $X = \textrm{vec}[X^{(1)},\ldots,X^{(p)}]$. Using the notation
introduced at the start of this section we have 
\begin{eqnarray*}
X_{t}^{(a)} - P_{\mathcal{H}-X_{t}^{(a)}}(X_{t}^{(a)})=X_{t}^{(a)}- \sum_{\tau\in
  \mathbb{Z}}\sum_{b=1}^{p}\beta_{(\tau,b)\shortarrow
  (t,a)}X_{\tau}^{(b)}   = \langle v, X\rangle.
\end{eqnarray*}
Since $P_{\mathcal{H}-X_{t}^{(a)}}(X_{t}^{(a)})$ minimises the mean
squared error $\Ex[X_{t}^{(a)} - Z]^{2}$ over all $Z\in
\mathcal{H}-X_{t}^{(a)}$ and 
$\Ex[X_{t}^{(a)} -
P_{\mathcal{H}-X_{t}^{(a)}}(X_{t}^{(a)})]^{2} = \sigma_{a,t}^{2}$ , this gives rise to the
normal equations
\begin{eqnarray*}
\cov[\langle v, X\rangle,X_{s}^{(c)}]=
\left\{
\begin{array}{cc}
0 & (c,s) \neq (a,t) \\
\sigma_{a,t}^{2} &  (c,s) = (a,t) \\
\end{array}
\right..
\end{eqnarray*}
Comparing the above with $D$, we observe that this proves $DC=I$ and
$CD=I$, thus
$D=C^{-1}$. To prove that $\|D\| = \lambda_{\inf}^{-1}$, we note that 
under Assumption \ref{assum:lambda}, $0<\lambda_{\inf} =
  \inf_{\|v\|_{2}=1, v\in \ell_{2,p}}\langle v, Cv \rangle \leq
  \sup_{\|v\|_{2}=1, v\in \ell_{2,p}}\langle v, Cv \rangle =
  \lambda_{\sup}<\infty$. Since $C$ is a self-adjoint operator,
  $\|C\|=\lambda_{\sup}$ and $\|D\| = \lambda_{\inf}^{-1}$. 

To prove that $\|D_{a,b}\|\leq \lambda_{\inf}^{-1}$ we first focus on the case $a=b$.
 Since $D_{a,a}$ are submatrices on the diagonal of $D$ and $0<\lambda_{\sup}^{-1} =
  \inf_{\|v\|_{2}=1, v\in \ell_{2,p}}\langle v, Dv \rangle \leq
  \sup_{\|v\|_{2}=1, v\in \ell_{2,p}}\langle v, Dv \rangle =
  \lambda_{\inf}^{-1}<\infty$, then it immediately follows that
\\*
$\lambda_{\sup}^{-1} \leq 
  \inf_{\|v\|_{2}=1, v\in \ell_{2}}\langle v, D_{a,a}v \rangle \leq
  \sup_{\|v\|_{2}=1, v\in \ell_{2}}\langle v, D_{a,a}v \rangle =
  \lambda_{\inf}^{-1}$. Therefore, since $D_{a,a}$ is self-adjoint (symmetric) we
   have $\|D_{a,a}\| \leq
\lambda_{\inf}^{-1}$. By a similar argument $\|D_{a,a}^{-1}\|\leq \lambda_{\sup}$. 

To prove the result for $a\neq b$ we focus on the sub-matrix
\begin{eqnarray*}
D_{\{a,b\}} = 
\left(
\begin{array}{cc}
D_{a,a} & D_{a,b}^{*} \\
D_{a,b} & D_{b,b} \\ 
\end{array}
\right).
\end{eqnarray*}
Using the same argument to prove that $\|D_{a,a}\|\leq
\lambda_{\inf}^{-1}$ it can be shown that $\|D_{\{a,b\}}\|\leq
\lambda_{\inf}^{-1}$. Thus for all
$v^{\prime}=(u^{(1)},u^{(2)})^{\prime}\in \ell_{2}^{2}$ we have
$\|D_{\{a,b\}}v\|\leq \lambda_{\inf}^{-1}\|v\|_{2}$. We use this bound
below.
We recall that an operator (matrix) $B$ is
bounded  if there exists a finite constant $K$ where for all $u\in
\ell_{2}$ we have $\|Bu\|\leq K\|u\|$, it follows that $\|B\|\leq
K$. Returning to $D_{a,b}$, we will show that $\|D_{a,b}u^{(1)}\|\leq \lambda_{\inf}^{-1}\|u^{(1)}\|_{2}$.
For all $u^{(1)}\in \ell_{2}$ we have 
\begin{eqnarray*}
\|D_{a,b}u^{(1)}\|_{2} \leq
  \sqrt{\|D_{a,a}u^{(1)}\|_{2}^{2}+\|D_{a,b}u^{(1)}\|_{2}^{2}}
 = \| D_{\{a,b\}}v\|_{2} \leq \|D_{\{a,b\}}\|\|v\|_{2}  = \|D_{\{a,b\}}\|\|u^{(1)}\|_{2}
\end{eqnarray*}
where $v^{\prime} = (u^{(1)},0)$. Thus $\|D_{a,b}\|  \leq
\lambda_{\inf}^{-1}$, as required.

Finally, to prove that 
$\sup_{t}\sum_{\tau\in \mathbb{Z}}\|{\bf  D}_{t,\tau}\|_{2}^{2}\leq p\lambda_{\inf}^{-2}$, we first prove that 
for every $t_0\in \mathbb{Z}$, we have 
$\sum_{a=1}^{p}\sum_{\tau\in
  \mathbb{Z}}[{\bf D}_{\tau,t_0}]_{a,1}^{2}\leq \lambda_{\inf}^{-2}$.
Define the sequence $v =\textrm{vec}
[u^{(1)},u^{(2)},\ldots,u^{(p)}]\in \ell_{2,p}$ where we set 
$[u^{(1)}]_{t_0} = 1$ and $[u^{(1)}]_{s} = 0$ for $s\neq
t_0$ and $u^{(a)} = 0$ (zero sequence) for all $a\neq 1$. Then by
definition of $v$ (which mainly consists of zeros except for one
non-zero entry) we have
\begin{eqnarray*}
Dv = \left(
\begin{array}{c}
D_{1,1}u^{(1)} \\
\vdots \\
D_{p,1}u^{(1)} \\ 
\end{array}
\right) =  \left(
\begin{array}{c}
([D_{1,1}]_{\tau,t_0};\tau\in \mathbb{Z}) \\
\vdots \\
([D_{p,1}]_{\tau,t_0}; \tau\in \mathbb{Z}) \\ 
\end{array}
\right). 
\end{eqnarray*}
Thus for every $t_0\in \mathbb{Z}$ we have $\|Dv\|_{2}^{2} = \sum_{a=1}^{p}\sum_{\tau\in
  \mathbb{Z}}[D_{a,1}]_{\tau,t_0}^{2}= \sum_{a=1}^{p}\sum_{\tau\in
  \mathbb{Z}}[{\bf D}_{\tau,t_0}]_{a,1}^{2}\leq \|D\|^{2}\|v\|_{2}^{2}\leq
\lambda_{\inf}^{-2}$. By the same argument for any $b\in
\{1,\ldots,p\}$ and $t\in \mathbb{Z}$
we have $\sum_{a=1}^{p}\sum_{\tau\in
  \mathbb{Z}}[{\bf D}_{\tau,t}]_{a,b}^{2}  \leq \|D\|^{2}\leq
\lambda_{\inf}^{-2}$, this gives $\sum_{\tau\in \mathbb{Z}}\|{\bf
  D}_{\tau,t}\|_{2}^{2}\leq p\lambda_{\inf}^{-2}$. This proves the
claim.
\hfill $\Box$

\vspace{1mm}
\noindent Many of the results in this section use the block
operator inversion identity  (see  \cite{b:tre-08}, page 35, and
\cite{p:ber-20}, Section 2.3)). For completeness we give the identity below.
As we are working with covariance matrix
operators  we focus on symmetric/self-adjoint matrices. Suppose 
\begin{eqnarray*}
G=\left(
\begin{array}{cc}
A & B \\
B^{*} & C
\end{array}
\right),
\end{eqnarray*}
and $G^{-1}$ exists, then
\begin{eqnarray}
\label{eq:inverseidentity}
G^{-1}=\left(
\begin{array}{cc}
P^{-1} & -P^{-1}BC^{-1} \\
-C^{-1}B^{*}P^{-1} & (C-B^{*}A^{-1}B)^{-1} \\
\end{array}
\right)
\end{eqnarray}
where $P = A-BC^{-1}B^{*}$. We mention that $P$ is the Schur
complement of $C$ of the matrix $G$. 

\subsection{Proof of results in Section \ref{sec:covD}}\label{sec:covDproofs}

There are different methods for proving the results in Section
\ref{sec:covD}. One method is to use the  properties of projections the
other is to use decompositions of the infinite dimensional matrix
$C$. In this section we take the matrix decomposition route, as similar matrix
decompositions form the core of the proofs in Section
\ref{sec:covtimeD}.  

As the results in Section \ref{sec:covD} concern the 
partial covariance between $X_{t}^{(a)}$ and $X_{\tau}^{(b)}$ given
\emph{all} the other random variables, we will consider a permuted version of
$C$ and its inverse, where we bring the covariance structure of
$(X_{t}^{(a)},X_{\tau}^{(b)})$ to the top left hand corner of the
matrix. To avoid introducing new notation we label these
permuted matrix operators as $C$ and $D$. 
The variance of $(X_{t}^{(a)},X_{\tau}^{(b)})$ is 
\begin{eqnarray*}
\var\left[ 
\begin{array}{c}
X_{t}^{(a)} \\
X_{\tau}^{(b)} \\
\end{array}
\right] = 
\widetilde{C}_{1,1} = 
\left(
\begin{array}{cc}
[C_{a,a}]_{t,t} & [C_{a,b}]_{t,\tau} \\
{} [C_{a,b}]_{t,\tau} & [C_{b,b}]_{\tau,\tau} \\
\end{array}
\right).
\end{eqnarray*}
This embeds in the top left hand side of the operator $C$, where 
\begin{eqnarray}
\label{eq:C11}
C &=& 
\left(
\begin{array}{cc}
\widetilde{C}_{1,1} & \widetilde{C}_{1,2} \\
 \widetilde{C}_{2,1} & \widetilde{C}_{2,2} \\
\end{array}
\right)
\end{eqnarray}
with $\widetilde{C}_{1,2} = \{[C_{c,e}]_{u,v};(c,u)\in \{(a,t),(b,\tau)\},
(e,v)\notin \{(a,t),(b,\tau)\}\}$, 
$\widetilde{C}_{2,1}=\widetilde{C}_{1,2}^{*}$ and  \\
$\widetilde{C}_{2,2} =\{[C_{c,e}]_{u,v};(c,u),(e,v)\notin \{(a,t),(b,\tau)\}\}$ (we have
used the tilde notation in
$\widetilde{C}_{i,j}$ to distinguish it from $C_{a,b}$). It is well
known that the Schur complement encodes the partial
covariance. Applying this to $(X_{t}^{(a)},X_{\tau}^{(b)})$, the Schur
complement of  $\widetilde{C}_{1,1}$ in $C$ is 
\begin{eqnarray}
\label{eq:PPPSec23}
\var\left[
\left(
\begin{array}{c}
X_{t}^{(a)} - P_{\mathcal{H}  -
  (X_{t}^{(a)},X_{\tau}^{(b)})}(X_{t}^{(a)}) \\
X_{\tau}^{(b)} - P_{\mathcal{H}  -
  (X_{t}^{(a)},X_{\tau}^{(b)})}(X_{\tau}^{(b)}) \\
\end{array}
\right)
\right] = \widetilde{C}_{1,1} - \widetilde{C}_{1,2}\widetilde{C}_{2,2}^{-1}\widetilde{C}_{2,1} = P.
\end{eqnarray}
The above matrix (which we label as $P$) forms an important part of
all the proofs in this section. As the entries in the variance matrix  on the left hand side
of (\ref{eq:PPPSec23}) is quite long we replace it with some shorter
notation for the conditional variances and covariances. Comparing (\ref{eq:PPPSec23})
with (\ref{eq:partial1}) we observe that
the off-diagonal is $\rho_{t,\tau}^{(a,b)} = \cov[X_{t}^{(a)} - P_{\mathcal{H}  -
  (X_{t}^{(a)},X_{\tau}^{(b)})}(X_{t}^{(a)}), X_{\tau}^{(b)} - P_{\mathcal{H}  -
  (X_{t}^{(a)},X_{\tau}^{(b)})}(X_{\tau}^{(b)})]$. However, the
diagonal of $P$ has not been defined in Section \ref{sec:covD}. As this is
the partial variance of $X_{t}^{(a)}$ after conditioning on everything
\emph{but} $X_{t}^{(a)}$ and $X_{\tau}^{(b)}$ we use the notation
\begin{eqnarray}
\label{eq:rhottconditional}
\rho_{t,t}^{(a,a)|\shortminus \{(a,t),(b,\tau)\}} = 
\var[X_{t}^{(a)} -  P_{\mathcal{H}  -
  (X_{t}^{(a)},X_{\tau}^{(b)})}(X_{t}^{(a)})].
\end{eqnarray}
To avoid confusion, we mention that this is different to the time series partial covariance
defined in Section \ref{sec:covtimeD}, where
$\rho_{t,t}^{(a,a)|\shortminus\{a,b\}} $ is the partial variance of
$X_{t}^{(a)}$ after conditioning on all the other time series but time
series $X^{(a)}$ and $X^{(b)}$. Using the new notation we have
\begin{eqnarray}
 \label{eq:var3a}
\var
\left(
\begin{array}{c}
X_{t}^{(a)} - P_{\mathcal{H}  -
  (X_{t}^{(a)},X_{\tau}^{(b)})}(X_{t}^{(a)}) \\
X_{\tau}^{(b)} - P_{\mathcal{H}  -
  (X_{t}^{(a)},X_{\tau}^{(b)})}(X_{\tau}^{(b)}) \\
\end{array}
\right) &=& 
\left(
\begin{array}{ll}
\rho_{t,t}^{(a,a)|\shortminus \{(a,t),(b,\tau)\}} & \rho_{t,\tau}^{(a,b)}\\
\rho_{\tau,t}^{(b,a)} & \rho_{\tau,\tau}^{(b,b)|\shortminus \{(a,t),(b,\tau)\}} \nonumber\\
\end{array}
\right) \\
&=& \widetilde{C}_{1,1} - \widetilde{C}_{1,2}\widetilde{C}_{2,2}^{-1}\widetilde{C}_{2,1} = P,
\end{eqnarray}
where $\rho_{t,\tau}^{(a,b)}$ and $\rho_{t,t}^{(a,a)|\shortminus
  \{(a,t),(b,\tau)\}}$  are defined in (\ref{eq:partial1}) and 
 (\ref{eq:rhottconditional}) respectively. 

Next we relate $P$ to the inverse $C^{-1}$.
Using the block operator inversion (see (\ref{eq:inverseidentity})) we have
\begin{eqnarray}
\label{eq:schurSmall}
D &=& 
\left(
\begin{array}{cccc}
P^{-1} & -P^{-1}\widetilde{C}_{1,2}\widetilde{C}_{2,2}^{-1} \\
-\widetilde{C}_{2,2}^{-1}\widetilde{C}_{2,1}P^{-1}
       &  
(\widetilde{C}_{2,2} -  \widetilde{C}_{2,1}\widetilde{C}_{1,1}^{-1}\widetilde{C}_{1,2})^{-1} \\
\end{array}
\right)
\end{eqnarray}
where $P$ is defined in (\ref{eq:PPPSec23}). Comparing $P^{-1}$ with
the upper left block of $D$, 
we connect the conditional variance of $(X_{t}^{(a)},X_{\tau}^{(b)})$ to the entries of $D$. In particular
\begin{eqnarray}
 \label{eq:var3}
\left(
\begin{array}{ll}
\rho_{t,t}^{(a,a)|\shortminus \{(a,t),(b,\tau)\}} & \rho_{t,\tau}^{(a,b)}\\
\rho_{\tau,t}^{(b,a)} & \rho_{\tau,\tau}^{(b,b)|\shortminus \{(a,t),(b,\tau)\}} \\
\end{array}
\right) &=&\left(
\begin{array}{cc}
[D_{a,a}]_{t,t} & [D_{a,b}]_{t,\tau} \\
{} [D_{a,b}]_{t,\tau} & [D_{b,b}]_{\tau,\tau} \\
\end{array}
\right)^{-1}.
\end{eqnarray} 
The identity (\ref{eq:var3}) forms an important component in the
proofs below. 

Before we state the next lemma, we require the following notation for the
partial correlation
between $X_{t}^{(a)}$ and $X_{\tau}^{(b)}$ 
 \begin{eqnarray}
\label{eq:partialA}
\phi_{t,\tau}^{(a,b)} = \cor\left[X_{t}^{(a)} - P_{\mathcal{H}  -
  (X_{t}^{(a)},X_{\tau}^{(b)})}(X_{t}^{(a)}), X_{t}^{(b)} - P_{\mathcal{H}  -
  (X_{t}^{(a)},X_{\tau}^{(b)})}(X_{\tau}^{(b)})\right].
\end{eqnarray} 
\begin{lemma}\label{lemma:partial}
Let $\beta_{(\tau,b)\shortarrow (t,a)}$, $\sigma^{2}_{t,a}$, $\rho_{t,t}^{(a,a)|\shortminus \{(a,t),(b,\tau)\}}$ and $\phi_{t,\tau}^{(a,b)}$ be
defined as in (\ref{eq:projectXat}), (\ref{eq:rhottconditional}) and 
(\ref{eq:partialA}).  Suppose Assumption \ref{assum:lambda} holds.  Then 
\begin{eqnarray}
\label{eq:varijvarxy}
\frac{\rho_{t,t}^{(a,a)|\shortminus
  \{(a,t),(b,\tau)\}}}{\rho_{\tau,\tau}^{(b,b)|\shortminus
  \{(a,t),(b,\tau)\}}} = \frac{\sigma^{2}_{t,a}}{\sigma^{2}_{\tau,b}}
\end{eqnarray}
and 
\begin{eqnarray}
\label{eq:betacorr}
\beta_{(\tau,b)\shortarrow (t,a)} 
&=&\phi_{t,\tau}^{(a,b)} \times \frac{\sigma_{a,t}}{\sigma_{b,\tau}}. 
\end{eqnarray}
\end{lemma}
PROOF. The proof of (\ref{eq:varijvarxy}) is based on comparing
$[D_{a,a}]_{t,t}/[D_{b,b}]_{\tau,\tau}$ and the ratio of the diagonal
entries of the conditional variance in (\ref{eq:var3})
\begin{eqnarray}
\label{eq:Ddiagonal}
&&\left(
\begin{array}{ll}
\rho_{t,t}^{(a,a)|\shortminus \{(a,t),(b,\tau)\}} & \rho_{t,\tau}^{(a,b)}\\
\rho_{\tau,t}^{(b,a)} & \rho_{\tau,\tau}^{(b,b)|\shortminus \{(a,t),(b,\tau)\}} \\
\end{array}
\right) \nonumber\\ 
&=& 
\frac{1}{[D_{a,a}]_{t,t}[D_{b,b}]_{\tau,\tau} - [D_{a,b}]_{t,\tau}^{2}}
\left(
\begin{array}{cc}
[D_{b,b}]_{\tau,\tau} & -[D_{a,b}]_{t,\tau} \\
{} -[D_{a,b}]_{t,\tau} & [D_{a,a}]_{t,t} \\
\end{array}
\right).
\end{eqnarray} 
We recall from (\ref{eq:Dttau}) that
\begin{eqnarray}
[D_{a,a}]_{t,t} = \frac{1}{\sigma_{a,t}^{2}} 
\textrm{ and }
 [D_{b,b}]_{\tau,\tau} = \frac{1}{\sigma_{b,\tau}^{2}}
\quad
\Rightarrow 
\frac{\sigma_{a,t}^{2}}{\sigma_{b,\tau}^{2}} = \frac{[D_{b,b}]_{\tau,
  \tau}}{[D_{a,a}]_{t,t}}. \label{eq:sigmaratio1}
\end{eqnarray}
Furthermore, by comparing the 
entries in (\ref{eq:Ddiagonal}) we have 
\begin{eqnarray*}
\rho_{t,t}^{(a,a)|\shortminus \{(a,t),(b,\tau)\}}
 &=&
     \frac{1}{[D_{a,a}]_{t,t}[D_{b,b}]_{\tau,\tau}-[D_{a,b}]_{t,\tau}^{2}}[D_{b,b}]_{\tau,
     \tau}
  \\
\textrm{ and }\rho_{\tau,\tau}^{(b,b)|\shortminus \{(a,t),(b,\tau)\}}
 &=&
     \frac{1}{[D_{a,a}]_{t,t}[D_{b,b}]_{\tau,\tau}-[D_{a,b}]_{t,\tau}^{2}}[D_{a,a}]_{t,t}.
\end{eqnarray*}
Thus evaluating ratio of the above gives
\begin{eqnarray}
\label{eq:sigmaratio2}
\frac{\rho_{t,t}^{(a,a)|\shortminus
  \{(a,t),(b,\tau)\}}}{\rho_{\tau,\tau}^{(b,b)|\shortminus
  \{(a,t),(b,\tau)\}}} =  \frac{[D_{b,b}]_{\tau,
  \tau}}{[D_{a,a}]_{t,t}}.
\end{eqnarray}
Comparing (\ref{eq:sigmaratio1}) and (\ref{eq:sigmaratio2}) gives (\ref{eq:varijvarxy}).

To prove (\ref{eq:betacorr}) we decompose the projection of
$P_{\mathcal{H}-X_{t}^{(a)}}(X_{t}^{(a)})$ in terms of its projections
onto the two space $\mathcal{H}-(X_{t}^{(a)},X_{\tau}^{(b)})$  and $\overline{\textrm{sp}}(X_{\tau}^{(b)} -
    P_{\mathcal{H}-(X_{t}^{(a)},X_{\tau}^{(b)})}(X_{\tau}^{(b)}))$. These
    two spaces are orthogonal and lead to a 
 simple expression for the coefficient $\beta_{(\tau,b)\shortarrow (t,a)}$;
\begin{eqnarray}
\label{eq:two-stage}
P_{\mathcal{H}-X_{t}^{(a)}}(X_{t}^{(a)}) &=&  \sum_{\tau\in \mathbb{Z}}\sum_{b=1}^{p}\beta_{(\tau,b)\shortarrow (t,a)}X_{\tau}^{(b)}\nonumber\\
&=& P_{\mathcal{H}-(X_{t}^{(a)},X_{\tau}^{(b)})}(X_{t}^{(a)}) +
    \beta_{(\tau,b)\shortarrow (t,a)}\left[X_{\tau}^{(b)} -
    P_{\mathcal{H}-(X_{t}^{(a)},X_{\tau}^{(b)})}(X_{\tau}^{(b)})\right].
\end{eqnarray}
Using the orthogonality of the two projections we have
\begin{eqnarray*}
\beta_{(\tau,b)\shortarrow (t,a)} = \frac{\cov\left[X_{t}^{(a)}, X_{\tau}^{(b)} -
  P_{\mathcal{H}-(X_{t}^{(a)},X_{\tau}^{(b)})}(X_{\tau}^{(b)})\right]}{
\var(X_{\tau}^{(b)} -
  P_{\mathcal{H}-(X_{t}^{(a)},X_{\tau}^{(b)})}(X_{\tau}^{(b)}))} = 
\frac{\rho_{t,\tau}^{(a,b)}}{\rho_{t,t}^{(a,a)|\shortminus
  \{(a,t),(b,\tau)\}}}.
\end{eqnarray*}
Replacing
  the covariance $\rho_{t,\tau}^{(a,b)}$ in $\beta_{(\tau,b)\shortarrow (t,a)}$
 with its correlation $\phi_{t,\tau}^{(a,b)}$ gives 
\begin{eqnarray}
\label{eq:alphabeta}
\beta_{(\tau,b)\shortarrow (t,a)} &=& \phi_{t,\tau}^{(a,b)}
\sqrt{
\frac{
\rho_{t,t}^{(a,a)|\shortminus \{(a,t),(b,\tau)\}}
}
{
\rho_{\tau,\tau}^{(b,b)|\shortminus  \{(a,t),(b,\tau)\}} 
}
}.
\end{eqnarray}
This links the partial correlation to the projection coefficients.
Finally, we substitute the identity (\ref{eq:varijvarxy}) into
(\ref{eq:alphabeta})  to give 
\begin{eqnarray*}
\beta_{(\tau,b)\shortarrow (t,a)} 
&=&\phi_{t,\tau}^{(a,b)} \times \frac{\sigma_{a,t}}{\sigma_{b,\tau}}.
\end{eqnarray*}
This proves (\ref{eq:betacorr}). 
\hfill $\Box$

\vspace{1mm}
\noindent We use the above to prove Lemma \ref{lemma:D}. 

\noindent {\bf PROOF of Lemma \ref{lemma:D}} 
By using  (\ref{eq:betacorr}) we connect $\phi_{t,\tau}^{(a,b)}$ to the precision
matrix. Since 
\begin{eqnarray*}
\beta_{(\tau,b)\shortarrow (t,a)} 
&=&\phi_{t,\tau}^{(a,b)} \times \frac{\sigma_{a,t}}{\sigma_{\tau,b}}
\end{eqnarray*}
and by definition of $D_{a,b}$ in (\ref{eq:Dab}) we have 
\begin{eqnarray*}
\rho_{t,\tau}^{(a,b)} = -\frac{[D_{a,b}]_{t,\tau}}{\sqrt{[D_{a,a}]_{t,t} [D_{b,b}]_{\tau,\tau}}}.
\end{eqnarray*}
This proves (\ref{eq:corpartial}). The proof of (\ref{eq:covpartial})
immediately follows from (\ref{eq:var3}). \hfill $\Box$

\vspace{2mm}
\noindent {\bf PROOF of Proposition \ref{lemma:networkpartial}} 
The proof hinges on the identity in (\ref{eq:var3}) for the separate
cases $a=b$ and $a\neq b$.
For the case $a\neq b$ and using (\ref{eq:var3}) 
it is clear that $\rho_{t,\tau}^{(a,b)}=0$ iff
$[D_{a,b}]_{t,\tau} = 0$. Thus $D_{a,b}=0$ iff for all $t$ and $\tau$, $\rho_{t,\tau}^{(a,b)}=0$,
this proves (i).

To prove (ii), we use  (\ref{eq:var3}) with $a=b$ and compare the entries of 
\begin{eqnarray}
&& \left(
\begin{array}{ll}
\rho_{t,t}^{(a,a)|\shortminus \{(a,t),(a,\tau)\}} & \rho_{t,\tau}^{(a,a)}\\
\rho_{\tau,t}^{(a,a)} & \rho_{\tau,\tau}^{(a,a)|\shortminus \{(a,t),(a,\tau)\}} \\
\end{array}
\right) \nonumber\\
&=& \frac{1}{[D_{a,a}]_{t,t}[D_{a,a}]_{\tau,\tau} -
     [D_{a,a}]_{t,\tau}^{2}}
\left(
\begin{array}{cc}
[D_{a,a}]_{\tau,\tau} & -[D_{a,a}]_{t,\tau} \\
{} -[D_{a,a}]_{t,\tau} & [D_{a,a}]_{t,t} \\
\end{array}
\right).\label{eq:DDDaaa}
\end{eqnarray}
We first show that  if $D_{a,a}$ is a symmetric, Toeplitz matrix, then 
$\rho_{t,\tau}^{(a,a)}$ is shift invariant (depends only on $t-\tau$).
If $D_{a,a}$ is a symmetric,  Toeplitz matrix, using that
$\rho_{t,t}^{(a,a)} =1/[D_{a,a}]_{t,t} = 1/[D_{a,a}]_{0,0}$ it is clear that $\rho_{t,t}^{(a,a)}$ does not depend on
$t$. We now study $\rho_{t,\tau}^{(a,a)}$ when $t\neq \tau$. 
To show that $\rho_{t,\tau}^{(a,a)}$ only depends on $|t-\tau|$
we use that $[D_{a,a}]_{t,\tau} = [D_{a,a}]_{0,\tau-t} =
[D_{a,a}]_{0,t-\tau}$ (due to $D_{a,a}$ being Toeplitz). Comparing the
off-diagonal entries on the left and right hand side of 
(\ref{eq:DDDaaa}) it follows that for all $t$ and $\tau$
$\rho_{t,\tau}^{(a,a)} = \rho_{0,t-\tau}^{(a,a)} = \rho_{0,\tau-t}^{(a,a)}$.

Next we show the converse, that is if for all $t$ and $\tau$; $\rho_{t,\tau}^{(a,a)} =
\rho_{0,t-\tau}^{(a,a)} = \rho_{0,\tau-t}^{(a,a)}$ then $D_{a,a}$ is a symmetric,
Toeplitz matrix. First the diagonal, 
since $[D_{a,a}]_{t,t}
= 1/\rho_{t,t}^{(a,a)} =  1/\rho_{0,0}^{(a,a)}$ it is clear that the diagonal $D_{a,a}$
does not depend on $t$. Next we show that if for all $t$ and $\tau$; $\rho_{t,\tau}^{(a,a)} =
\rho_{0,t-\tau}^{(a,a)} = \rho_{0,\tau-t}^{(a,a)}$, then 
\begin{itemize}
\item[(a)] $[D_{a,a}]_{t,\tau}$ only
depends on $|t-\tau|$. 
\item[(b)] The conditional variance 
$\rho_{t,t}^{(a,a)|\shortminus \{(a,t),(a,\tau)\}}$ only depends on
$|t-\tau|$.
Note that this is not in the statement of the theorem, but is a
useful by product of the proof. 
\end{itemize}
Comparing the entries of the matrices in (\ref{eq:DDDaaa}) we have 
\begin{eqnarray}
\rho_{t,t}^{(a,a)|\shortminus \{(a,t),(a,\tau)\}} &=& \frac{[D_{a,a}]_{\tau,\tau}}{[D_{a,a}]_{t,t}[D_{a,a}]_{\tau,\tau} -
     [D_{a,a}]_{t,\tau}^{2}}, \label{eq:rho111}\\
 \rho_{t,\tau}^{(a,a)} &=& \frac{-[D_{a,a}]_{t,\tau}}{[D_{a,a}]_{t,t}[D_{a,a}]_{\tau,\tau} -
     [D_{a,a}]_{t,\tau}^{2}},\label{eq:rho112}\\
\textrm{ and }0 &<& [D_{a,a}]_{t,t}[D_{a,a}]_{\tau,\tau} -     [D_{a,a}]_{t,\tau}^{2}
                    \textrm{ (since this is the determinant).} \label{eq:rho113}
\end{eqnarray}
We first show that $[D_{a,a}]_{t,\tau}$ only depends on $|t-\tau|$. To
reduce notation, we set the entries on the diagonal of $D_{a,a}$ to
$\theta=[D_{a,a}]_{t,t}$  and let
$\theta_{t-\tau} = \rho_{t,\tau}^{(a,a)}$. Substituting this into
(\ref{eq:rho112}) gives 
\begin{eqnarray*}
\theta_{t-\tau}(\theta^{2}- [D_{a,a}]_{t,\tau}^{2}) &=& -[D_{a,a}]_{t,\tau}.
\end{eqnarray*}
The above is quadratic equation in $[D_{a,a}]_{t,\tau}$. 
Thus we can express $[D_{a,a}]_{t,\tau}$ in terms of $\theta$ and $\theta_{t-\tau}$;
\begin{eqnarray*}
[D_{a,a}]_{t,\tau} = \frac{-1 + \sqrt{1+4\theta_{t-\tau}^{2}\theta^{2}}}{2\theta_{t-\tau}}.
\end{eqnarray*}
Note that $-1 + \sqrt{1+4\theta_{t-\tau}^{2}\theta^{2}}$ is part of
the solution and not $-1 - \sqrt{1+4\theta_{t-\tau}^{2}\theta^{2}}$ due to
the positivity condition in (\ref{eq:rho113}).
This proves that $D_{a,a}$ is a symmetric,
Toeplitz matrix. This proves (a) and (ii) in the lemma.
To prove (b), we use (\ref{eq:rho111}) and observe that the 
right hand side depends only on $|t-\tau|$, thus proving 
that $\rho_{t,t}^{(a,a)|\shortminus \{(a,t),(a,\tau)\}}$ only depends on
$|t-\tau|$.

To prove (iii) we use (\ref{eq:var3}) with $a\neq b$.
From (\ref{eq:var3}) it immediately follows that if $D_{a,a},$ $D_{b,b}$ and
$D_{a,b}$ are Toeplitz, then $\rho_{t,\tau}^{(a,b)}$ only depends on the lag
$(t-\tau)$.  

Conversely, to prove that $D_{a,b}$ is Toeplitz given that for all $t$ and $\tau$; $\rho_{t,\tau}^{(a,a)} =
\rho_{0,t-\tau}^{(a,a)} = \rho_{0,\tau-t}^{(a,a)}$.
$\rho_{t,\tau}^{(b,b)} =
\rho_{0,t-\tau}^{(b,b)} = \rho_{0,\tau-t}^{(b,b)}$, 
$\rho_{t,\tau}^{(a,b)} =
\rho_{0,t-\tau}^{(a,b)}$,
we use the same strategy used to prove (ii). This yields the solution
\begin{eqnarray*}
[D_{a,b}]_{t,\tau} &=& \frac{-1 +
                      \sqrt{1+4(\rho_{0,t-\tau}^{(a,b)})^{2}\sigma_{a}^{-2}\sigma_{b}^{-2}}}{2
                      \rho_{0,t-\tau}^{(a,b)}},
\end{eqnarray*}
which proves that $D_{a,b}$ is Toeplitz.  Thus proving the result. \hfill $\Box$

\subsection{Proof of results in Section \ref{sec:covtimeD}}\label{sec:covtimeDproofs}

To prove the results in Section \ref{sec:covtimeD} we follow 
a similar strategy to the proofs of Section \ref{sec:covD}, but
permute the submatrices $\{C_{a,b}\}$ in $C$ rather than the individual entries in
$C$. The proofs in this section are less technical than those in Section \ref{sec:covD}.

Define the two non-intersecting sets 
$\mathcal{S} = \{\alpha_{1},\ldots,\alpha_{r}\}$ and its complement
$\mathcal{S}^{\prime} = \{\beta_{1},\ldots,\beta_{s}\}$ where 
$\mathcal{S}\cup \mathcal{S}^{\prime} = \{1,2,\ldots,p\}$.
We now obtain an expression for the covariance of  
$\{X_{t}^{(c)};c\in\mathcal{S}\}$ after removing their linear dependence
on $\{X_{s}^{(c)};s\in \mathbb{Z}, c\in\mathcal{S}^{\prime}\}$.
To do so, we define the submatrix $C_{\mathcal{S},\mathcal{S}} =
(C_{a_{1},b_{1}};a_{1},b_{1}\in \mathcal{S})$ where we note that
\begin{eqnarray*}
\var[X_{t}^{(a)};t\in \mathbb{Z},a\in
\mathcal{S}] = C_{\mathcal{S},\mathcal{S}}. 
\end{eqnarray*}
A block permuted version of $C$ with $C_{\mathcal{S},\mathcal{S}}$ in
the top left hand corner is 
\begin{eqnarray*}
C =
\left(
\begin{array}{cccc}
C_{\mathcal{S},\mathcal{S}} & C_{\mathcal{S},\mathcal{S}^{\prime}} \\
C_{\mathcal{S},\mathcal{S}^{\prime}}^{*} & C_{\mathcal{S}^{\prime},\mathcal{S}^{\prime}} \\
\end{array}
\right),
\end{eqnarray*}
where $C_{\mathcal{S},\mathcal{S}^{\prime}} = (C_{a_{1},b_{2}};a_{1}\in \mathcal{S}
\textrm{ and }b_{2}\in \mathcal{S}^{\prime})$, $C_{\mathcal{S}^{\prime},\mathcal{S}^{\prime}} = (C_{a_{2},b_{2}};a_{2}\in \mathcal{S}^{\prime}
\textrm{ and }b_{2}\in \mathcal{S}^{\prime})$. By using
standard results, the conditional variance of $(X_{t}^{(a)};t\in \mathbb{Z},a\in
\mathcal{S})$ given $\{X^{(b)};b\in \mathcal{S}^{\prime}\}$ is the
Schur complement of $C_{\mathcal{S}^{\prime},\mathcal{S}^{\prime}}$ of $C$:
\begin{eqnarray}
\label{eq:project2}
\var\left[X_{t}^{(a)} -  
P_{\mathcal{H}-(X^{(c)};c\in \mathcal{S}^{\prime})}(X_{t}^{(a)});t\in
  \mathbb{Z},a\in \mathcal{S}\right]  = 
C_{\mathcal{S},\mathcal{S}} -
  C_{\mathcal{S},\mathcal{S}^{\prime}}C_{\mathcal{S}^{\prime},\mathcal{S}^{\prime}}^{-1}
C_{\mathcal{S}^{\prime},\mathcal{S}} = P.
\end{eqnarray}
Using the above,  entrywise for all $a,b\in \mathcal{S}$ and
$t,\tau\in \mathbb{Z}$, we have
\begin{eqnarray*}
\cov\left[X_{t}^{(a)} -
  P_{\mathcal{H}-(X^{(c)};c\in \mathcal{S}^{\prime})}(X_{t}^{(a)}), 
X_{\tau}^{(b)} -
  P_{\mathcal{H}-(X^{(c)};c\in \mathcal{S}^{\prime})}(X_{\tau}^{(b)})\right] 
= [P_{a,b}]_{t,\tau}.
\end{eqnarray*}


We now relate the conditional variance $P$ (defined in (\ref{eq:project2})) to the matrix $D$. Using 
the block operator inversion identity in (\ref{eq:inverseidentity}) we have 
\begin{eqnarray}
\label{eq:schur}
D &=& 
\left(
\begin{array}{cccc}
P^{-1} & -P^{-1}C_{\mathcal{S},\mathcal{S}^{\prime}}C_{\mathcal{S}^{\prime},\mathcal{S}^{\prime}}^{-1} \\
-C_{\mathcal{S}^{\prime},\mathcal{S}^{\prime}}^{-1}C_{\mathcal{S}^{\prime},\mathcal{S}}P^{-1}
       &  
(C_{\mathcal{S}^{\prime},\mathcal{S}^{\prime}} -
         C_{\mathcal{S}^{\prime},\mathcal{S}}C_{\mathcal{S},\mathcal{S}}^{-1}
C_{\mathcal{S},\mathcal{S}^{\prime}})^{-1} \\
\end{array}
\right),
\end{eqnarray}
where $P$ is defined in (\ref{eq:project2}).

We use (\ref{eq:schur}) to prove the results in Section \ref{sec:covtimeD}.
\vspace{2mm}

\noindent {\bf PROOF of Theorem \ref{lemma:GGMstGGM}} 
To prove the result we use (\ref{eq:project2}), where we set $\mathcal{S} = \{a\}$ and
$\mathcal{S} = \{a,b\}$.

To prove (i) we let $\mathcal{S} = \{a\}$. By using
  (\ref{eq:project2}) we have 
\begin{eqnarray}
\label{eq:project3}
\var\left[X_{t}^{(a)} -  
P_{\mathcal{H}-X^{(a)}}(X_{t}^{(a)});t\in \mathbb{Z}\right]  = D_{a,a}^{-1}.
\end{eqnarray}
Thus entrywise by definition we have 
$\rho_{t,\tau}^{(a,a)|\shortminus\{a\}}= [D_{a,a}^{-1}]_{t,\tau}$,
this proves (i).

To prove (ii) we let $\mathcal{S} = \{a,b\}$.
By using  (\ref{eq:project2}) and (\ref{eq:inverseidentity}) we have 
\begin{eqnarray*}
&&\var\left[X_{t}^{(a)|\shortminus\{a,b\}};t\in \mathbb{Z}, c\in \{a,b\}\right]
 = 
\left(
\begin{array}{cc}
D_{a,a} & D_{a,b} \\
D_{b,a} & D_{b,b} \\
\end{array}
\right)^{-1} \\
&=&
\left(
\begin{array}{cc}
 (D_{a,a} - D_{a,b}D_{b,b}^{-1}D_{b,a})^{-1} & -(D_{a,a} - D_{a,b}D_{b,b}^{-1}D_{b,a})^{-1}D_{a,b}D_{b,b}^{-1} \\
-D_{b,b}^{-1}D_{b,a}(D_{a,a} - D_{a,b}D_{b,b}^{-1}D_{b,a})^{-1}
       &  
(D_{b,b} - D_{b,a}D_{a,a}^{-1}D_{a,b})^{-1} \\
\end{array}
\right),
\end{eqnarray*}
where the above follows from  (\ref{eq:inverseidentity}). Comparing entries in the above matrix gives
\begin{eqnarray*}
\rho_{t,\tau}^{(a,a)|\shortminus\{a,b\}}&=&  [-(D_{a,a} - D_{a,b}D_{b,b}^{-1}D_{b,a})^{-1}D_{a,b}D_{b,b}^{-1}
                                             ]_{t,\tau} \\
\rho_{t,\tau}^{(a,b)|\shortminus\{a,b\}}&=&[(D_{a,a} -
                                             D_{a,b}D_{b,b}^{-1}D_{b,a})^{-1}]_{t,\tau} \\
\textrm{ and }
 \rho_{t,\tau}^{(b,b)|\shortminus\{a,b\}}&=&  [(D_{b,b} -
                                             D_{b,a}D_{a,a}^{-1}D_{a,b})^{-1}]_{t,\tau}.
\end{eqnarray*}
This proves (ii).
\hfill $\Box$

\vspace{3mm}

\vspace{2mm}
\noindent {\bf PROOF of Theorem \ref{lemma:conditional}} Before we prove
the result, we note the  following invariance properties of
(infinite dimension) Toeplitz operators. If $A$
and $B$ are bounded Toeplitz operators then (a) $AB$ is
Toeplitz (b) if $A$ is Toeplitz and has a bounded inverse, then
$A^{-1}$ is Toeplitz; these results are a consequence of Toeplitz
Theorem, \cite{p:toe-11}.
It is important to mention that these results only hold if the
Toeplitz operators are bi-infinite in the sense the entries of $A$ are $A_{t,\tau} =
A_{t-\tau}$ for all $t,\tau\in \mathbb{Z}$. 
The same results do not hold if the
Toeplitz operators are semi-infinite where $A$ is defined as $A_{t,\tau} =
A_{t-\tau}$ for all $t,\tau\in \mathbb{Z}^{+}$. 

We recall from the proof of Theorem \ref{lemma:GGMstGGM} that 
 \begin{eqnarray}
\label{eq:covAAcond}
\var\left[X_{t}^{(a)|\shortminus\{a\}};t\in \mathbb{Z}, c\in
   \{a\}\right]  = D_{a,a}^{-1}
\end{eqnarray}
and 
\begin{eqnarray}
&&\var\left[X_{t}^{(a)|\shortminus\{a,b\}};t\in \mathbb{Z}, c\in
   \{a,b\}\right] \nonumber\\
&=&
\left(
\begin{array}{cc}
 (D_{a,a} - D_{a,b}D_{b,b}^{-1}D_{b,a})^{-1} & -(D_{a,a} - D_{a,b}D_{b,b}^{-1}D_{b,a})^{-1}D_{a,b}D_{b,b}^{-1} \\
-D_{b,b}^{-1}D_{b,a}(D_{a,a} - D_{a,b}D_{b,b}^{-1}D_{b,a})^{-1}
       &  
(D_{b,b} - D_{b,a}D_{a,a}^{-1}D_{a,b})^{-1} \label{eq:covABcond}\\
\end{array}
\right)
\end{eqnarray}
we use this to prove the result. 

We first prove (i). If $\{X_{t}^{(a)|\shortminus\{a,b\}},X_{t}^{(b)|\shortminus\{a,b\}} \}_{t}$  is conditionally
noncorrelated then $D_{a,b} = D_{b,a}=0$. From (\ref{eq:covABcond}) we
have 
\begin{eqnarray*}
&&\var\left[X_{t}^{(a)|\shortminus\{a,b\}};t\in \mathbb{Z}, c\in
   \{a,b\}\right]  = 
\left(
\begin{array}{cc}
D_{a,a}^{-1} & 0 \\
0 & D_{b,b}^{-1} \\
\end{array}
\right).
\end{eqnarray*}
Thus $\rho_{t,\tau}^{(a,b)|\shortminus\{a,b\}}=\cov[X_{t}^{(a)|\shortminus\{a,b\}},X_{\tau}^{(b)|\shortminus\{a,b\}}]=0$ for all $t$ and
$\tau$. Conversely, if \\
$\rho_{t,\tau}^{(a,b)|\shortminus\{a,b\}}=\cov[X_{t}^{(a)|\shortminus\{a,b\}},X_{\tau}^{(b)|\shortminus\{a,b\}}]=0$
for all $t$ and $\tau$, then using   (\ref{eq:covABcond}) we have $D_{a,b} = 0$.   This proves (i).
 
To prove (ii) we use (\ref{eq:covAAcond}).
If $D_{a,a}$ is Toeplitz, then
$D_{a,a}^{-1}$ is Toeplitz and $\rho_{t,\tau}^{(a,a)|\shortminus\{a\}}=
[D_{a,a}^{-1}]_{t,\tau} =\rho_{0,t-\tau}^{(a,a)|\shortminus\{a\}}$ for
$t$ and $\tau$ (thus $\rho_{t,\tau}^{(a,a)|\shortminus\{a\}}$ is shift
invariant).
Conversely, if for all $t$ and $\tau$, there exists a sequence 
$\{\rho_{r}^{(a,a)|\shortminus\{a\}}\}_{r}$ where
$\rho_{t-\tau}^{(a,a)|\shortminus\{a\}} =\rho_{t,\tau}^{(a,a)|\shortminus\{a\}}$,
then since $\rho_{t,\tau}^{(a,a)|\shortminus\{a\}}=
[D_{a,a}^{-1}]_{t,\tau}$ this implies 
$D_{a,a}^{-1}$ is Toeplitz. Thus $D_{a,a}$ is Toeplitz. This proves (ii).

To prove (iii) we use (\ref{eq:covABcond}).
If $D_{a,a}, D_{a,b}$ and
$D_{b,b}$ are Toeplitz, then 
by the inverse properties of Toeplitz operators (described at the start of the proof)
$(D_{a,a}-D_{a,b}D_{b,b}D_{ba})^{-1}$, $ -(D_{a,a}-D_{a,b}D_{b,b}D_{ba})^{-1}D_{a,b}D_{b,b}^{-1}$
and $(D_{b,b}-D_{b,a}D_{a,a}D_{a,b})^{-1}$
are Toeplitz. Thus 
the conditional covariances
$\{\rho_{t,\tau}^{(a,a)|\shortminus\{a,b\}}\}$,
$\{\rho_{t,\tau}^{(a,b)|\shortminus\{a,b\}}\}$ and
$\{\rho_{t,\tau}^{(a,b)|\shortminus\{a,b\}}\}$ are shift invariant. 
 Conversely, suppose 
\begin{eqnarray*}
\var\left[X_{t}^{(c)|\shortminus\{a,b\}};t\in \mathbb{Z}, c\in \{a,b\}\right]
 &=& \left(
\begin{array}{cc}
E_{a,a} & E_{a,b} \\
E_{a,b}^{*} & E_{b,b} \\
\end{array}
\right)
\end{eqnarray*}
where $E_{a,a}$, $E_{a,b}$ and $E_{b,b}$ are Toeplitz.
Then by using the relation 
\begin{eqnarray*}
\left(
\begin{array}{cc}
E_{a,a} & E_{a,b} \\
E_{a,b}^{*} & E_{b,b} \\
\end{array}
\right)^{-1}=
\left(
\begin{array}{cc}
D_{a,a} & D_{a,b} \\
D_{a,b}^{*} & D_{b,b} \\
\end{array}
\right),
\end{eqnarray*}
and (\ref{eq:inverseidentity}),
we have that $D_{a,a}, D_{a,b}$ and $D_{b,b}$ are Toeplitz. 
This proves (iii). \hfill $\Box$

\vspace{2mm}
\noindent {\bf PROOF of Corollary \ref{cor:nodes}} 
The result
follows immediately from (\ref{eq:project2}) where 
\begin{eqnarray*}
&&\var\left[X_{t}^{(a)} -  
P_{\mathcal{H} - (X^{(c);}c\in
  \mathcal{S}^{\prime})}(X_{t}^{(a)});t\in \mathbb{Z},a\in
  \mathcal{S}\right]  \\
&=&
\left(
\begin{array}{cccc}
D_{\alpha_{1},\alpha_{1}} & D_{\alpha_{1},\alpha_{2}} & \ldots &
                                                                 D_{\alpha_{1},\alpha_{r}} 
  \\
 D_{\alpha_{2},\alpha_{1}} & D_{\alpha_{2},\alpha_{2}} & \ldots &
                                                                 D_{\alpha_{2},\alpha_{r}} \\
\vdots & \vdots & \ddots & \vdots \\
D_{\alpha_{r},\alpha_{1}} & D_{\alpha_{r},\alpha_{2}} & \ldots &
                                                                 D_{\alpha_{r},\alpha_{r}} \\
\end{array}
\right)^{-1}.
\end{eqnarray*}
Thus proving the result. \hfill $\Box$
\vspace{2mm}

%% file: Appendix_Section3.tex
\section{Proofs for Section \ref{sec:fourier}}\label{sec:fourierproof}

\subsection{Proof of results in Section \ref{sec:fourierback}}\label{sec:fourierproofgen}



We start by reviewing some of the relationships between the bounded matrix
operator $A: \ell_{2}\rightarrow \ell_2$ (where $A = (A_{t,\tau}; t,\tau\in \mathbb{Z})$) 
and  the corresponding integral kernel of $F^{*}AF$, which is $\sum_{t\in \mathbb{Z}}\sum_{\tau \in
\mathbb{Z}}A_{t,\tau}e^{it\omega - i\tau \lambda}$. We mention that if
the entries of $A$ were the covariance of a time series and
$\sum_{t,\tau\in \mathbb{Z}}A_{t,\tau}^{2}<\infty$, then
$A(\omega,\lambda) = \sum_{t\in \mathbb{Z}}\sum_{\tau \in
\mathbb{Z}}A_{t,\tau}e^{it\omega - i\tau \lambda}$ (the Loeve
dual-frequency spectrum)  is a well defined
function in $L_{2}[0,2\pi)^{2}$ (see, for example, \cite{p:gor-19} and \cite{p:ast-19}).

The $j$th row of $A$ can be extracted from $A$ using $A^{\prime}u_{j}$, where
$u_{j}\in \ell_{2}$ with $u_{j}  = (\ldots,0,1,0,0,\ldots)$ with $1$
at the $j$th entry. It is clear that $A^{\prime}u_{j} =
(A_{j,\cdot})^{\prime}$ (the $j$th row of $A$) and $\{A^{\prime}u_{j}\}_{j\in
  \mathbb{Z}}$ reproduces all the rows of $A$. 
We now find the parallel to $A^{\prime}u_{j}$ for $F^{*}AF$.
Since $F$ is an isomorphism from $\ell_{2}$ to $L_{2}[0,2\pi)$
the equivalent of $u_{j}$ in $L_{2}[0,2\pi)$ is 
$F^{*}u_{j} = \exp(-ij\omega)$ (inverting back gives
$[F\exp(-ij\cdot)]_{t}=[u_{j}]_{t}$, the $t$th entry in the vector
$u_{j}$). Therefore, if $E
= F^{*}AF$ has integral kernel $A(\omega,\lambda)$, then  
\begin{eqnarray*}
[EF^{*}u_{j}](\lambda) =
  \int_{0}^{\pi}A(\omega,\lambda)\exp(-ij\omega)d\omega = A_{j}(\lambda), 
\end{eqnarray*}
where $A_{j}(\lambda) = \sum_{\tau\in
  \mathbb{Z}}A_{j,\tau}\exp(-i\tau\lambda) 
\in L_{2}[0,2\pi)$ and forms the building blocks of
$A(\omega,\lambda)$ (since $A(\omega,\lambda) = \sum_{t\in \mathbb{Z}}A_{t}(\lambda)\exp(it\omega)$). 
 $(FEF^{*})u_{j}$ yields
the $j$th row of the infinite dimensional matrix
$(FEF^{*})$ and the $(j,s)th$ entry of $A=(FEF^{*})$ is
\begin{eqnarray}
\label{eq:FEFstar}
[(FEF^{*})u_{j}]_{s} = \frac{1}{(2\pi)^{2}}\int_{0}^{2\pi}\int_{0}^{2\pi}
  A(\omega,\lambda)\exp(-ij\omega)\exp(is\lambda)d\omega d\lambda.
\end{eqnarray}
The above gives the relationship between $A(\omega,\lambda)$ and $A$.

The proof of Lemma \ref{lemma:diagonal} 
follows from \cite{p:toe-11}
(see \cite{b:bot-00}, Theorem 1.1). However, for completeness and to
explicitly connect the result to $A(\omega,\lambda)$ we give a proof
below (it is based on the discussion above).

\vspace{2mm}

\noindent {\bf PROOF of Lemma \ref{lemma:diagonal}} 
We  prove that the infinite dimensional
Toeplitz matrix $A$ leads to a
diagonal kernel of the form $\delta_{\omega,\lambda}A(\omega)$. Suppose that $A$ is a bounded operator that is a
Toeplitz matrix with entries $\{a_{j}\}_{j}$. Then the integral kernel is
\begin{eqnarray*}
A(\omega,\lambda) &=& \sum_{t\in \mathbb{Z}}\sum_{\tau\in
                      \mathbb{Z}}a_{t-\tau}\exp(it\omega-i\tau\lambda)
  \\
 &=&  \sum_{\tau\in \mathbb{Z}}\exp(-i\tau(\lambda-\omega))\sum_{r\in
     \mathbb{Z}}a_{r}\exp(ir\omega) = \delta_{\omega,\lambda}A(\omega)
\end{eqnarray*}
where $A(\omega)  = \sum_{r\in \mathbb{Z}}a_{r}\exp(ir\omega)$. Since
$\|A\|<\infty$, defining the infinite sequence $v=\{v_{j}\}$ where
$v_{j}=0$ for all $j\neq 0$ and $v_{0}=1$ we have
$\sum_{j\in \mathbb{Z}}a_{j}^{2} =\|Av\|_{2}\leq \|A\|\|v\|_{2}\leq
\|A\|$, thus $A(\cdot)\in L_{2}[0,2\pi)$. 

We now use (\ref{eq:FEFstar}) to prove the converse. 
Substituting $A(\omega,\lambda) =
\delta_{\omega,\lambda}A(\omega)$ into (\ref{eq:FEFstar}) gives 
\begin{eqnarray*}
[(FEF^{*})u_{j}]_{s} &=& \frac{1}{(2\pi)^{2}}\int_{0}^{2\pi}\int_{0}^{2\pi}A(\omega)\delta_{\omega,\lambda}
 \exp(-ij\omega)\exp(is\lambda)d\omega d\lambda  = a_{j-s}
\end{eqnarray*}
where
$a_{r}=(2\pi)^{-1}\int_{0}^{2\pi}A(\omega)\exp(-ir\omega)d\omega$. Thus the
$j$th column of $FEF^{*}$ is $\{a_{s-j}\}_{s\in \mathbb{Z}}$, which
proves that the matrix defined by $FEF^{*}$ is Toeplitz. 
\hfill $\Box$

\vspace{2mm}

\noindent {\bf PROOF of Lemma \ref{lemma:inversegeneral}}  
Since $A$ is block Toeplitz it follows from Lemma
\ref{lemma:diagonal} that ${\bf A}(\omega,\lambda) = {\bf
  A}(\omega)\delta_{\omega,\lambda}$. 

To derive an expression for the inverse, we 
first consider the case that $d=1$. By
definition $AA^{-1} = I$ (where $I$ denotes the infinite dimension
identity matrix), thus $ F^{*}F =  (F^{*}AF)(F^{*}A^{-1}F)$. By Lemma \ref{lemma:diagonal},
the kernel operator of $F^{*}AF$ is $A(\omega)\delta_{\omega,\lambda}$
and the kernel operator of $F^{*}A^{-1}F$ (since $A^{-1}$ is Toeplitz) is
$B(\omega)\delta_{\omega,\lambda}$. Since for all $g\in L_{2}[0,2\pi)$
we have
\begin{eqnarray*}
g(\omega) &=& [F^{*}F(g)](\omega) =
  \frac{1}{(2\pi)^{2}}\int_{0}^{2\pi}B(\omega)\delta_{\omega,u}\int_{0}^{2\pi}A(u)\delta_{u,\lambda}g(\lambda)
  d\lambda du \\
 &=&
     \frac{1}{(2\pi)}\int_{0}^{2\pi}B(\omega)A(u)g(u)\delta_{\omega,u}du
     = B(\omega)A(\omega)g(\omega),
\end{eqnarray*}
then $B(\omega) = A(\omega)^{-1}$. This proves the result for all
$d=1$. 

The proof for  $d>1$ uses the following invariance properties. If $A$ and
$B$ are bounded Toeplitz matrix operators with kernels $A(\omega)\delta_{\omega,\lambda}$ and 
$B(\omega)\delta_{\omega,\lambda}$ respectively, then $A+B$ and $AB$
are Toeplitz with kernels
$[A(\omega)+B(\omega)]\delta_{\omega,\lambda}$ and 
$A(\omega)B(\omega)\delta_{\omega,\lambda}$ respectively.
Using these 
properties together with the block operator inversion identity
(in (\ref{eq:inverseidentity})) we will show, below, that the Lemma \ref{lemma:inversegeneral}
holds for $d\geq 2$.  We focus on $d=2$ (the proof for
$d>2$ follows by induction). Let 
\begin{eqnarray*}
G=\left(
\begin{array}{cc}
A & B \\
B^{*} & C
\end{array}
\right)
\end{eqnarray*}
where $G$ is a bounded operator and $A,B$ and $C$ are Toeplitz operators on
$\ell_{2}$, with integral kernels $A(\omega)\delta_{\omega,\lambda}$,
$B(\omega)\delta_{\omega,\lambda}$  and $C(\omega)\delta_{\omega,\lambda}$.
Then by (\ref{eq:inverseidentity}) 
\begin{eqnarray*}
FG^{-1}F^{*}=
\left(
\begin{array}{cc}
FP^{-1}F^{*} & -FP^{-1}BC^{-1}F^{*} \\
-FC^{-1}B^{*}P^{-1}F^{*} & F(C-B^{*}A^{-1}B)^{-1}F^{*} \\
\end{array}
\right),
\end{eqnarray*}
where $P = A-BC^{-1}B^{*}$. By the Toeplitz invariance properties described
above, the integral kernel of $FPF^{*}$ is $P(\omega)\delta_{\omega,\lambda}$
where
\begin{eqnarray*}
P(\omega) = [A(\omega) - |B(\omega)|^{2}C(\omega)]. 
\end{eqnarray*}
Thus by the proof for $d=1$, the integral kernel of $FP^{-1}F^{*}$
(the top left hand side of $FG^{-1}F^{*}$) is
 $P(\omega)^{-1}\delta_{\omega,\lambda}$. 
A similar result holds for the other entries in $FG^{-1}F^{*}$. Therefore,  the  integral kernel of $FG^{-1}F^{*}$ is 
\begin{eqnarray*}
&&\left(
\begin{array}{cc}
P(\omega)^{-1} & P(\omega)^{-1}B(\omega)C(\omega)^{-1}\\
C(\omega)^{-1}B(\omega)^*P(\omega)^{-1} & (C(\omega)-|B(\omega)|^{2}A(\omega)^{-1})^{-1} \\
\end{array}
\right)\delta_{\omega,\lambda} \\
&=& 
\left(
\begin{array}{cc}
A(\omega) & B(\omega) \\
B(\omega)^* & C(\omega) \\
\end{array}
\right)^{-1}\delta_{\omega,\lambda}.
\end{eqnarray*}
This proves the result for $d=2$. By induction the result
can be proved for $d>2$.
\hfill $\Box$

\subsection{Proof of results in Sections \ref{sec:precisionfourier}
  and \ref{sec:coherence}}\label{appendix:preccoh}

\noindent {\bf PROOF of Theorem \ref{lemma:spectral}} Under Assumption
\ref{assum:lambda} and  by using Lemma \ref{lemma:inverse}
for all  $1\leq a,b\leq p$, $D_{a,b}$ are bounded operators. Thus the proof is a
straightforward application of Lemma \ref{lemma:diagonal}. We
summarize the main points below.

To prove (i) we note that 
$D_{a,b}=0$ is a special case of Toeplitz matrix, thus  $F^{*}D_{a,b}F
= 0\cdot\delta_{\omega,\lambda}=0$. Conversely, if $F^{*}D_{a,b}F = 0$, then $D_{a,b}=0$.

The proof of (ii) and (iii) immediately follow from Definition
\ref{def:network} and Lemma \ref{lemma:diagonal}. \hfill $\Box$

\vspace{3mm}
We now prove the results in Section \ref{sec:coherence}. We first
consider the Fourier transform of the rows of $D_{a,b}$ and $D_{a,b}$ in the case
a node or edge is conditionally stationary.
\vspace{2mm}

\noindent {\bf PROOF of Theorem \ref{lemma:spectral2}}
We first prove (i). If the node $a$ is conditionally stationary then 
 $D_{a,a}$ is Toeplitz and its entries are 
determined by the row $\{[D_{a,a}]_{0,r}\}_{r}$. By using Lemma
\ref{lemma:diagonal} we have $\Gamma^{(a,a)}_{t}(\omega)=\Gamma^{(a,a)}(\omega) = 
\sum_{r=-\infty}^{\infty}[D_{a,a}]_{0,r}\exp(ir\omega)$.  To
understand the meaning of this quantity, we note that from Proposition
\ref{lemma:networkpartial}  for all $t\neq \tau$
\begin{eqnarray*} 
\phi_{0,t-\tau}^{(a,a)} = \frac{\rho_{0,t-\tau}^{(a,a)}}{
\sqrt{\rho_{0,0}^{(a,a)|\shortminus
  \{(a,0),(a,t-\tau)\}}\rho_{0,0}^{(a,a)|\shortminus
  \{(a,0),(a,t-\tau)\}} }} = -\frac{[D_{a,a}]_{t,\tau}}{\sqrt{[D_{a,a}]_{t,t}[D_{a,a}]_{\tau,\tau}}},
\end{eqnarray*}
where $\rho_{0,t-\tau}^{(a,a)}$ and $\rho_{0,0}^{(a,a)|\shortminus
  \{(a,0),(a,t-\tau)\}}$ is defined in  (\ref{eq:partial1}) and 
 (\ref{eq:rhottconditional}) respectively. Thus, we have 
 $[D_{a,a}]_{t,\tau} =
\sqrt{[D_{a,a}]_{t,t}[D_{a,a}]_{\tau,\tau}}\phi_{0,t-\tau}^{(a,a)}
$. Further, we know that
$[D_{a,a}]_{t,t}= 1/\sigma^{2}_{a} = 1/\rho_{0,0}^{(a,a)}$. Together
this gives 
\begin{eqnarray*}
\Gamma^{(a,a)}(\omega) = 
\frac{1}{\rho_{0,0}^{(a,a)}}\left[ 1-
  \sum_{r\in \mathbb{Z}\backslash\{0\}}\phi_{0,r}^{(a,a)}\exp(ir\omega)\right].
\end{eqnarray*}
This proves (i). 

The proof of (ii) is identical to (i), thus we omit the
details. \hfill $\Box$
\vspace{2mm}

\noindent {\bf PROOF of Theorem \ref{lemma:spectral3}}
We first prove (i). By using Lemma
\ref{lemma:inverse}, the integral kernel of $F^{*}D_{a,a}^{-1}F$ is 
$\delta_{\omega,\lambda}[\Gamma^{(a,a)}(\omega)]^{-1}$. We recall that
$D_{a,a}^{-1}$ contains the time series partial covariances and by conditional
stationarity and Theorem \ref{lemma:GGMstGGM} we have 
$[D_{a,a}^{-1}]_{t,\tau} = \rho_{0,t-\tau}^{(a,a)|\shortminus\{a\}}$. Using
this it is easily seen that the partial spectrum for the nodal time series partial covariance
is 
\begin{eqnarray*}
\Gamma^{(a,a)}(\omega)^{-1} = \sum_{r\in \mathbb{Z}}\cov[X_{0}^{(a,a)|\shortminus\{a\}}, X_{r}^{(a,a)|\shortminus\{a\}}]\exp(ir\omega)=
\sum_{r\in \mathbb{Z}}\rho_{r}^{(a,a)|\shortminus\{a\}}\exp(ir\omega).
\end{eqnarray*}
To prove (ii) we use that  $(a,b)$ is a conditionally stationary
edge and define the suboperator block Toeplitz matrix
$D_{\{a,b\}}:\ell_{2,2}\rightarrow \ell_{2,2}$, where  $D_{\{a,b\}} =
(D_{e,f};e,f\in \{a,b\})$. The integral kernel of $F^{*}D_{\{a,b\}}F$ is
$\Gamma_{\{a,b\}}(\omega)\delta_{\omega,\lambda}$ where 
\begin{eqnarray}
\Gamma_{\{a,b\}}(\omega)&=&\left(
\begin{array}{cc}
 \Gamma^{(a,a)}(\omega)& \Gamma^{(a,b)}(\omega) \\
\Gamma^{(a,b)}(\omega)^{*} & \Gamma^{(b,b)}(\omega) \\
\end{array}
\right).
\end{eqnarray}
The time series partial covariances are contained within the inverse
$[D_{\{a,b\}}]^{-1}$ (which is block Toeplitz) (see Theorem
\ref{lemma:GGMstGGM}). By using 
Lemma \ref{lemma:inversegeneral} the
kernel of $F^{*}[D_{\{a,b\}}]^{-1}F$ is
$\Gamma_{\{a,b\}}(\omega)^{-1}\delta_{\omega,\lambda}$.  
Using this together with equation (\ref{eq:covABcond}) we have
\begin{eqnarray*}
\Gamma_{\{a,b\}}(\omega)^{-1} 
&=& 
\sum_{r\in \mathbb{Z}}
\left(
\begin{array}{cc}
\rho^{(a,a)|\shortminus\{a,b\}}_{0,r} & \rho^{(a,b)|\shortminus\{a,b\}}_{0,r}\\
\rho^{(b,a)|\shortminus\{a,b\}}_{0,r} &  \rho^{(b,b)|\shortminus\{a,b\}}_{0,r} \\
\end{array}
\right)
\exp(ir\omega) \\
&=& 
\frac{1}{\det[\Gamma_{\{a,b\}}(\omega)]}\left(
\begin{array}{cc}
\Gamma^{(b,b)}(\omega) & -\Gamma^{(a,b)}(\omega) \\
-\Gamma^{(a,b)}(\omega)^{*} & \Gamma^{(a,a)}(\omega) \\
\end{array}
\right).
\end{eqnarray*}
This proves the result. 
\hfill $\Box$

\subsection{Proof of the results in Section \ref{sec:noderegression}}\label{appendix:regression}

\noindent {\bf PROOF of Proposition \ref{lemma:expectation}}
To connect the regression coefficients to entries in $D$ we use
the identity in (\ref{eq:schur}) where we set $\mathcal{S}=\{a\}$ and
$\mathcal{S}^{\prime}=\{1,\ldots,p\}\backslash\{a\}$. This gives 
\begin{eqnarray}
\label{eq:TD}
D &=& 
\left(
\begin{array}{cccc}
D_{a,a} & -D_{a,a}H_{a}^{\prime}G_{a}^{-1} \\
-G^{-1}_{a}H_{a}D_{a,a} &  (G_{a} - H_{a}^{*}C_{a,a}^{-1}H_{a})^{-1} \\
\end{array}
\right)
\end{eqnarray}
where $H_{a} = C_{a,\mathcal{S}^{\prime}}$ and $G_{a} =
C_{\mathcal{S}^{\prime},\mathcal{S}^{\prime}}$.
We recall that $D = (D_{e,f}, e,f\in \{1,\ldots,p\})$. Therefore,
comparing the blocks on the left and right hand side of (\ref{eq:TD})
gives the block vector
\begin{eqnarray}
\label{eq:DAB20}
(D_{a,b};b\neq a) =   -D_{a,a}H_{a}^{*}G_{a}^{-1}.
\end{eqnarray} 
Furthermore, by comparision, it is clear that the prediction
coefficients $B_{b\shortarrow a}$  satisfy 
\begin{eqnarray}
\label{eq:DAB30}
H_{a}^{*}G_{a}^{-1} = \left(B_{b\shortarrow a};b\neq a \right).
\end{eqnarray} 
Comparing (\ref{eq:DAB20}) and (\ref{eq:DAB30})
 for $b\neq a$ we have $-D_{a,a}B_{b\shortarrow a} = D_{a,b}$. Using that
$D_{a,a}$ has an inverse yields the identity
\begin{eqnarray*}
B_{b\shortarrow a} = -D_{a,a}^{-1}D_{a,b}.
\end{eqnarray*} 
This gives the result. \hfill $\Box$

\vspace{3mm}

\noindent {\bf PROOF of Theorem \ref{lemma:regressiontoeplitz}} To prove
(i) and (ii)
we use Lemma \ref{lemma:expectation} where 
\begin{eqnarray*}
B_{b\shortarrow a} = -D_{a,a}^{-1}D_{a,b}.
\end{eqnarray*} 
To prove (i) we note that under Assumption \ref{assum:lambda}, the
null space of  $D_{a,a}^{-1}$ is $0$. Therefore, $D_{a,b}=0$ iff 
$B_{b\shortarrow a} = 0$. This proves (i).

To prove (ii) we note by the invariance properties of infinite
 Toeplitz matrix operators if $D_{a,b}$ and $B_{b\shortarrow a} $
are Toeplitz, then $D_{a,b} = -D_{a,b}B_{b\shortarrow a} $ is
Toeplitz. Conversely, if $D_{a,b}$ and $D_{b,a} $
are Toeplitz, then $B_{b\shortarrow a} = -D_{a,a}^{-1}D_{a,b}$ is
Toeplitz. This proves (ii). \hfill $\Box$

\vspace{3mm}

\noindent {\bf PROOF of Corollary \ref{cor:regfourier}}
To prove (i) we use that  (a) under Assumption \ref{assum:lambda} that $\|D_{a,b}^{-1}\|\leq
\lambda_{\sup}$ and (b)  from Lemma \ref{lemma:inverse},  $D_{a,b}$ is a bounded
operator. Thus, since $B_{b\shortarrow a} =
D_{a,b}^{-1}D_{a,b}$, we have   $\|B_{b\shortarrow a}\|\leq
\|D_{a,b}^{-1}\|\|D_{a,b}\|<\infty$, thus proving (i).

The proofs of (ii) and (iii) are 
similar to the proof of Theorem \ref{lemma:spectral}, thus we omit the
details. \hfill $\Box$

%% file: Appendix_Section4.tex
\section{Proof of Section \ref{sec:finitelocal}}

\subsection{Proof of results in Section \ref{sec:finite}}\label{sec:prooffinite}

We break the proof of Theorem \ref{lemma:meyer3a} into a few
steps. To bound the difference between the rows of $\widetilde{D}_{n}=C_{n}^{-1}$
and $D_{n}$ (the submatrix of $D$), we use that the entries of $\widetilde{D}_{n}$
and $D_{n}$ are the entries of coefficients in a regression. This
allows us to use the Baxter inequality methods developed in
\cite{p:kre-17} to bound the difference between projections on finite
dimensional spaces and infinite dimensional spaces. 
The infinite and finite dimensional spaces we will use are 
$\mathcal{H} = \overline{\textrm{sp}}(X_{t}^{(c)};t\in \mathbb{Z},1\leq
c\leq p)$ and
$\mathcal{H}_{n} = \overline{\textrm{sp}}(X_{\tau}^{(b)};1\leq \tau\leq n,1\leq b \leq  p)$. 

We recall from (\ref{eq:projectXat}) that
\begin{eqnarray}
\label{eq:PX}
P_{\mathcal{H}-X_{t}^{(a)}}(X_{t}^{(a)}) = \sum_{\tau\in \mathbb{Z}}\sum_{b=1}^{p}\beta_{(\tau,b)\shortarrow (t,a)}X_{\tau}^{(b)}
\end{eqnarray}
with $\beta_{(t,a)\shortarrow (t,a)}=0$. Similarly, projecting
$X_{t}^{(a)}$ onto the finite dimensional space $\mathcal{H}_{n}$ is 
\begin{eqnarray}
\label{eq:PX1}
P_{\mathcal{H}_{n}-X_{t}^{(a)}}(X_{t}^{(a)}) = \sum_{\tau = 1}^{n}\sum_{b=1}^{p}
\theta_{(\tau,b)\shortarrow (t,a),n}X_{\tau}^{(b)}, 
\end{eqnarray}
with $\theta_{(t,a)\shortarrow (t,a),n}=0$. Let 
\begin{eqnarray}
\sigma_{a,t}^{2} &=&
\Ex[X_{t}^{(a)} - P_{\mathcal{H}-X_{t}^{(a)}}(X_{t}^{(a)})]^{2} \nonumber\\
\textrm{ and }
\widetilde{\sigma}_{a,t,n}^{2} &=& \Ex[X_{t}^{(a)} - P_{\mathcal{H}_n-X_{t}^{(a)}}(X_{t}^{(a)})]^{2}.  \label{eq:sigmaphiat}
\end{eqnarray}
Define the $(np-1)$-dimensional
vectors  
\begin{eqnarray}
\label{eq:BD}
\underline{B}^{(a,t)}_{n} &=& \{\beta_{(\tau,b)\shortarrow(t,a)};1\leq
  \tau \leq n \textrm{ and }1\leq b\leq p, \textrm{ not }(\tau,b)=(t,a)\} \nonumber\\
\textrm{ and } 
\underline{\Theta}^{(a,t)}_{n} &=& \{\theta_{(\tau,b)\shortarrow(t,a),n};1\leq
  \tau \leq n \textrm{ and }1\leq b\leq p,\textrm{ not }(\tau,b)=(t,a)\}.
\end{eqnarray}
To minimise notation we will drop the $n$, and let 
$\theta_{(\tau,b)\shortarrow (t,a)}=\theta_{(\tau,b)\shortarrow
  (t,a),n}$ and $\widetilde{\sigma}_{a,t}^{2}
=\widetilde{\sigma}_{a,t,n}^{2}$. But we should keep in mind that both
$\theta$ and $\sigma$ depend on $n$. Since the coefficients of a precision matrix are
closely related to the coefficients in a regression it is clear that
the $t$th ``row'' of the matrix $\widetilde{D}_{n}=C_{n}^{-1}$ at node $a$ (which is the
$((a-1)n+t)$th row of $\widetilde{D}_{n}$) is the  rearranged vector
\begin{eqnarray}
\widetilde{\underline{\Theta}}^{(a,t)}_{n} =
  \frac{1}{\widetilde{\sigma}_{a,t}^{2}}[1,-\underline{\Theta}^{(a,t)}_{n}]. \label{eq:Thetatilde}
\end{eqnarray}
The $t$th row of matrix $D_{n}$ at node $a$ is the similarly
rearranged vector
\begin{eqnarray*}
\widetilde{\underline{B}}^{(a,t)}_{n} = \frac{1}{\sigma_{a,t}^{2}}[1,-\underline{B}^{(a,t)}_{n}].
\end{eqnarray*}
Thus the difference between
$\widetilde{\underline{\Theta}}^{(a,t)}_{n}$ and
$\widetilde{\underline{B}}^{(a,t)}_{n}$ is 
\begin{eqnarray*}
&&\widetilde{\underline{\Theta}}^{(a,t)}_{n}
  -\widetilde{\underline{B}}^{(a,t)}_{n} \\
&=& \left[\frac{1}{\widetilde{\sigma}_{a,t}^{2}} -
  \frac{1}{\sigma_{a,t}^{2}}\right]\left(1,- \underline{\Theta}^{(a,t)}_{n}\right)
+\frac{1}{\sigma_{a,t}^{2}}\left[0,\left(\underline{B}^{(a,t)}_{n}-\underline{\Theta}^{(a,t)}_{n}\right)\right].
\end{eqnarray*}
Since both $\widetilde{\underline{\Theta}}^{(a,t)}_{n}$ and
$\widetilde{\underline{B}}^{(a,t)}_{n}$ are the (same) rearranged rows of
$[\widetilde{D}_{n}]_{(a-1)n+t,\cdot}$ and $[D_{n}]_{(a-1)n+t,\cdot}$,
The $\ell_{1}$-difference between the $((a-1)n+t)$th row of
$D_{n}$ and $\widetilde{D}_{n}$ is 
\begin{eqnarray}
\label{eq:rowdiffs}
&& \left\| [D_{n}]_{(a-1)n+t,\cdot} -  [\widetilde{D}_{n}]_{(a-1)n+t,\cdot}\right\|_{1} 
=\left\|\widetilde{\underline{\Theta}}^{(a,t)}_{n}
  -\widetilde{\underline{B}}^{(a,t)}_{n} \right\|_{1} \nonumber\\
&\leq& \frac{|\sigma_{a,t}^{2}-\widetilde{\sigma}_{a,t}^{2}|}{\widetilde{\sigma}_{a,t}^{2}\sigma_{a,t}^{2}} 
\left(1+ \|\underline{\Theta}^{(a,t)}_{n}\|_{1}\right)
+\frac{1}{\sigma_{a,t}^{2}}\left\|\underline{B}^{(a,t)}_{n}-\underline{\Theta}^{(a,t)}_{n}\right\|_{1}.
\end{eqnarray}
In the two lemmas below we obtain a bound for the differences $|\sigma_{a,t}^{2}-\widetilde{\sigma}_{a,t}^{2}|$
and
$\|\underline{B}^{(a,t)}_{n}-\underline{\Theta}^{(a,t)}_{n}\|_{1}$. These
two bounds will prove Theorem \ref{lemma:meyer3a}.

\begin{lemma}\label{lemma:meyer}
Suppose Assumptions \ref{assum:lambda}  and \ref{assum:invcovarianceK} hold.
Let $\underline{B}^{(a,t)}_{n}$ and $\underline{\Theta}^{(a,t)}_{n}$ be
defined as in (\ref{eq:BD}). 
Then 
\begin{eqnarray*}
\left\| \underline{B}^{(a,t)}_{n}-\underline{\Theta}^{(a,t)}_{n} \right\|_{2}
&\leq& \lambda_{\inf}^{-1}\lambda_{\sup}\sum_{b=1}^{p}\sum_{\tau\notin
       \{1,\ldots,n\}}|\beta_{(\tau,b)\shortarrow(t,a)}|.
\end{eqnarray*}
\end{lemma}
{\bf PROOF}  The proof is based on the innovative technique developed in \cite{p:kre-17} (who
used the method to obtain Baxter bounds for stationary spatial processes). 
We start by deriving the normal equations corresponding to
(\ref{eq:PX}) and (\ref{eq:PX1}) for $1\leq s\leq
n$ and $c=1,\ldots,p$ (excluding $(c,s)=(a,t)$). For equation (\ref{eq:PX}) this
gives the normal equations
\begin{eqnarray}
\label{eq:PXnormal}
\cov(X_{t}^{(a)},X_{s}^{(c)}) &=&
                                     \sum_{b=1}^{p}\sum_{\tau=1}^{n}\beta_{(\tau,b)\shortarrow
                                     (t,a)}\cov(X_{\tau}^{(b)}, X_{s}^{(c)}) +\nonumber\\
&&\sum_{b=1}^{p}\sum_{\tau\notin \{1,\ldots,n\}}\beta_{(\tau,b)\shortarrow (t,a)}
\cov(X_{\tau}^{(b)}, X_{s}^{(c)}) 
\end{eqnarray}
and for (\ref{eq:PX1}) this gives 
\begin{eqnarray}
\label{eq:PX1normal}
\cov(X_{t}^{(a)},X_{s}^{(c)}) &=&  \sum_{b=1}^{p}\sum_{\tau=1}^{n}\theta_{(\tau,b)\shortarrow (t,a)}
\cov(X_{\tau}^{(b)}, X_{s}^{(c)}).
\end{eqnarray}
Taking the difference between (\ref{eq:PXnormal}) and
(\ref{eq:PX1normal}) we have
\begin{eqnarray*}
\sum_{b=1}^{p}\sum_{\tau=1}^{n}\left[\beta_{(\tau,b)\shortarrow
                                     (t,a)}-\theta_{(\tau,b)\shortarrow (t,a)}
\right]\cov(X_{\tau}^{(b)}, X_{s}^{(c}) 
&=&-\sum_{b=1}^{p}\sum_{\tau\notin \{1,\ldots,n\}}
\beta_{(\tau,b)\shortarrow (t,a)}
\cov(X_{\tau}^{(b)}, X_{s}^{(c)}). 
\end{eqnarray*}
As the above holds for all $1\leq s \leq n$ and $1\leq c\leq p$
(excluding $X_{t}^{(a)}$) we can write the above as a vector equation
\begin{eqnarray}
&&\sum_{b=1}^{p}\sum_{\tau\in \{1,\ldots,n\}}\left[\beta_{(\tau,b)\shortarrow
                                     (t,a)}-\theta_{(\tau,b)\shortarrow (t,a)}\right]\cov(X_{\tau}^{(b)},\underline{Y}_{n}) \nonumber\\
&=&-\sum_{b=1}^{p}\sum_{\tau\notin \{1,\ldots,n\}}\beta_{(\tau,b)\shortarrow (t,a)}\cov(X_{\tau}^{(b)}, \underline{Y}_{n}), \label{eq:vectoreq}
\end{eqnarray}
where $\underline{Y}_{n} = (X_{s}^{(c)};1\leq s\leq n,1\leq
c\leq p, (c,s)\neq (a,t)\})$. We observe that the LHS of the above
can be expressed as 
\begin{eqnarray}
&&\sum_{b=1}^{p}\sum_{\tau\in \{1,\ldots,n\}}\left[\beta_{(\tau,b)\shortarrow
                                     (t,a)}-\theta_{(\tau,b)\shortarrow
   (t,a)}\right]\cov(X_{\tau}^{(b)},\underline{Y}_{n}) \nonumber\\
&=&
\cov\left( \sum_{b=1}^{p}\sum_{\tau\in \{1,\ldots,n\}}\left[\beta_{(\tau,b)\shortarrow
                                     (t,a)}-\theta_{(\tau,b)\shortarrow
   (t,a)}\right]X_{\tau}^{(b)},\underline{Y}_{n}\right) \nonumber\\
&=&
    \cov\left(\left[\underline{B}^{(a,t)}_{n}-\underline{\Theta}^{(a,t)}_{n} \right]^{\prime}\underline{Y}_{n}^{},\underline{Y}_{n}\right), \label{eq:cov100}
\end{eqnarray}
where the last line of the above is due to
\begin{eqnarray*}
\sum_{b=1}^{p}\sum_{\tau\in \{1,\ldots,n\}}\left[\beta_{(\tau,b)\shortarrow
                                     (t,a)}-\theta_{(\tau,b)\shortarrow
  (t,a)}\right]X_{\tau}^{(b)} = 
\left[\underline{B}^{(a,t)}_{n}-\underline{\Theta}^{(a,t)}_{n} \right]^{\prime}\underline{Y}_{n}^{}.
\end{eqnarray*}
Substituting (\ref{eq:cov100}) into the LHS of (\ref{eq:vectoreq}) gives the vector equation
\begin{eqnarray*}
\left(\var[\underline{Y}_{n}]\right)[\underline{B}^{(a,t)}_{n}-\underline{\Theta}^{(a,t)}_{n}]
=-\sum_{b=1}^{p}\sum_{\tau\notin \{1,\ldots,n\}} \beta_{(\tau,b)\shortarrow
                                     (t,a)}\cov(X_{\tau}^{(b)}, \underline{Y}_{n}^{}).
\end{eqnarray*}
Therefore 
\begin{eqnarray*}
\left[\underline{B}^{(a,t)}_{n}-\underline{\Theta}^{(a,t)}_{n} \right]
=-\var[\underline{Y}_{n}^{}]^{-1}\sum_{b=1}^{p}\sum_{\tau\notin
  \{1,\ldots,n\}}\beta_{(\tau,b)\shortarrow(t,a)}\cov(X_{\tau}^{(b)}, \underline{Y}_{n}).
\end{eqnarray*}
Now taking the $\ell_{2}$-norm of the above we have 
\begin{eqnarray*}
\left\| \underline{B}^{(a,t)}_{n}-\underline{\Theta}^{(a,t)}_{n} \right\|_{2}
&\leq&\left\|\var[\underline{Y}_{n}^{}]^{-1}\right\|\sum_{b=1}^{p}\sum_{\tau\notin  \{1,\ldots,n\}}
|\beta_{(\tau,b)\shortarrow(t,a)}|\left\|\cov(X_{\tau}^{(b)},
       \underline{Y}_{n}^{})\right\|_{2} \\
&\leq&\left\|\var[\underline{Y}_{n}^{}]^{-1}\right\|\sum_{b=1}^{p}
\sum_{\tau\notin \{1,\ldots,n\}}|\beta_{(\tau,b)\shortarrow(t,a)}|\left\|\cov(X_{\tau}^{(b)},
       \underline{Y}_{n}^{})\right\|_{2} \\
&\leq&
       \left\|\var[\underline{Y}_{n}]^{-1}\right\|\left(\sup_{\tau,b}[\sum_{c=1}^{p}\sum_{s=-\infty}^{\infty}
\cov(X_{\tau}^{(b)},X_{s}^{(c)})^{2}]^{1/2}\right)
\sum_{b=1}^{p}\sum_{\tau\notin
       \{1,\ldots,n\}}|\beta_{(\tau,b)\shortarrow(t,a)}|,
\end{eqnarray*}
where $\|\cdot\|_{}$ denotes the (spectral) matrix norm. By
using Assumption \ref{assum:lambda} we have 
$\lambda_{\min}(\var[\underline{Y}_{n}])\geq \lambda_{\inf}$, thus
$\left\|\var[\underline{Y}_{n}^{}]^{-1}\right\|\leq
\lambda_{\inf}^{-1}$. Again by Assumption \ref{assum:lambda} we have 
\\*
$\sup_{\tau,b}[\sum_{c=1}^{p}\sum_{s=-\infty}^{\infty}\cov(X_{\tau}^{(b)},X_{s}^{(c)})^{2}]^{1/2}\leq
\lambda_{\sup}$. Substituting these two bounds into the above, gives 
\begin{eqnarray*}
\left\| \underline{B}^{(a,t)}_{n}-\underline{\Theta}^{(a,t)}_{n} \right\|_{2}
&\leq& \lambda_{\inf}^{-1}\lambda_{\sup}\sum_{b=1}^{p}\sum_{\tau\notin
       \{1,\ldots,n\}}|\beta_{(\tau,b)\shortarrow(t,a)}|.
\end{eqnarray*}
This proves the result. \hfill $\Box$

\vspace{2mm}
\noindent The above gives a bound for the $\ell_{2}$-norm. To obtain a bound on the
$\ell_{1}$-norm we use the Cauchy-Schwarz inequality 
to give
\begin{eqnarray}
\label{eq:BTheta}
\left\|\underline{B}^{(a,t)}_{n}-\underline{\Theta}^{(a,t)}_{n}
   \right\|_{1}  &\leq&(np)^{1/2}\lambda_{\inf}^{-1}\lambda_{\sup}\sum_{b=1}^{p}\sum_{\tau\notin
       \{1,\ldots,n\}}|\beta_{(\tau,b)\shortarrow(t,a)}|.
\end{eqnarray}
Next we bound the difference between $ \widetilde{\sigma}_{a,t}^{2}$  and
$\sigma_{a,t}^{2}$. 

\begin{lemma}\label{lemma:meyer2}
Suppose Assumptions \ref{assum:lambda}  and \ref{assum:invcovarianceK} hold.
Let $\sigma_{a,t}^{2}$ and
$\widetilde{\sigma}_{a,t,n}^{2}$ be defined as in (\ref{eq:sigmaphiat}). 
Then 
\begin{eqnarray*}
0\leq  \widetilde{\sigma}_{a,t,n}^{2} - \sigma_{a,t}^{2}
&\leq& \left[\lambda_{\inf}^{-1}\lambda_{\sup}^{2}+\lambda_{\sup}\right]\sum_{b=1}^{p}\sum_{\tau\notin
       \{1,\ldots,n\}}|\beta_{(\tau,b)\shortarrow(t,a)}|. 
\end{eqnarray*}
\end{lemma}
{\bf PROOF} First we note that since
$[\mathcal{H}_{n}-X_{t}^{a}]\subseteq  [\mathcal{H}-X_{t}^{(a)}]$, then $\widetilde{\sigma}_{a,t,n}^{2} \geq 
\sigma_{a,t}^{2}$
and $0\leq  \widetilde{\sigma}_{a,t,n}^{2} - \sigma_{a,t}^{2}$.
To prove the result, we  recall if $P_{\mathcal{G}}(Y)$ is
the projection of $Y$ onto $\mathcal{G}$, then 
\begin{eqnarray*}
\var[Y - P_{\mathcal{G}}(Y)] = \var[Y] - \cov[Y, P_{\mathcal{G}}(Y)].
\end{eqnarray*}
Using the above, with $\mathcal{G}_{n} = \mathcal{H}_{n}-X_{t}^{(a)}$ and $\mathcal{G} = \mathcal{H}-X_{t}^{(a)}$ and taking differences gives 
\begin{eqnarray*}
\widetilde{\sigma}_{a,t,n}^{2} - \sigma_{a,t}^{2}
  &=&\var[X_{t}^{(a)}-P_{\mathcal{H}_{n}-X_{t}^{(a)}}(X_{t}^{(a)})] 
- \var[X_{t}^{(a)}-P_{\mathcal{H}-X_{t}^{(a)}}(X_{t}^{(a)})] \\
 &=& \cov\left[X_{t}^{(a)}, P_{\mathcal{H}-X_{t}^{(a)}}(X_{t}^{(a)}) -P_{\mathcal{H}_{n}-X_{t}^{a}}(X_{t}^{(a)})
     \right].
\end{eqnarray*}
Substituting the expressions for  
$P_{\mathcal{H}-X_{t}^{(a)}}(X_{t}^{(a)})$ and $P_{\mathcal{H}_{n}-X_{t}^{(a)}}(X_{t}^{(a)})$ in (\ref{eq:PX}) and (\ref{eq:PX1}) into
the above gives 
\begin{eqnarray*}
 \widetilde{\sigma}_{a,t,n}^{2} - \sigma_{a,t}^{2} &=& \cov\left[X_{t}^{(a)}, \sum_{b=1}^{p}\sum_{\tau=1}^{n}\left(\beta_{(\tau,b)\shortarrow
                                     (t,a)}-\theta_{(\tau,b)\shortarrow
    (t,a)}\right)X_{\tau}^{(b)}\right] + \\
 &&\cov\left[X_{t}^{(a)},
    \sum_{b=1}^{p}\sum_{\tau\notin
    \{1,\ldots,n\}}\beta_{(\tau,b)\shortarrow(t,a)}X_{\tau}^{(b)}\right]\\
 &=& \sum_{b=1}^{p}\sum_{\tau=1}^{n}\left(\beta_{(\tau,b)\shortarrow
                                     (t,a)}-\theta_{(\tau,b)\shortarrow
    (t,a)}\right)\cov[X_{\tau}^{(b)},X_{t}^{(a)}] + \\
 &&    \sum_{b=1}^{p}\sum_{\tau\notin
    \{1,\ldots,n\}}\beta_{(\tau,b)\shortarrow(t,a)}\cov[X_{\tau}^{(b)},X_{t}^{(a)}].
\end{eqnarray*}
Applying the Cauchy-Schwarz inequality to the above gives
\begin{eqnarray*}
\widetilde{\sigma}_{a,t,n}^{2} - \sigma_{a,t}^{2} 
&\leq& \left[\sum_{b=1}^{p}\sum_{\tau=1}^{n}\left(\beta_{(\tau,b)\shortarrow
                                     (t,a)}-\theta_{(\tau,b)\shortarrow
    (t,a)}\right)^{2}\right]^{1/2}\left[\sum_{b=1}^{p}\sum_{\tau=1}^{n}\cov[X_{\tau}^{(b)},X_{t}^{(a)}]^{2}
       \right]^{1/2} \\
&& + \left[\sum_{b=1}^{p}\sum_{\tau\notin
    \{1,\ldots,n\}}\beta_{(\tau,b)\shortarrow(t,a)}^{2}\right]^{1/2}\left[\sum_{b=1}^{p}\sum_{\tau\notin
    \{1,\ldots,n\}}\cov[X_{\tau}^{(b)},X_{t}^{(a)}]^{2}\right]^{1/2}.
\end{eqnarray*}
Applying the bound in Lemma \ref{lemma:meyer}  to the first term on
the RHS and using that the sum of the covariances squared are bounded
by $\lambda_{\sup}$ we have 
\begin{eqnarray*}
\widetilde{\sigma}_{a,t,n}^{2} - \sigma_{a,t}^{2} 
&\leq& \lambda_{\inf}^{-1}\lambda_{\sup}^{2}\sum_{b=1}^{p}\sum_{\tau\notin
       \{1,\ldots,n\}}|\beta_{(\tau,b)\shortarrow(t,a)}|+ \lambda_{\sup}\left[\sum_{b=1}^{p}\sum_{\tau\notin
    \{1,\ldots,n\}}\beta_{(\tau,b)\shortarrow(t,a)}^{2}\right]^{1/2}
  \\
&\leq& \left[\lambda_{\inf}^{-1}\lambda_{\sup}^{2}+\lambda_{\sup}\right]\sum_{b=1}^{p}\sum_{\tau\notin
       \{1,\ldots,n\}}|\beta_{(\tau,b)\shortarrow(t,a)}|,
\end{eqnarray*}
where the above follows from the fact that the $\ell_{2}$-norm of a
vector is bounded from above by the $\ell_{1}$-norm. This gives the required result. \hfill $\Box$

\vspace{2mm}
\noindent We use Lemmas \ref{lemma:meyer}, \ref{lemma:meyer2} and
equation (\ref{eq:BTheta}) to prove Theorem \ref{lemma:meyer3a}.

\vspace{2mm}

\noindent {\bf PROOF of Theorem \ref{lemma:meyer3a}} We first obtain a
bound for the sum of the regression coefficients
\begin{eqnarray*}
\sum_{b=1}^{p}\sum_{\tau\notin
       \{1,\ldots,n\}}|\beta_{(\tau,b)\shortarrow(t,a)}| 
&=& \sum_{b=1}^{p}\left[\sum_{\tau=n+1}^{\infty} + \sum_{\tau=-\infty}^{0}\right]
|\beta_{(\tau,b)\shortarrow(t,a)}|.
\end{eqnarray*}
Under Assumption \ref{assum:invcovarianceK} we have
\begin{eqnarray*}
\sum_{b=1}^{p}\sum_{\tau\notin
       \{1,\ldots,n\}}|\beta_{(\tau,b)\shortarrow(t,a)}| 
&\leq& \lambda_{\inf}^{-1}\sum_{b=1}^{p}\left[\sum_{\tau=n+1}^{\infty} +
       \sum_{\tau=-\infty}^{0}\right]\frac{1}{\ell(\tau-t)} \\
&=& \lambda_{\inf}^{-1}\sum_{b=1}^{p}\sum_{\tau=n+1}^{\infty} \frac{1}{\ell(\tau-t)}+\sum_{b=1}^{p}\sum_{\tau=-\infty}^{0}
       \frac{1}{\ell(\tau-t)} \\
&=& \lambda_{\inf}^{-1}\sum_{b=1}^{p}\sum_{j=n+1-t}^{\infty} \frac{1}{\ell(j)} + \sum_{b=1}^{p}\sum_{\tau=-\infty}^{-t}
       \frac{1}{\ell(j)}.
\end{eqnarray*}
Under  Assumption  \ref{assum:invcovarianceK} we have  for $r>0$
  $\lambda_{\inf}^{-1}\sum_{j=r}^{\infty}\ell(j) \leq
  \lambda_{\inf}^{-1} r^{-K}\sum_{j=r}^{\infty}j^{K}\ell(j)^{-1} \leq C_{\ell}r^{-K}$, where $C_{\ell} =
  \lambda_{\inf}^{-1}\sum_{j\in \mathbb{Z}}\ell(j)^{-1}$.  And by a
  similar argument for $r<0$, $\lambda_{\inf}^{-1}\sum_{j=-\infty}^{r}\ell(j) \leq
  C_{\ell}|r|^{-K}$. Applying these two bounds to the above we have 
\begin{eqnarray}
\label{eq:betabound}
\sum_{b=1}^{p}\sum_{\tau\notin
       \{1,\ldots,n\}}|\beta_{(\tau,b)\shortarrow(t,a)}| \leq  2pC_{\ell} \min(|n+1-t|,|t|)^{-K}
\end{eqnarray}
We use this inequality to prove the result. 

We return to  (\ref{eq:rowdiffs}) which gives the bound
\begin{eqnarray*}
 \left\| [\widetilde{D}_{n}]_{(a-1)n+t,\cdot} -  [D_{n}]_{(a-1)n+t,\cdot}\right\|_{1} 
&\leq& \frac{|\sigma_{a,t}^{2}-\widetilde{\sigma}_{a,t}^{2}|}{\widetilde{\sigma}_{a,t}^{2}\sigma_{a,t}^{2}} 
\left(1+ \|\underline{\Theta}^{(a,t)}_{n}\|_{1}\right)
+\frac{1}{\sigma_{a,t}^{2}}\left\|\underline{B}^{(a,t)}_{n}-\underline{\Theta}^{(a,t)}_{n}\right\|_{1}.
\end{eqnarray*}
Substituting  the bounds in (\ref{eq:BTheta}) and Lemma
\ref{lemma:meyer2} into the above gives 
\begin{eqnarray*}
&& \left\| [\widetilde{D}_{n}]_{(a-1)n+t,\cdot} -
  [D_{n}]_{(a-1)n+t,\cdot}\right\|_{1} \\ 
&\leq&
       \frac{1}{\widetilde{\sigma}_{a,t}^{2}\sigma_{a,t}^{2}}\left[\lambda_{\inf}^{-1}\lambda_{\sup}^{2}+\lambda_{\sup}\right]
\left(1+ \|\underline{\Theta}^{(a,t)}_{n}\|_{1}\right)
\sum_{b=1}^{p}\sum_{\tau\notin
       \{1,\ldots,n\}}|\beta_{(\tau,b)\shortarrow(t,a)}| \\ 
&& 
+\frac{1}{\sigma_{a,t}^{2}}(np)^{1/2}\lambda_{\inf}^{-1}\lambda_{\sup}\sum_{b=1}^{p}\sum_{\tau\notin
       \{1,\ldots,n\}}|\beta_{(\tau,b)\shortarrow(t,a)}| \\
&\leq& 
\left( \frac{1}{\widetilde{\sigma}_{a,t}^{2}\sigma_{a,t}^{2}}\left[\lambda_{\inf}^{-1}\lambda_{\sup}^{2}+\lambda_{\sup}\right]
\left(1+ \|\underline{\Theta}^{(a,t)}_{n}\|_{1}\right)+
\frac{1}{\sigma_{a,t}^{2}}(np)^{1/2}\lambda_{\inf}^{-1}\lambda_{\sup}
\right) \sum_{b=1}^{p}\sum_{\tau\notin
       \{1,\ldots,n\}}|\beta_{(\tau,b)\shortarrow(t,a)}|.
\end{eqnarray*}
We now bound $\|\underline{\Theta}^{(a,t)}_{n}\|_{1}$,
$\widetilde{\sigma}_{a,t}^{-2}$ and $\sigma_{a,t}^{-2}$ in terms of
the eigenvalues of $C$. 
By using (\ref{eq:Thetatilde})  we have $\widetilde{\underline{\Theta}}^{(a,t)}_{n} =
  \frac{1}{\widetilde{\sigma}_{a,t}^{2}}[1,-\underline{\Theta}^{(a,t)}_{n}]$,
  this gives the inequality
\begin{eqnarray*}
\|\underline{\Theta}^{(a,t)}_{n}\|_{1} \leq \widetilde{\sigma}_{a,t}^{2}\|\widetilde{\underline{\Theta}}^{(a,t)}_{n}\|_{1}.
\end{eqnarray*}
Since $\widetilde{\underline{\Theta}}^{(a,t)}_{n}$ are the
(rearranged) rows of $\widetilde{D}_{n} = C_{n}^{-1}$ and the smallest eigenvalue of
$C_{n}$ is bounded from below by $\lambda_{\inf}$ we have that 
\begin{eqnarray*}
\|\underline{\Theta}^{(a,t)}_{n}\|_{1} \leq \widetilde{\sigma}_{a,t}^{2}\lambda_{\inf}^{-1}.
\end{eqnarray*}
 Since $\widetilde{\sigma}_{a,t}^{2}\leq \var[X_{t}^{(a)}] = [{\bf C}_{t,t}]_{a,a}$,
 and $[{\bf C}_{t,t}]_{a,a} \leq \sum_{\tau}\|[{\bf
   C}_{t,\tau}]_{a,\cdot}\|_{2}^{2}\leq \lambda_{\sup}$
then $\widetilde{\sigma}_{a,t}^{2} \leq \lambda_{\sup}$, thus
$\|\underline{\Theta}^{(a,t)}_{n}\|_{1} \leq
\lambda_{\inf}^{-1}\lambda_{\sup}$.

By using (\ref{eq:sigmaatbound})  (from the start of Appendix
\ref{sec:covarianceproofs}) we have $\sigma_{a,t}^{-2}\leq
\lambda_{\inf}^{-1}$. Furthermore, by using the same arguments used to
show that $\sigma_{a,t}^{-2}\leq
\lambda_{\inf}^{-1}$ we can also show $\widetilde{\sigma}_{a,t}^{-2}\leq
\lambda_{\inf}^{-1}$. 

Altogether, these bounds with (\ref{eq:betabound}) give
\begin{eqnarray*}
&&\left\| [\widetilde{D}_{n}]_{(a-1)n+t,\cdot} -
  [D_{n}]_{(a-1)n+t,\cdot}\right\|_{1} \\
&\leq& 
\left( \frac{1}{\lambda_{\inf}^{2}}\left[\lambda_{\inf}^{-1}\lambda_{\sup}^{2}+\lambda_{\sup}\right]
\left(1+ \lambda_{\inf}^{-1}\lambda_{\sup}\right)+
\frac{1}{\lambda_{\inf}}(np)^{1/2}\lambda_{\inf}^{-1}\lambda_{\sup}
\right) 2pC_{\ell} \min(|n+1-t|,|t|)^{-K} \\
&=& O\left( (np)^{1/2}\min(|n+1-t|,|t|)^{-K}\right),
\end{eqnarray*}
where the constants above only depend on $\lambda_{\inf}$,
$\lambda_{\sup}$, $p$ and $C_{\ell}=\lambda_{\inf}^{-1}\sum_{j\in
  \mathbb{Z}}\ell(j)^{-1}$. Thus proving the result. \hfill $\Box$

\vspace{2mm}
\noindent {\bf PROOF of Proposition \ref{lemma:conditionalstat}} By
definition we have 
\begin{eqnarray*}
 [{\bf K}_{n}(\omega_{k_{1}},\omega_{k_{2}})]_{a,b} =[F_{n}^{*}\widetilde{D}_{a,b;n}F_{n}]_{k_{1},k_{2}} =
  \frac{1}{n}\sum_{\tau,t=1}^{n}[\widetilde{D}_{a,b;n}]_{t,\tau}\exp(-it\omega_{k_1}+i\tau\omega_{k_{2}})
\end{eqnarray*}
Replacing $\widetilde{D}_{a,b;n}$ with $D_{a,b;n}$ and using Theorem
\ref{lemma:meyer3a} gives 
\begin{eqnarray*}
 |([F_{n}^{*}(\widetilde{D}_{n}-D_{n})F_{n}]_{k_{1},k_{2}})_{a,b}| &\leq&
\frac{1}{n}\sum_{t=1}^{n}\sum_{\tau=1}^{n}\left|[\widetilde{D}_{a,b;n} - D_{a,b;n}]_{t,\tau}\right| \nonumber\\
&\leq&
       \frac{1}{n}\sum_{t=1}^{n}\frac{(np)^{1/2}}{\min(|t-n+1|,|t|)^{K}}
= O\left(\frac{(np)^{1/2}}{n^{K}}\right), 
\end{eqnarray*}
where the above holds for $K>1$. 
This gives 
\begin{eqnarray}
\label{eq:Kabbound}
 [{\bf K}_{n}(\omega_{k_{1}},\omega_{k_{2}})]_{a,b} &=&
  [F_{n}^{*}D_{a,b;n}F_{n}]_{k_{1},k_{2}} +
                                                        O\left(\frac{(np)^{1/2}}{n^{K}}\right)
  \nonumber\\
&=& \frac{1}{n}\sum_{t=1}^{n}\sum_{\tau=1}^{n}[D_{a,b}]_{t,\tau}
     \exp(-it\omega_{k_1}+i\tau\omega_{k_2}) + 
O\left(\frac{(np)^{1/2}}{n^{K}}\right).
\end{eqnarray}
Now we obtain an expression for the leading term in the RHS of the
above in terms of $\Gamma_{t}^{(a,b)}(\omega)$;
\begin{eqnarray*}
&&
   \frac{1}{n}\sum_{t=1}^{n}\sum_{\tau=1}^{n}[D_{a,b}]_{t,\tau}\exp(-it\omega_{k_1}+i\tau\omega_{k_2})
  \\
&=&
    \frac{1}{n}\sum_{t=1}^{n}\sum_{\tau=1}^{n}[D_{a,b}]_{t,\tau}\exp(-it(\omega_{k_1}-\omega_{k_2}))
\exp(i(\tau-t)\omega_{k_2}) \quad \textrm{ let }r = \tau-t\\
&=&\frac{1}{n}\sum_{t=1}^{n}\exp(-it(\omega_{k_1}-\omega_{k_2}))\sum_{r=1-t}^{n-t}[D_{a,b}]_{t,t+r}\exp(ir\omega_{k_2}) \\
&=&
    \frac{1}{n}\sum_{t=1}^{n}\exp(-it(\omega_{k_1}-\omega_{k_2}))\sum_{r=-\infty}^{\infty}[D_{a,b}]_{t,t+r}\exp(ir\omega_{k_2}) +
O\left( \frac{1}{n}\sum_{t=1}^{n}\left[\sum_{r=-\infty}^{1-t}+\sum_{r=n-t+1}^{\infty}\right]
\frac{1}{\ell(r)}\right) \\
&=& \frac{1}{n}\sum_{t=1}^{n}\exp(-it(\omega_{k_1}-\omega_{k_2}))\sum_{r=-\infty}^{\infty}[D_{a,b}]_{t,t+r}
 \exp(ir\omega_{k_2}) + O\left( \frac{1}{n}\sum_{r\in \mathbb{Z}}\frac{|r|}{\ell(r)}\right)\\
&=&
    \frac{1}{n}\sum_{t=1}^{n}\exp(-i(k_1-k_2)\omega_{t})\Gamma_{t}^{(a,b)}(\omega_{k_2})
    + O\left( n^{-1}\right).
\end{eqnarray*}
By a similar argument we can show that 
\begin{eqnarray*}
&&
   \frac{1}{n}\sum_{t=1}^{n}\sum_{\tau=1}^{n}[D_{a,b}]_{t,\tau}\exp(-it\omega_{k_1}+i\tau\omega_{k_2})
  \\
&=&
    \frac{1}{n}\sum_{\tau=1}^{n}\exp(i(k_2-k_1)\omega_{\tau})\Gamma_{\tau}^{(b,a)}(\omega_{k_1})^{*}
    + O\left( n^{-1}\right) \\
&=& \left[\frac{1}{n}\sum_{\tau=1}^{n}\exp(-i(k_2-k_1)\omega_{\tau})\Gamma_{\tau}^{(b,a)}(\omega_{k_1})\right]^{*}
    + O\left( n^{-1}\right)
\end{eqnarray*}
Therefore, since $O((np)^{1/2}/n^{K}) = O(1/n)$ when $K\geq 3/2$ and
$p$ is fixed, substituting the above into (\ref{eq:Kabbound}) we have
\begin{eqnarray*}
[{\bf K}_{n}(\omega_{k_1},\omega_{k_2})]_{a,b} &=&
   \frac{1}{n}\sum_{t=1}^{n}\exp(-i(k_1-k_2)\omega_{t})\Gamma_{t}^{(a,b)}(\omega_{k_2})+                                               
O\left(\frac{1}{n}\right)
\end{eqnarray*}
and 
\begin{eqnarray*}
[{\bf K}_{n}(\omega_{k_1},\omega_{k_2})]_{a,b} &=&
   \left[\frac{1}{n}\sum_{t=1}^{n}\exp(-i(k_2-k_1)\omega_{t})\Gamma_{t}^{(b,a)}(\omega_{k_1})\right]^{*}+                                               
O\left(\frac{1}{n}\right)
\end{eqnarray*}
this proves (\ref{eq:Gammatbound}).

To prove (\ref{eq:Gammatboundstat}) (under conditional stationarity) we use that
$\Gamma_{t}^{(a,b)}(\omega) = \Gamma^{(a,b)}(\omega)$ for all
$t$. Substituting this into (\ref{eq:Gammatbound}) gives 
\begin{eqnarray*}
[{\bf K}_{n}(\omega_{k_1},\omega_{k_2})]_{a,b} &=&\Gamma^{(a,b)}(\omega_{k_2})
   \frac{1}{n}\sum_{t=1}^{n}\exp(-i(k_1-k_2)\omega_{t})+
O\left(\frac{1}{n}\right),
\end{eqnarray*}
Now by using that 
\begin{eqnarray*}
\frac{1}{n}\sum_{t=1}^{n}\exp\left(-it\omega_{k_1-k_2,n}\right)
 = \left\{
\begin{array}{cc}
 0 & k_1-k_2 \notin n\mathbb{Z}\\
1 & k_1-k_2 \in n\mathbb{Z}\\
\end{array}
\right.
\end{eqnarray*}
immediately proves (\ref{eq:Gammatboundstat}).  \hfill $\Box$

\subsection{Proofs for Section \ref{sec:KLS}}\label{sec:LSproof}

\noindent{\bf PROOF of Proposition \ref{lemma:NetworkK}}
The proof follows from the definition of $K_{r}^{(a,b)}(\omega)$
\begin{eqnarray*}
K_{r}^{(a,b)}(\omega) &=&\int_{0}^{1}e^{-2\pi
    iru}\Gamma^{(a,b)}(u;\omega)du.
\end{eqnarray*}
\vspace{2mm}

\noindent {\bf PROOF of Proposition \ref{lemma:Hrab}}
To prove the result we use (\ref{eq:Gammatbound}) in Proposition 
\ref{lemma:conditionalstat} to give 
\begin{eqnarray*}
[{\bf K}_{n}(\omega_{k_{1}},\omega_{k_{2}})]_{a,b} = 
  \frac{1}{n}\sum_{t=1}^{n}\Gamma_{t}^{(a,b)}(\omega_{k_{1}})\exp(-it(\omega_{k_{1}}-\omega_{k_{2}})) + O\left(\frac{1}{n} \right)
\end{eqnarray*}
We replace $\Gamma_{t}^{(a,b)}(\omega_{k_2})$ with
$\Gamma^{(a,b)}(t/n,\omega_{k_2})$. Using the locally stationary
approximation bound in (\ref{eq:GLS}) we have 
\begin{eqnarray}
\label{eq:Kdef0}
[{\bf K}_{n}(\omega_{k_1},\omega_{k_2})]_{a,b} &=&
   \frac{1}{n}\sum_{t=1}^{n}\exp(-i(k_1-k_2)\omega_{t})\Gamma^{(a,b)}\left(\frac{t}{n},\omega_{k_2}\right)+                           
O\left(\frac{1}{n}\right).
\end{eqnarray}
This proves (\ref{eq:Kdef}). 

In order to prove (\ref{eq:Kdef2}) we
study the smoothness of $\Gamma^{(a,b)}(u,\omega_{})$ over $u$ and its
corresponding Fourier coefficients (keeping
$\omega$ fixed). We first observe that under Assumption \ref{assum:LS}
 we have that 
\begin{eqnarray*}
\frac{\partial \Gamma^{(a,b)}(u,\omega)}{\partial u} = \sum_{j\in
  \mathbb{Z}}\frac{d[{\bf D}_{j}(u)]_{a,b}}{du}
\exp(ij\omega).
\end{eqnarray*}
This leads to the bound
\begin{eqnarray}
\label{eq:gammaderivative}
\sup_{u,\omega}\left|\frac{\partial \Gamma^{(a,b)}(u,\omega)}{\partial u}\right| = \sup_{u}\sum_{j\in
  \mathbb{Z}}\left|\frac{d[{\bf D}_{j}(u)]_{a,b}}{du}\right|
\leq \sum_{j\in
  \mathbb{Z}}\ell(j)^{-1} <\infty.
\end{eqnarray}
We use this bound below. To simplify notation, we drop the $(a,b)$ and
$\omega$ in $ \Gamma^{(a,b)}(u,\omega)$ (as they do not play a role in
the bound).  In order to understand the rate of decay
of the Fourier coefficients of $\Gamma(\cdot)$ we note that  $\Gamma$
is a piecewise continuous $1$-periodic function (where
$\Gamma(u)=\Gamma(u+n)$ for all $n\in \mathbb{Z}$). Define the
Fourier coefficient
\begin{eqnarray*}
K_{r} = \int_{0}^{1}\Gamma(u)\exp(-ir2\pi u)du.
 \end{eqnarray*}
By using (\ref{eq:gammaderivative}) the derivative of
$\Gamma(\cdot)$ is bounded on the interior $(0,1)$ (it is unlikely to exist
at $0$ and $1$ since typically $\Gamma(0)\neq \Gamma(1)$). Thus
by  integration by parts we have the bound 
\begin{eqnarray}
\label{eq:HHR}
|K_{r}|\leq C|r|^{-1} \textrm{ for }r\neq
0. 
\end{eqnarray}
We now obtain the limit for
\begin{eqnarray*}
\frac{1}{n}\sum_{k=1}^{n}\Gamma(k/n)\exp(-ir\omega_{k}) \qquad  |r|\leq n/2.
\end{eqnarray*}
In particular, we show that 
\begin{eqnarray*}
\sup_{|r|\leq n/2}\left|\frac{1}{n}\sum_{k=1}^{n}\Gamma(k/n)\exp\left(-ir\frac{2\pi k}{n}\right) -
  \int_{0}^{1}\Gamma(u)\exp(-i2\pi u r)du \right| =O\left(\frac{1}{n}\right).
\end{eqnarray*}
Using the mean value theorem a crude bound for the above is 
 $O((|r|+1)/n)$. To obtain a uniform $O(1/n)$ bound  over $|r|\leq
 n/2$ requires a more subtle techique which we describe below.

Taking difference between the sum and
integral gives 
\begin{eqnarray*}
&&\sum_{k=1}^{n}\int_{(k-1)/n}^{k/n}\left[\Gamma(k/n)\exp\left(-ir\omega_{k}\right) - \Gamma(u)\left(-ir 2\pi u\right) \right]du \\
&=& \sum_{k=1}^{n}\int_{(k-1)/n}^{k/n}\left[\Gamma(k/n)-\Gamma(u)\right]\exp\left(-ir\omega_k\right)du + \\
&& \sum_{k=1}^{n}\int_{(k-1)/n}^{k/n}
\Gamma(u)\left[\exp\left(-ir\omega_k\right)-\exp\left(-ir 2\pi u\right) \right]du \\
 &=& I_{1} + I_{2}.
\end{eqnarray*}
It is clear by the Lipschitz continuity of $\Gamma$ and 
$|\exp(i2\pi u)|\leq 1$ that $I_{1} = O(1/n)$ uniformly over all
$r$. To obtain a similar bound for the second term we exploit the symmetries
of the $\cos$ and $\sin$ functions that make up
$\exp(-ir\omega_{k})$. 

We separate $I_{2}$ into its sin and cosine transforms 
\begin{eqnarray*}
I_{2} = I_{2,C} - iI_{2,S}
\end{eqnarray*}
where
\begin{eqnarray*}
I_{2,C} &=& \sum_{k=1}^{n-1}\int_{(k-1)/n}^{k/n}\Gamma(u)\left[\cos\left(r\frac{2\pi
    k}{n}\right)-\cos\left(r 2\pi u\right) \right] du \\
I_{2,S} &=& \sum_{k=1}^{n-1}\int_{(k-1)/n}^{k/n}\Gamma(u)\left[\sin\left(r\frac{2\pi
    k}{n}\right)-\sin\left(r 2\pi u\right) \right] du.
\end{eqnarray*}
We focus on the cosine transform 
\begin{eqnarray*}
I_{2,C} &=& \sum_{k=1}^{n-1}\int_{(k-1)/n}^{k/n}\Gamma(u)\left[\cos\left(r\frac{2\pi
    k}{n}\right)-\cos\left(r 2\pi u\right) \right] du \\
 &=&
     \sum_{k=1}^{n}\int_{0}^{1/n}\Gamma\left(u+\frac{k}{n}\right)\left[\cos\left(2\pi r\frac{
    k}{n}\right)-\cos\left(2\pi r\left[u+\frac{k}{n}\right]\right) \right] du.
\end{eqnarray*}
Applying the mean value theorem to the above term would give the bound
$O(|r|/n)$. Instead we turn the above integral into the differences of
 cosines \emph{and} $\Gamma$. We show that the resulting product of
 differences cancel the unwanted 
$|r|$ term. We split the sum $\sum_{k=1}^{n}f_{k}$ into a double sum 
$\sum_{j=0}^{r-1} \sum_{k=1}^{n/(2r)}f_{jn/(2r)+k} + \sum_{j=0}^{r-1}
\sum_{k=1}^{n/(2r)}f_{jn/(2r)+n/(2r)+ k}$. 
This gives the double sum
\begin{eqnarray*}
I_{2,C}  &=& \sum_{j=0}^{r-1}
     \sum_{k=1}^{n/(2r)}\int_{0}^{1/n}\Gamma\left(u+\frac{k+jn/r}{n}\right)\times \\
&& \bigg\{\cos\left(2\pi r\frac{
    k+jn/r}{n}\right)-\cos\left[2\pi r\left(u+\frac{k+jn/r}{n}\right)\right] \bigg\} du
\\
&& + \sum_{j=0}^{r-1}
     \sum_{k=1}^{n/(2r)}\int_{0}^{1/n}\Gamma\left(u+\frac{k+jn/r+n/(2r)}{n}\right)\times\\
&& \left[\cos\left(2\pi r\left[\frac{
    k+jn/r+n/(2r)}{n}\right]\right)-\cos\left(2\pi r\left[u+\frac{k+jn/r+n/(2r)}{n}\right]\right)
   \right] du 
\end{eqnarray*}
Now we use that 
\begin{eqnarray*}
&&\cos\left(2\pi r\frac{
    k+jn/r+ n/(2r)}{n}\right)-\cos\left(2\pi r[u+\frac{k+jn/r+n/(2r)}{n}]\right) \\
&=&
-\left[\cos\left(2\pi r\frac{k+jn/r}{n}\right)-\cos\left(2\pi
r[u+\frac{k+jn/r}{n}]\right)\right]
\end{eqnarray*}
and substitute this into the above to give
\begin{eqnarray*}
I_{2,C}&=& \sum_{j=0}^{r-1}
     \sum_{k=1}^{n/(2r)}\int_{0}^{1/n}\left[\Gamma\left(u+\frac{k+jn/r}{n}\right)
- \Gamma\left(u+\frac{k+jn/r}{n}+\frac{n/(2r)}{n}\right)
\right]\times \\
&&\left[\cos\left(2\pi r\frac{
    k}{n}\right)-\cos\left(2\pi r\left[u+\frac{k}{n}\right]\right) \right] du.
\end{eqnarray*}
Observe that $I_{2,C}$ is expressed as a double difference. We bound
both these differences using the  Lipschitz continuity of
$\Gamma(\cdot)$ and $\cos(r\cdot)$;
$|\Gamma(u)-
\Gamma(v)|\leq \sup|\Gamma^{\prime}(u)|\cdot|u-v|$
and $|\cos(ru)-\cos(rv)|\leq r|u-v|$. This yields the bound
\begin{eqnarray*}
I_{2,C} &\leq&
               \sup_{u}|\Gamma^{\prime}(u)| \sum_{j=0}^{r-1}
     \sum_{k=1}^{n/(2r)} \frac{1}{r}\times
               \frac{r}{n}\times\frac{1}{n} =  \sup_{u}|\Gamma^{\prime}(u)|n^{-1}
\end{eqnarray*}
which is a uniform bound for all $|r|\leq n/2$. The same bound holds
for the sin transform $I_{2,S}$. Altogether, the bounds for $I_{1}$,
$I_{2,C}$ and $I_{2,S}$ give
\begin{eqnarray}
\label{eq:T10}
\sup_{\omega}\sup_{|r|\leq n/2}
\left|\frac{1}{n}\sum_{k=1}^{n}\Gamma^{(a,b)}(k/n,\omega)\exp\left(ir\frac{2\pi k}{n}\right) -
  K_{r}^{(a,b)}(\omega)\right| = O(n^{-1}).
\end{eqnarray}
Thus for $|k_1-k_2|\leq n/2$ we have 
\begin{eqnarray*}
[{\bf K}_{n}(\omega_{k_1},\omega_{k_2})]_{a,b} &=& 
K_{k_{1}-k_{2}}^{(a,b)}(\omega_{k_2}) + O(n^{-1}).
\end{eqnarray*}
For $n/2<k_{1}-k_{2}<n$ we return to (\ref{eq:Kdef0}) 
\begin{eqnarray*}
[{\bf K}_{n}(\omega_{k_1},\omega_{k_2})]_{a,b} &=&
   \frac{1}{n}\sum_{t=1}^{n}\exp(i(k_1-k_2)\omega_{t})\Gamma^{(a,b)}\left(\frac{t}{n},\omega_{k_2}\right)+   O(n^{-1})                                     \\
&=& \frac{1}{n}\sum_{t=1}^{n}\exp(i(k_1-k_2-n)\omega_{t})
\Gamma^{(a,b)}\left(\frac{t}{n},\omega_{k_2}\right) + O(n^{-1})\\
&=& K_{k_{1}-k_{2}-n}^{(a,b)}(\omega_{k_2}) + O(n^{-1}) 
\end{eqnarray*}
where we use that $|k_{1}-k_{2}-n|<n/2$ and (\ref{eq:T10}). By a
similar argument for $-n<k_{1}-k_{2}<n/2$ we have
\begin{eqnarray*}
[{\bf K}_{n}(\omega_{k_1},\omega_{k_2})]_{a,b} 
&=& K_{k_{1}-k_{2}+n}^{(a,b)}(\omega_{k_2}) + O(n^{-1}), 
\end{eqnarray*}
this proves (\ref{eq:Kdef2}). \hfill $\Box$

\vspace{2mm}
\noindent {\bf PROOF of Proposition \ref{lemma:smoothness}}
To prove (\ref{eq:Hrdecay}) we use that under Assumption \ref{assum:LS}
 $\Gamma^{(a,b)}(\cdot;\omega)\in L_{2}[0,1]$. Thus
$\sum_{r}|K_{r}^{(a,b)}(\omega)|^{2}<\infty$, this immediately gives
(\ref{eq:Hrdecay}).  The bound $\sup_{\omega}|K_{r}^{(a,b)}(\omega)|\leq
C|r|^{-1}$ follows immediately from (\ref{eq:HHR}).

To prove (\ref{eq:Hrsmooth}) we use the mean value theorem 
\begin{eqnarray*}
|K_{r}^{(a,b)}(\omega_{1}) -
K_{r}^{(a,b)}(\omega_{2})| \leq
\sup_{\omega}|dK_{r}^{(a,b)}(\omega)/d\omega|\cdot|\omega_{1}-\omega_{2}|. 
\end{eqnarray*}
To
bound $\sup_{\omega}|dK_{r}^{(a,b)}(\omega)/d\omega|$ we use that 
\begin{eqnarray}
\label{eq:Djbound}
\left|\frac{d}{d\omega}K_{r}^{(a,b)}(\omega)\right| \leq \sum_{j\in \mathbb{Z}}(1+|j|)
\left|\int_{0}^{1}e^{-2\pi iru}[{\bf D}_{j}(u)]_{a,b}du\right|.
\end{eqnarray}
To bound the integral in the above we 
use integration by parts, this together with  Assumption
\ref{assum:LS} gives 
\begin{eqnarray*}
\left|\int_{0}^{1}e^{-2\pi iru}D_{j}^{(a,b)}(u)du\right| \leq 
\left\{
\begin{array}{cc}
C\ell(j)^{-1} & r = 0\\
C|r|^{-1}\ell(j)^{-1} & r \neq 0\\
\end{array}
\right..
\end{eqnarray*}
Substituting this into (\ref{eq:Djbound}) gives
\begin{eqnarray}
\label{eq:Djbound}
\left|\frac{d}{d\omega}K_{r}^{(a,b)}(\omega)\right| = 
\left\{
\begin{array}{cc}
C \sum_{j\in
  \mathbb{Z}}(1+|j|)\ell(j)^{-1} & r = 0\\
C|r|^{-1}\sum_{j\in
  \mathbb{Z}}(1+|j|)\ell(j)^{-1} & r \neq 0\\
\end{array}
\right.,
\end{eqnarray}
this immediately leads to the 
required result. \hfill $\Box$

%% file: Appendix_Section5.tex
\section{Assumptions
  \ref{assum:lambda}, \ref{assum:invcovarianceK} and \ref{assum:LS} 
and the tvVAR process}\label{sec:proofAR}

We show that under certain conditions the tvVAR process satisfies 
Assumptions  \ref{assum:lambda}, \ref{assum:invcovarianceK} and \ref{assum:LS}.
Then in Appendix \ref{example:1} we consider the inverse time-varying
spectral density of a tvVAR$(1)$ model.

\subsection{Assumptions and the tvVAR}

\subsubsection*{tvVAR and Assumption \ref{assum:lambda}}

We first show that Assumption \ref{assum:lambda} holds for the model $\underline{X}_{t} =
{\bf A}(t)\underline{X}_{t-1}+\underline{\varepsilon}_{t}$ where 
$\sup_{t}\|{\bf A}(t)\|<1-\delta$. We will show that both the largest
eigenvalues of $C$ and $D$ are finite (which proves  Assumption
\ref{assum:lambda}). We prove the result by showing the absolute sum
of each row of $C_{a,b}$ and $D_{a,b}$ is bounded for each $1\leq
a,b\leq p$.

We first obtain a bound for the
largest eigenvalue of $C$ in terms of the covariances. 
Since $\sup_{t}\|{\bf A}(t)\|<1-\delta$, $\underline{X}_{t}$ almost surely has the causal
solution
$\underline{X}_{t}=\sum_{\ell=0}^{\infty}[\prod_{j=0}^{\ell}{\bf A}(t-j)]\underline{\varepsilon}_{t-\ell}$
Using this expansion and $\sup_{t}\|{\bf A}(t)\|<1-\delta$ it is easily
shown that $|\cov[X_{t}^{(a)}, X_{\tau}^{(b)}]| \leq
K(1-\delta)^{|t-\tau|}$ for some finite constant $K$. Thus by using 
Gerschgorin Circle  Theorem we have 
\begin{eqnarray*}
\lambda_{\sup}(C) \leq Kp\sum_{r\in \mathbb{Z}}(1-\delta)^{|r|}. 
\end{eqnarray*}
Next we show that $\lambda_{\sup}(D)<\infty$. 
Under $\sup_{t}\|{\bf A}(t)\|<1-\delta$, the rows of ${\bf A}(t)$ are such that 
\begin{eqnarray*}
\sup_{t,a}\|[{\bf A}(t)]_{a\cdot}\|_{1} \leq p^{1/2}\sup_{t,a}\|[{\bf
  A}(t)]_{a\cdot}\|_{2} \leq p^{1/2}(1-\delta)^{}.
\end{eqnarray*}
Therefore,  by using the above, the representation of $D_{a,b}$ in (\ref{eq:Dabttau3}) together with 
Gerschgorin Circle  Theorem we have $\lambda_{\sup}(D) <\infty$. Thus
Assumptions \ref{assum:lambda} is satisfied.

We mention that 
$\sup_{t}\|{\bf A}(t)\|<1-\delta$ is a sufficient condition. It can be
relaxed to allow for a contraction on the spectral radius of ${\bf
  A}(t)$ and smoothness conditions on ${\bf A}(t)$ (see \cite{p:keu-95}).
The above result can be extended to finite order tvVAR$(d)$ models, by
rewriting the $p$-dimensional tvVAR$(d)$ model as a $pd$-dimension tvVAR$(1)$ model 
and placing similar conditions on the corresponding tvVAR$(1)$ matrix. 


\subsubsection*{tvVAR and Assumption \ref{assum:invcovarianceK}}
Suppose that $\{\underline{X}_{t}\}$ has a tvVAR$(\infty)$ representation where
$\sup_{t}\|{\bf A}_{j}(t)\|_2\leq \ell(j)^{-1}$ and $\{\ell(j)\}$ is a
monotonically increasing sequence as 
$|j|\rightarrow \infty$. If $\sum_{j\in
  \mathbb{Z}}|j|^{K+1}\ell(j)^{-1}<\infty$ (for some $K\geq 1$), then we
show below that Assumption 
\ref{assum:invcovarianceK} is satisfied.

To show this we require the following lemma.
 
\begin{lemma}\label{lemma:monotonic}
Suppose $\ell(j)^{-1}$ is monotonically decreasing as $|j|\rightarrow \infty$ with $\sum_{j\in
  \mathbb{Z}}|j|^{K}\ell(j)^{-1}<\infty$ (for $K\geq 2$). Then for all $r\in \mathbb{Z}$
\begin{eqnarray*}
\sum_{s=-\infty}^{\infty}\frac{1}{\ell(s)\ell(s+r)} \leq  \frac{1}{\widetilde{\ell}(r)}
\textrm{ where } 
\widetilde{\ell}(j) = \left[3\frac{1}{\ell(\lfloor|j|/2\rfloor)}\sum_{s\in\mathbb{Z}}\frac{1}{\ell(s)}\right]^{-1}
\end{eqnarray*}
and $\sum_{j\in \mathbb{Z}}|j|^{K}\widetilde{\ell}(j)^{-1}<\infty$.
\end{lemma}
PROOF. Without loss of generality we prove the result for $r\geq 0$.
We partition the sum $\sum_{s=-\infty}^{\infty}$
into three terms
\begin{eqnarray*}
\sum_{s=-\infty}^{\infty}\frac{1}{\ell(s)\ell(s+r)} &=& 
\sum_{s=0}^{\infty}\frac{1}{\ell(s)\ell(s+r)}  +
 \sum_{s=-\infty}^{-r}\frac{1}{\ell(s)\ell(s+r)}
+\sum_{s=-r+1}^{-1}\frac{1}{\ell(s)\ell(s+r)}
  \\
&=& II_{1} + II_{2}+II_{3}.
\end{eqnarray*}
Using that $1/\ell(j)$ is monotonically decreasing  as $j\rightarrow
\infty$ is is easily seen that
\begin{eqnarray}
\label{eq:III1}
II_{1} \leq
              \frac{1}{\ell(r)}\sum_{s=0}^{\infty}\frac{1}{\ell(s)}
              \textrm{ and }
II_{2} \leq \frac{1}{\ell(r)}\sum_{s=-\infty}^{0}\frac{1}{\ell(s)}.
\end{eqnarray}
To bound $II_{3}$ we use that for $-r/2\leq s \leq -1$ that
 $\ell(s)\leq\ell(\lfloor r/2\rfloor)$ 
and for $-r+1\leq s \leq -r/2$ then
$\ell(s+r)\leq\ell(\lfloor r/2\rfloor)$. Altogether this gives
the bound
\begin{eqnarray*}
II_{3} &\leq& \left[\sum_{s=-r/2}^{-1}+
               \sum_{s=-r+1}^{-r/2}\right]\frac{1}{\ell(s)}\frac{1}{\ell(s+r)}
 \leq  \frac{1}{\ell(\lfloor r/2\rfloor)}\sum_{s\in\mathbb{Z}}\frac{1}{\ell(s)}.
\end{eqnarray*}
The above bound together with (\ref{eq:III1}) (noting that
$\ell(r)>\ell(\lfloor r/2 \rfloor)$) gives
\begin{eqnarray}
\label{eq:DDD3}
\sum_{s=-\infty}^{\infty}\frac{1}{\ell(s)\ell(s+r)}
 \leq
  3\frac{1}{\ell(\lfloor(r/2\rfloor)}
\sum_{s\in\mathbb{Z}}\frac{1}{\ell(s)}.
\end{eqnarray}
For all $j$ define 
\begin{eqnarray*}
\widetilde{\ell}(j) = \left[3\frac{1}{\ell(\lfloor|j|/2\rfloor)}\sum_{s\in\mathbb{Z}}\frac{1}{\ell(s)}\right]^{-1}.
\end{eqnarray*}
Then from (\ref{eq:DDD3}) we have the bound
\begin{eqnarray*}
\sum_{s=-\infty}^{\infty}\frac{1}{\ell(s)\ell(s+r)}\leq \frac{1}{\widetilde{\ell}(r)}.
\end{eqnarray*}
Since by assumption $\sum_{j\in
  \mathbb{Z}}|j|^{K}\ell(|j|)^{-1}<\infty$, it is immediately clear from
the definition of $\widetilde{\ell}(j) $ that
$\sum_{j\in \mathbb{Z}}|j|^{K}\widetilde{\ell}(j)^{-1}<\infty$. This
proves the result. \hfill $\Box$

\vspace{2mm}
\begin{lemma}
Suppose $\underline{X}_{t}$ has a tvVAR$(\infty)$ representation
that satisfies (\ref{eq:tvVAR}) and ${\bf D}_{t,\tau}$ be defined as
in (\ref{eq:DARt}). If the  time-varying AR matrices satisfy
$\sup_{t}\|{\bf A}_{j}(t)\|_2\leq \ell(j)^{-1}$ where $\ell(j)^{-1}$ is
monotonically decreasing as $|j|\rightarrow \infty$ and $\sum_{j\in
  \mathbb{Z}}|j|^{K}\ell(j)^{-1}<\infty$, then $\sup_{t}\|{\bf
  D}_{t,t+j}\|_{1} \leq \widetilde{\ell}(j)^{-1}$ where 
$\sup_{t}\sum_{j\neq  0}|j|^{K}\widetilde{\ell}(j)^{-1}<\infty$.
\end{lemma}
\noindent {\bf PROOF} By using (\ref{eq:DARt}) we have
\begin{eqnarray*}
\|{\bf D}_{t,t+j}\|_{1} &\leq& \sup_{t}
\sum_{s=-\infty}^{\infty}\left\|\widetilde{{\bf
  A}}_{s}(t+s)^{\prime}\widetilde{{\bf
                                                     A}}_{j+s}(t+s)\right\|_{1} \\
 &\leq& \sum_{s=-\infty}^{\infty}\left\|\widetilde{{\bf
  A}}_{s}(t+s)\right\|_{2}\left\|\widetilde{{\bf
                                                     A}}_{j+s}(t+s)\right\|_{2}
\leq \sum_{s=-\infty}^{\infty}\frac{1}{\ell(s)\ell(s+j)}.
\end{eqnarray*}
Finally, from the above and by using Lemma \ref{lemma:monotonic} we have 
\begin{eqnarray*}
\|{\bf D}_{t,t+j}\|_{1} &\leq& \widetilde{\ell}(j)^{-1},
\end{eqnarray*}
this proves the result. \hfill $\Box$

\subsubsection*{tvVAR and Assumption \ref{assum:LS}}

We now show that under certain conditions on $\{{\bf A}_{j}(t)\}$ the tvVAR process satisfies 
Assumption \ref{assum:LS}. Define the matrices ${\bf
  A}_{j}:[0,1]\rightarrow \mathbb{R}^{p\times p}$, which are Lipschitz in the sense that 
\begin{eqnarray}
\label{eq:Asmooth}
\|{\bf A}_{j}(u) - {\bf A}_{j}(v)  \|_{1}\leq \frac{1}{\ell(j)}|u-v|
\end{eqnarray}
where $\ell(|j|)^{-1}$ is
monotonically decreasing as $|j|\rightarrow \infty$
with $\sum_{j\in
  \mathbb{Z}}j^{2}\ell(j)^{-1}<\infty$ 
and 
\begin{eqnarray}
\label{eq:Au}
\sup_{u}\|{\bf A}_{j}(u)\|_{1}\leq \ell(j)^{-1}.
\end{eqnarray}
Following \cite{p:dah-00a} we define the locally stationary tvVAR model as 
\begin{eqnarray}
\label{eq:tvVARLS}
\underline{X}_{t,n} = \sum_{j=1}^{\infty}{\bf A}_{j}\left(\frac{t}{n}\right)\underline{X}_{t-j,n} + \underline{\varepsilon}_{t},
\end{eqnarray}
where  $\{\underline{\varepsilon}_{t}\}_{t}$ are i.i.d random variables
with $\var[\underline{\varepsilon}_{t}]={\boldsymbol \Sigma}$ 
($0<\lambda_{\min}({\boldsymbol \Sigma})\leq
\lambda_{\max}({\boldsymbol \Sigma})<\infty$). 
To define the suitable ${\bf D}_{j}(u)$ (as given in Assumption
\ref{assum:LS}), we first define the auxillary, stationary process corresponding to $\underline{X}_{t,n}$;
\begin{eqnarray}
\label{eq:statARLS}
\underline{X}_{t}(u) = \sum_{j=1}^{\infty}{\bf A}_{j}(u)\underline{X}_{t-j}(u) + \underline{\varepsilon}_{t}.
\end{eqnarray}
The inverse covariance of $\{\underline{X}_{t}(u)\}_{t}$ is 
$D(u) = (D_{a,b}(u);1\leq a,b\leq p)$ with 
$[D_{a,b}(u)]_{t,\tau} = [{\bf D}_{t-\tau}(u)]_{a,b}$ and 
\begin{eqnarray}
\label{eq:Dju1}
{\bf D}_{t-\tau}(u)=\sum_{\ell=-\infty}^{\infty}\widetilde{{\bf
  A}}_{\ell}\left(u\right)^{\prime}\Sigma^{-1}\widetilde{{\bf
  A}}_{\ell+(\tau-t)}(u).
\end{eqnarray}
The spectral density matrix corresponding to 
 $\{\underline{X}_{t}(u)\}_{t}$ is 
\begin{eqnarray}
\label{eq:statARLSspec}
\boldsymbol{\Sigma}(u;\omega) = [I_{p}-\sum_{j=1}^{\infty}{\bf
  A}_{j}(u)\exp(-ij\omega)]^{-1}{\boldsymbol \Sigma}
([I_{p}-\sum_{j=1}^{\infty}{\bf A}_{j}(u)\exp(-ij\omega)]^{-1})^{*}.
\end{eqnarray}
Thus the  time-varying spectral precision
matrix associated with $\{X_{t,n}\}_{t}$ is
$\boldsymbol{\Gamma}(t/n;\omega)
=\boldsymbol{\Sigma}(t/n;\omega)^{-1}$ (see Section \ref{sec:KLS}).
In the following lemma we show that the time series $\{\underline{X}_{t,n}\}_{t}$
satisfies Assumption \ref{assum:LS}.

\begin{lemma}\label{lemma:ARsmooth}
Suppose that the time series $\{\underline{X}_{t,n}\}_{t}$ has the
representation in (\ref{eq:tvVARLS}), where the tvVAR matrices satisfy conditions
(\ref{eq:Asmooth}) and (\ref{eq:Au}). 
Let ${\bf D}_{t,\tau}$ and ${\bf D}_{t-\tau}(u)$ be defined as in
(\ref{eq:DARt}) and (\ref{eq:Dju1}). 
Then 
\begin{eqnarray}
\label{eq:DttauDu}
\left\|{\bf D}_{t,\tau} - {\bf D}_{t-\tau}\left(\frac{t+\tau}{2n}\right)
  \right\|_{1} \leq  \frac{|t-\tau|+1}{n\widetilde{\ell}(t-\tau)} \textrm{ and }
\left\| {\bf D}_{t-\tau}(u) -  {\bf D}_{t-\tau}(v)  \right\|_{1} \leq  \frac{|u-v|}{n\widetilde{\ell}(t-\tau)} 
\end{eqnarray}
and $\sup_{u}\sum_{j\in \mathbb{Z}}j^{2}\|{\bf D}_{j}(u)\|_{1}<\infty$, 
where $\widetilde{\ell}(|j|)$ is monotonically increasing as
$|j|\rightarrow \infty$ and 
$\sum_{j\in \mathbb{Z}}j^{2}\cdot\widetilde{\ell}(|j|)^{-1}<\infty$. 

Further, if $\sup_{u}\|d{\bf A}_{j}(u)/du\|_{1}\leq 1/\ell(j)$ then
\begin{eqnarray}
\label{eq:derD}
\sup_{u}\left\|\frac{d{\bf D}_{t-\tau}(u)}{du}\right\|_{1} = \frac{1}{\widetilde{\ell}(t-\tau)}.
\end{eqnarray}
\end{lemma}
{\bf PROOF} We first prove (\ref{eq:DttauDu}). We recall that 
\begin{eqnarray*}
{\bf D}_{t,\tau}=\sum_{\ell=-\infty}^{\infty}\widetilde{{\bf
  A}}_{\ell}\left(t+s\right)^{\prime}\Sigma^{-1}\widetilde{{\bf
  A}}_{(\tau-t)+\ell}(t+s),
\end{eqnarray*}
and 
\begin{eqnarray}
\label{eq:Ajuu}
\|{\bf A}_{j}(u) - {\bf A}_{j}(v)  \|_{1}&\leq& \frac{1}{\ell(j)}|u-v| 
\textrm{ and } 
\sup_{u}\|{\bf A}_{j}(u)\|_{1} \leq \ell(j)^{-1}.
\end{eqnarray}
To simplify the notation (the proof does not change), we set ${\bf H}=I_{p}$.
Using the above and evaluating the difference gives 
\begin{eqnarray*}
{\bf D}_{t,\tau} - {\bf D}_{t-\tau}(u)
&=&  \sum_{s=-\infty}^{\infty}
\left(\widetilde{{\bf  A}}_{s}((t+s)/n)^{\prime}
\widetilde{{\bf A}}_{s+(\tau-t)}((t+s)/n) - 
\widetilde{{\bf  A}}_{s}(u)^{\prime} 
\widetilde{{\bf A}}_{s+(\tau-t)}(u)
\right) = I_{1} + I_{2},
\end{eqnarray*}
where
\begin{eqnarray*}
I_{1}&=& \sum_{s=-\infty}^{\infty}
\left(\widetilde{{\bf  A}}_{s}((t+s)/n)^{\prime} 
\widetilde{{\bf A}}_{s+(\tau-t)}((t+s)/n) - 
\widetilde{{\bf  A}}_{s}((t+s)/n)^{\prime} 
\widetilde{{\bf A}}_{s+(\tau-t)}(u)
\right)  \\
I_{2}& = & \sum_{s=-\infty}^{\infty}
\left(\widetilde{{\bf  A}}_{s}((t+s)/n)^{\prime} 
\widetilde{{\bf A}}_{s+(\tau-t)}(u) - 
\widetilde{{\bf  A}}_{s}(u)^{\prime} 
\widetilde{{\bf A}}_{s+(\tau-t)}(u)
\right).
\end{eqnarray*}
The two bounds are very similar, we focus on obtaining a bound for $I_{1}$.
By the Cauchy-Schwarz inequality and that $\|\cdot\|_{2}\leq
\|\cdot\|_{1}$ we have
\begin{eqnarray*}
\|I_{1}\|_{1} &\leq& \sum_{s=-\infty}^{\infty}
\left\|\widetilde{{\bf  A}}_{s}((t+s)/n)^{\prime} \left[
                     \widetilde{{\bf
                     A}}_{s+(\tau-t)}((t+s)/n)-\widetilde{{\bf
                     A}}_{s+(\tau-t)}(u)\right]^{}\right\|_{1}
  \\
&\leq& \sum_{s=-\infty}^{\infty}
\left\|\widetilde{{\bf  A}}_{s}(t/n) \right\|_{2}\left\|
                     \widetilde{{\bf
                     A}}_{s+(\tau-t)}((t+s)/n)-\widetilde{{\bf
                     A}}_{s+(\tau-t)}(u)\right\|_{2}\\
&\leq& \sum_{s=-\infty}^{\infty}
\left\|\widetilde{{\bf  A}}_{s}(t/n) \right\|_{1}\left\|
                     \widetilde{{\bf
                     A}}_{s+(\tau-t)}((t+s)/n)-\widetilde{{\bf
                     A}}_{s+(\tau-t)}(u)\right\|_{1} \\
\end{eqnarray*}
Substituting the bounds for ${\bf A}_{j}(u)$ given in (\ref{eq:Ajuu})
into the above we have 
\begin{eqnarray*}
\|I_{1}\|_{1}&\leq & \sum_{s=-\infty}^{\infty}\frac{1}{\ell(s)}\frac{1}{\ell(s+\tau-t)}\left|\frac{t+s}{n}-u\right|.
\end{eqnarray*}
By the same argument we have 
\begin{eqnarray*}
\|I_{2} \|_1\leq \sum_{s=-\infty}^{\infty}\frac{1}{\ell(s)}\frac{1}{\ell(s+\tau-t)}\left|\frac{t+s}{n}-u\right|.
\end{eqnarray*}
Thus 
\begin{eqnarray*}
\|{\bf D}_{t,\tau} - {\bf D}_{t-\tau}(u)\|_{1} \leq
  2\sum_{s=-\infty}^{\infty}\frac{1}{\ell(s)}\frac{1}{\ell(s+\tau-t)}\cdot|\frac{t+s}{n}-u|.
\end{eqnarray*}
Setting $u=(t+\tau)/(2n)$ gives 
\begin{eqnarray*}
\left\|{\bf D}_{t,\tau} - {\bf
  D}_{t-\tau}\left(\frac{t+\tau}{2n}\right)\right\|_{1} &\leq& 
  \sum_{s=-\infty}^{\infty}\frac{1}{\ell(s)}\frac{1}{\ell(s+\tau-t)}\left(\frac{|\tau-t|}{2n}
                                                               +
                                                               \frac{|s|}{n}\right) \\
&=&\frac{|t-\tau|}{2n}\sum_{s=-\infty}^{\infty}\frac{1}{\ell(s)}\frac{1}{\ell(s+\tau-t)}
    +\frac{1}{n}\sum_{s=-\infty}^{\infty}\frac{|s|}{\ell(s)}\frac{1}{\ell(s+\tau-t)}
\end{eqnarray*}
Now by using Lemma \ref{lemma:monotonic} we have
we have the bound
\begin{eqnarray}
\label{eq:DW1}
\left\|{\bf D}_{t,\tau} - {\bf
  D}_{t-\tau}\left(\frac{t+\tau}{2n}\right)\right\|_{1}  \leq
  \frac{|t-\tau|}{2n\widetilde{\ell}(t-\tau)} +  \frac{1}{n\widetilde{\ell}(t-\tau)}.
\end{eqnarray}
Since by assumption $\sum_{j\in
  \mathbb{Z}}(j^{2}+1)/\ell(|j|)<\infty$, it is immediately clear from
the definition of $\widetilde{\ell}(j) $ that
$\sum_{j\in \mathbb{Z}}(j^{2}+1)\widetilde{\ell}(j)^{-1}<\infty$. 
Under the stated assumptions in (\ref{eq:Ajuu}) and using Lemma \ref{lemma:monotonic} 
we can show
\begin{eqnarray}
\label{eq:DW2}
\left\| {\bf D}_{t-\tau}(u) -  {\bf D}_{t-\tau}(v)  \right\|_{1} \leq
  \frac{|u-v|}{n\widetilde{\ell}(t-\tau)} 
\end{eqnarray}
and $\sup_{u}\sum_{j\in \mathbb{Z}}j^{2}\|{\bf
  D}_{j}(u)\|_{1}<\infty$. (\ref{eq:DW1}) and (\ref{eq:DW2}) together 
 prove (\ref{eq:DttauDu}).

We now prove (\ref{eq:derD}). The elementwise derivative of ${\bf D}_{j}(u)$ is
\begin{eqnarray*}
\frac{d}{du}{\bf D}_{j}(u)=\sum_{s=-\infty}^{\infty}\frac{d}{du}\widetilde{{\bf
  A}}_{s}\left(u\right)^{\prime}\widetilde{{\bf
  A}}_{(\tau-t)+s}(u) + \sum_{s=-\infty}^{\infty}\widetilde{{\bf
  A}}_{s}\left(u\right)^{\prime}\frac{d}{du}\widetilde{{\bf
  A}}_{(\tau-t)+s}(u).
\end{eqnarray*}
Using the conditions in (\ref{eq:Ajuu}), $\sup_{u}\|d{\bf A}_{j}(u)/du\|_{1}\leq 1/\ell(j)$
and Lemma \ref{lemma:monotonic} we can show that
\begin{eqnarray*}
\left\|\frac{d}{du}{\bf D}_{j}(u)\right\|_{1} \leq
2\sum_{s\in \mathbb{Z}}\frac{1}{\ell(s)\ell(s+j)},
\end{eqnarray*}
this gives (\ref{eq:derD}). \hfill $\Box$

\subsection{Example: Locally stationary time-varying VAR$(1)$}\label{example:1}

We consider the locally stationary time-varying VAR$(1)$ model
\begin{eqnarray*}
\underline{X}_{t,n} = {\bf
  A}\left(\frac{t}{n}\right)\underline{X}_{t-1,n}+\underline{\varepsilon}_{t}\quad t \in \mathbb{Z},
\end{eqnarray*}
where $\{\underline{\varepsilon}_{t}\}$ are i.i.d. with $\var[\underline{\varepsilon}_{t}]=I_{p}$. 
We assume that the matrix ${\bf A}(u)$ satisfies (\ref{eq:Asmooth}) and
(\ref{eq:Au}).  Using (\ref{eq:statARLSspec}) the time-varying
spectral precision matrix corresponding to $\{\underline{X}_{t,n}\}_{}$ is 
${\boldsymbol \Gamma}(u;\omega) = [I_{p}- {\bf
  A}(u)\exp(-i\omega)]^{*}[I_{p}- {\bf A}(u)\exp(-i\omega)]$.  
We partition ${\boldsymbol \Gamma}(u;\omega)$ into the conditional
stationary and nonstationary matrices.

 Define the set $\mathcal{S}\subseteq
\{1,\ldots,p\}$ where for all $a\in \mathcal{S}$ the
columns $[{\bf A}(t)]_{\cdot,a}$ do not depend on $t$. 
Then $\{X_{t}^{(a)};a\in \mathcal{S}\}$ 
is a conditionally stationary subgraph and $D_{\mathcal{S}} = \{D_{a,b};a,b\in \mathcal{S}\}$
is a block Toeplitz matrix (see Corollary
\ref{cor:nodes}).  By  Lemma \ref{lemma:diagonal}
the integral kernel associated with  $D_{\mathcal{S}}$ is 
${\boldsymbol
  \Gamma}_{\mathcal{S},\mathcal{S}}(\omega)\delta_{\omega,\lambda}$.
We obtain an expression for ${\boldsymbol  \Gamma}_{\mathcal{S},\mathcal{S}}(\omega)$ 
below. 

We denote the set $\mathcal{S}$ as  $\mathcal{S} =
\{a_{1},\ldots,a_{|\mathcal{S}|}\}$, where $|\mathcal{S}|$ denotes the cardinality of
$\mathcal{S}$. Define the
$p\times |\mathcal{S}|$ matrix 
${\bf A}_{\mathcal{S}}$ where ${\bf A}_{\mathcal{S}} = ([{\bf
  A}_{\mathcal{S}}]_{\cdot,r} = [{\bf A}(0)]_{\cdot,a_{r}};a_{r}\in
\mathcal{S}, 1\leq r \leq |\mathcal{S}|)$,
let $I_{p,|\mathcal{S}|}$ denote the
$p\times |\mathcal{S}|$ ``indicator'' matrix which is comprised of
zeros except at  the entries $\{(r,a_{r});a_{r}\in \mathcal{S}\}$ where
$[I_{p,|\mathcal{S}|}]_{r,a_{r}} = 1$. Then 
\begin{eqnarray}
\label{eq:Gamma2}
{\boldsymbol \Gamma}_{\mathcal{S}, \mathcal{S}}(\omega) = \left[I_{p,|\mathcal{S}|} -
  {\bf A}_{\mathcal{S}}\exp(-i\omega)\right]^{*}
  \left[I_{p,|\mathcal{S}|} - 
 {\bf A}_{\mathcal{S}}\exp(-i\omega)\right].
\end{eqnarray}
Using ${\boldsymbol \Gamma}_{\mathcal{S}, \mathcal{S}}(\omega)$ we can
deduce the  the partial spectral coherence for conditionally stationary nodes and pairs
(see Lemma \ref{lemma:spectral3}).

We now obtain the nonstationary submatrices in ${\boldsymbol \Gamma}(u;\omega)$. Let 
$\mathcal{S}^{\prime} = \{b_{1},\ldots,b_{|\mathcal{S}^{\prime}|}\}$
denote the complement of $\mathcal{S}$. 
Analogous to ${\bf
  A}_{\mathcal{S}}$ and $I_{p,|\mathcal{S}|}$,  we define the 
$p\times |\mathcal{S}^{\prime}|$ dimensional matrices
${\bf  A}_{\mathcal{S}^{\prime}}(u) = ([{\bf
  A}_{\mathcal{S}^{\prime}}(u)]_{\cdot,r}=[{\bf A}(u)]_{\cdot,b_{r}};b_{r}\in
\mathcal{S}^{\prime},1\leq r\leq |\mathcal{S}_{r}|)$  and 
$I_{p\times |\mathcal{S}^{\prime}|}$ which is comprised of
zeros except at  the entries $\{(r,b_{r});b_{r}\in \mathcal{S}^{\prime}\}$
where $[I_{p,|\mathcal{S}^{\prime}|}]_{r,b_r} = 1$. A rearranged
version of ${\boldsymbol \Gamma}(u;\omega)$ (which for simplicity we
call ${\boldsymbol \Gamma}(u;\omega)$) is 
\begin{eqnarray*}
{\boldsymbol \Gamma}(u;\omega) = 
\left(
\begin{array}{cc}
\Gamma_{\mathcal{S}, \mathcal{S}}(\omega)
  &\Gamma_{\mathcal{S},\mathcal{S}^{\prime}}(u;\omega) \\
\Gamma_{\mathcal{S},\mathcal{S}^{\prime}}(u;\omega)^{*} &
                                                          \Gamma_{\mathcal{S}^{\prime}, \mathcal{S}^{\prime}}(u;\omega) \\
\end{array}
\right),
\end{eqnarray*}
where
\begin{eqnarray*}
\Gamma_{\mathcal{S}, \mathcal{S}}(\omega) &=& \left[I_{p,|\mathcal{S}|} -
  {\bf A}_{\mathcal{S}}\exp(-i\omega)\right]^{*}
  \left[I_{p,|\mathcal{S}|} - 
 {\bf A}_{\mathcal{S}}\exp(-i\omega)\right], \\
\Gamma_{\mathcal{S},\mathcal{S}^{\prime}}(u;\omega) &=& \left[I_{p,|\mathcal{S}|} -
  {\bf A}_{\mathcal{S}}\exp(-i\omega)\right]^{*}
  \left[I_{p,|\mathcal{S}^{\prime}|} - 
 {\bf A}_{\mathcal{S}^{\prime}}(u)\exp(-i\omega)\right] \\
\textrm{ and } \Gamma_{\mathcal{S}^{\prime}, \mathcal{S}^{\prime}}(u;\omega) &=& \left[I_{p,|\mathcal{S}^{\prime}|} -
  {\bf A}_{\mathcal{S}^{\prime}}(u)\exp(-i\omega)\right]^{*}
  \left[I_{p,|\mathcal{S}^{\prime}|}- 
 {\bf A}_{\mathcal{S}^{\prime}}(u)\exp(-i\omega)\right]. 
\end{eqnarray*}
Using the above we can deduce $K_{r}^{(a,b)}(\omega)$ and thus
approximations to the entries of ${\bf
  K}_{n}(\omega_{k_1},\omega_{k_2})$. 
The system for Example \ref{example:running} is described in detail in
Appendix \ref{sec:runningexample}.

%% file: Appendix_example.tex
\section{Examples}\label{sec:runningexample}

In the following two sections we study the running 
time-varying AR$(1)$ example (introduced in Example \ref{example:running}). 
In Appendix \ref{app:global} we compare the nonstationary graph of a
piecewise stationary time series with its piecewise stationary gaphs.

\subsection{The tvVAR model and corresponding local spectral precision matrix}

For all $1\leq t\leq n$ the model is defined as

\begin{minipage}{0.8\textwidth}
\noindent
\begin{eqnarray*}
\left(
\begin{array}{c}
X_{t}^{(1)} \\
X_{t}^{(2)} \\
X_{t}^{(3)} \\
X_{t}^{(4)} \\
\end{array}
\right) = 
\left(
\begin{array}{cccc}
\alpha_{1}(t/n) & 0 & \alpha_{3} & 0 \\
\beta_{1} & \beta_{2} & 0 & \beta_{4} \\ 
 0 & 0 & \gamma_{3}(t/n) & 0 \\
0 & \nu_{2} & 0 & \nu_{4} \\ 
\end{array}
\right)
\left(
\begin{array}{c}
X_{t-1}^{(1)} \\
X_{t-1}^{(2)} \\
X_{t-1}^{(3)} \\
X_{t-1}^{(4)} \\
\end{array}
\right)
+
\left(
\begin{array}{c}
\varepsilon_{t}^{(1)} \\
\varepsilon_{t}^{(2)} \\
\varepsilon_{t}^{(3)} \\
\varepsilon_{t}^{(4)} \\
\end{array}
\right)
\end{eqnarray*}
\end{minipage}
\begin{minipage}{0.2\textwidth}
\includegraphics[scale = 0.4]{plotspaper/corrected_graph.png}
\end{minipage}
\vspace{1mm}

\noindent and $\{\underline{\varepsilon}_{t}\}_{t}$ are iid random vectors
with $\var[\underline{\varepsilon}_{t}]$ and
$\alpha_{1}(\cdot),\gamma_{3}(\cdot)\in L_{2}[0,1]$ and are Lipschitz continuous.
By using the results in Section \ref{sec:tvVAR} we obtain the network
on the right.

\noindent We now obtain the time-varying conditional spectral density
${\boldsymbol \Gamma}(u;\omega)$. Let 
\begin{eqnarray*}
[I_{2,4} - {\bf A}_{\mathcal{S}}\exp(-i\omega)] &=& 
\left(
\begin{array}{cc}
0 & 0 \\
1 - \beta_{2}e^{-i\omega} & -\beta_{4}e^{-i\omega} \\
0 & 0 \\
-\nu_{2}e^{-i\omega} & 1 - \nu_{4}e^{-i\omega} \\
\end{array}
\right)\\
\textrm{ and }
[I_{1,3} - {\bf A}_{\mathcal{S}^{\prime}}(u)\exp(-i\omega)] &=& \left(
\begin{array}{cccc}
1 - \alpha_{1}(u)e^{-i\omega} & -\alpha_{3}e^{-i\omega} \\
  -\beta_{1}e^{-i\omega} & 0\\
0 & 1-\gamma_{3}(u)e^{-i\omega} \\
0 & 0 \\
\end{array}
\right).
\end{eqnarray*}
Then 
\begin{eqnarray*}
{\boldsymbol \Gamma}(u;\omega) = 
\left(
\begin{array}{cc}
\Gamma_{(2,4), (2,4)}(\omega)  & \Gamma_{(2,4),(1,3)}(u;\omega)\\
\Gamma_{(2,4),(1,3)}(u;\omega)^{*} & \Gamma_{(1,3), (1,3)}(u;\omega) \\
\end{array}
\right)
\end{eqnarray*}
where 
\begin{eqnarray*}
\Gamma_{(2,4), (2,4)}(\omega) &=& \left[I_{2,4} -
  {\bf A}_{\mathcal{S}}\exp(-i\omega)\right]^{*}
  \left[I_{2,4} - 
 {\bf A}_{\mathcal{S}}\exp(-i\omega)\right] \\
&=& \left(
\begin{array}{cc}
|1 - \beta_{2}e^{-i\omega}|^{2}+\nu_{2}^{2} & -\beta_{4}e^{i\omega}(1
                                              - \beta_{2}e^{-i\omega})-\nu_{2}e^{-i\omega}(1-\nu_{4}e^{i\omega})\\
-\beta_{4}e^{-i\omega}(1 - \beta_{2}e^{i\omega}) -\nu_{2}e^{i\omega}(1-\nu_{4}e^{-i\omega})& |1 - \nu_{4}e^{-i\omega}|^{2}+\beta_{4}^{2}
\end{array} 
\right) \\
\Gamma_{(2,4),(1,3)}(u;\omega) &=& \left[I_{2,4} -
  {\bf A}_{\mathcal{S}}\exp(-i\omega)\right]^{*}
  \left[I_{1,3} - 
 {\bf B}_{\mathcal{S}^{\prime}}(u)\exp(-i\omega)\right] \\
&=& 
\left(
\begin{array}{cc}
-\beta_{1}e^{i\omega}[1-\beta_{2}e^{-i\omega}]  & 0 \\
\beta_{1}\beta_{4} & 0 \\
\end{array}
\right) \\
\Gamma_{(1,3), (1,3)}(u;\omega) &=& \left[I_{1,3} -
  {\bf A}_{\mathcal{S}^{\prime}}(u)\exp(-i\omega)\right]^{*}
  \left[I_{1,3} - 
 {\bf A}_{\mathcal{S}^{\prime}}(u)\exp(-i\omega)\right] \\
&=& 
\left(
\begin{array}{cc}
| 1 - \alpha_{1}(u)e^{-i\omega}|^{2} + \beta_{1}^{2} & -\alpha_{3}e^{i\omega}[1-\alpha_{1}(u)e^{-i\omega}] \\
-\alpha_{3}e^{-i\omega}[1-\alpha_{1}(u)e^{i\omega}] & |1-\gamma_{3}(u)e^{-i\omega}|^{2}+\alpha_{3}^{2}\\
\end{array}
\right). 
\end{eqnarray*}

\subsection{The partial spectral coherence}

Based on the results in Section \ref{sec:coherence},
we use $\Gamma_{(2,4), (2,4)}(\omega)$ to define the partial spectral
coherence for the conditionally stationary nodes and edges ($2$ and $4$). 
We observe that $\Gamma_{(2,4), (2,4)}(\omega)$ resembles the spectral density matrix of a stationary
vector moving average model of order one (or equivalently the inverse of a
vector autoregressive of order one). Using $\Gamma_{(2,4), (2,4)}(\omega)$, the 
 partial spectra for the conditionally stationary nodes
$2$ and $4$ are
\begin{eqnarray*}
\Gamma^{(2,2)}(\omega)^{-1} = \frac{1
  }{|1 - \beta_{2}e^{-i\omega}|^{2}+\nu_{2}^{2}}
  \textrm{ and }\Gamma^{(4,4)}(\omega)^{-1} = \frac{1}{|1 - \nu_{4}e^{-i\omega}|^{2}+\beta_{4}^{2}}.
\end{eqnarray*}
Furthermore, 
by using (\ref{eq:Rabomega}), 
the partial spectral coherence for the conditionally stationary edge $(2,4)$ is 
\begin{eqnarray*}
R_{2,4}(\omega) = 
-\frac{-\beta_{4}e^{i\omega}(1 - \beta_{2}e^{-i\omega}) -
\nu_{2}e^{-i\omega}(1-\nu_{4}e^{i\omega})}{\sqrt{\Gamma^{(2,2)}(\omega) \Gamma^{(4,4)}(\omega)}}
\end{eqnarray*}

\subsection{A comparision of networks in a piecewise stationary
  VAR$(1)$ model}\label{app:global}

Below we consider both the directed and undirected graphs for a
piece-wise stationary VAR$(1)$ model where $p=4$.  
 We suppose
that $\underline{X}_{t}$ is piece-wise stationary in the sense that 
\begin{eqnarray*}
\underline{X}_{t} = 
\left\{
\begin{array}{cc}
A_{1}\underline{X}_{t-1} + \underline{\varepsilon}_{t} & 1\leq t \leq
                                                         n/2 \\
A_{2}\underline{X}_{t-1} + \underline{\varepsilon}_{t} & n/2+1\leq t \leq
                                                         n \\
\end{array}
\right.
\end{eqnarray*}
where $\var[\underline{\varepsilon}_{t}]=I_{3}$ and
$\{\underline{\varepsilon}_{t}\}_{t}$ are random vectors,
\begin{eqnarray*}
A_{1}  = \left(
\begin{array}{cccc}
\alpha_1 & 0 & 0 & 0 \\
\beta_{1} & \beta_{2} & 0 & \beta_{4} \\
0 & 0 & \gamma_{3} & 0 \\
0 & \nu_2 & 0 & \nu_4\\
\end{array}
\right) \quad \textrm{and}\quad
A_{2} = \left(
\begin{array}{cccc}
\widetilde{\alpha}_1 & 0 & \alpha_3 & 0 \\
\beta_{1} & \beta_{2} & 0 & \beta_{4} \\
0 & 0 & \widetilde{\gamma}_{3} & 0 \\
0 & \nu_2 & 0 & \nu_4\\
\end{array}
\right).
\end{eqnarray*}
note that $\alpha_1\neq \widetilde{\alpha}_1$ and 
$\gamma_3\neq \widetilde{\gamma}_3$.
For $1\leq t \leq n$ the above model can be written as the
nonstationary model 
\begin{eqnarray*}
\label{eq:example1}
\left(
\begin{array}{c}
X_{t}^{(1)} \\
X_{t}^{(2)} \\
X_{t}^{(3)} \\
X_{t}^{(4)} \\
\end{array}
\right) = 
\left(
\begin{array}{cccc}
\alpha(t) & 0 & \alpha_3(t) & 0 \\
\beta_{1} & \beta_{2} & 0 & \beta_{4} \\
0 & 0 & \gamma_{}(t) & 0 \\
0 & \nu_2 & 0 & \nu_4 \\
\end{array}
\right) \left(
\begin{array}{c}
X_{t-1}^{(1)} \\
X_{t-1}^{(2)} \\
X_{t-1}^{(3)} \\
X_{t-1}^{(4)} \\
\end{array}
\right) + \underline{\varepsilon}_{t}
= A(t)\underline{X}_{t-1}+\underline{\varepsilon}_{t} 
\end{eqnarray*}

Below we will give the directed and undirected graphs for each
component (this is in the stationary and can be deduced from standard
results and definitions in the literature; see the definition of moralized
 k-complex, Figures 2 and 3 in \citet{andersson2001alternative}, 
 and section 3.3 in \citet{dahlhaus2003causality}), then we give the combined
nonstationary graph.

\begin{minipage}{0.4\textwidth}
For $t\in \{1,\ldots,n/2\}$ the graph is based on 
\begin{eqnarray*}
A_{1}  = \left(
\begin{array}{cccc}
\alpha_1 & 0 & 0 & 0 \\
\beta_{1} & \beta_{2} & 0 & \beta_{4} \\
0 & 0 & \gamma_{3} & 0 \\
0 & \nu_2 & 0 & \nu_4\\
\end{array}
\right) 
\end{eqnarray*}
\end{minipage}
\begin{minipage}{0.6\textwidth}
\includegraphics[scale =0.4]{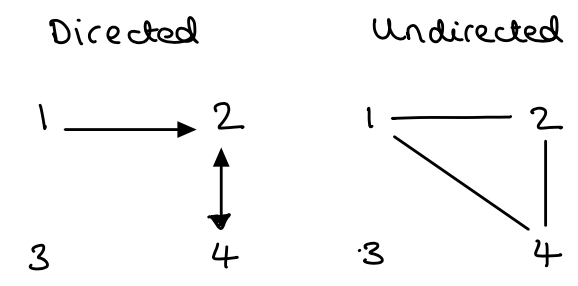}

Observe that the undirected graph adds an additional edge (1,4), 
and converts the arrows in the directed graphs to lines. This is
 the process of \textit{moralizing} a k-complex ($k=1$) in the 
 original directed graph \citep{andersson2001alternative}.
\end{minipage}

\begin{minipage}{0.4\textwidth}
For $t\in \{n/2+1,\ldots,n\}$ the graph is based on 
\begin{eqnarray*}
A_{2}  = \left(
\begin{array}{cccc}
\widetilde{\alpha}_1 & 0 & \alpha_3 & 0 \\
\beta_{1} & \beta_{2} & 0 & \beta_{4} \\
0 & 0 & \widetilde{\gamma}_{3} & 0 \\
0 & \nu_2 & 0 & \nu_4\\
\end{array}
\right) 
\end{eqnarray*}
\end{minipage}
\begin{minipage}{0.6\textwidth}
\includegraphics[scale =0.4]{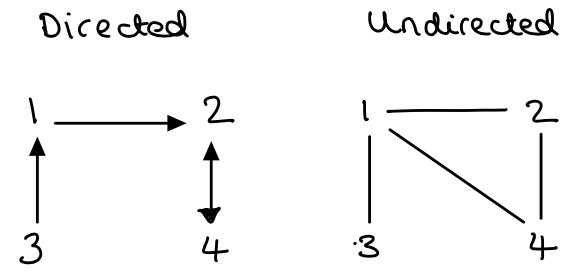}
\end{minipage}

On the other hand the nonstationary undirected graph is 

\begin{minipage}{0.5\textwidth}
For $t\in \{1,\ldots,n/2\}$ the graph is based on 
\begin{eqnarray*}
A(t)=\left(
\begin{array}{cccc}
\alpha(t) & 0 & \alpha_3(t) & 0 \\
\beta_{1} & \beta_{2} & 0 & \beta_{4} \\
0 & 0 & \gamma_{}(t) & 0 \\
0 & \nu_2 & 0 & \nu_4 \\
\end{array}
\right)
\end{eqnarray*}
\end{minipage}
\begin{minipage}{0.5\textwidth}
\includegraphics[scale =0.4]{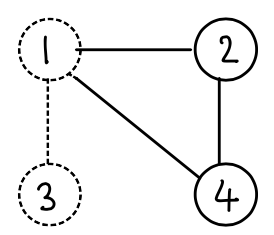}
\end{minipage}

The two stationary graphs are what we mean by local information. 
Comparing the two stationary undirected graphs with the nonstationary
directed graphs we observe that the nonstationary directed graph is
summarising all the information in the two stationary undirected
graphs. This is what we mean by global information, it tells us that
something is changing in nodes (1) and (3) (which correspond to
columns (1) and (3) in the transition matrix).

%% file: Appendix_stationary.tex
\section{Connection to graphical models for stationary time series}\label{sec:stat}

We now apply the results above to stationary multivariate time
series. This gives an alternative derivation for the
partial spectral coherency of stationary time series (see
\cite{p:bri-96} and \cite{p:dah-00b}) which is usually based on the
Wiener filter.

Suppose that $\{\underline{X}_{t}\}_{t}$ is a $p$-dimension second order stationary time
series, with spectral density matrix ${\boldsymbol \Sigma}(\omega) =
\sum_{r\in \mathbb{Z}}{\bf C}_{r}\exp(-ir\omega)$. 
By Lemmas \ref{lemma:diagonal} and \ref{lemma:inverse}
${\bf C}(\omega,\lambda)$ and ${\bf K}(\omega,\lambda)$ are diagonal
kernels where ${\bf C}(\omega,\lambda)={\boldsymbol \Sigma}(\omega) \delta_{\omega,\lambda}$ and ${\bf
  K}(\omega,\lambda) = {\boldsymbol \Gamma}(\omega)\delta_{\omega,\lambda}$
where ${\boldsymbol \Gamma}(\omega) = {\boldsymbol \Sigma}(\omega)^{-1}$. Let
$\Gamma^{(a,b)}(\omega)$ denote the $(a,b)$th entry of
${\boldsymbol \Gamma}(\omega)$. Our aim is to interprete the entries of
${\boldsymbol \Gamma}(\omega)$ in terms of the partial correlation and partial
spectral coherence. We keep in mind that since
$\{\underline{X}_{t}\}_{t}$ is second order stationary time series all
the nodes and edges of its corresponding network are conditionally stationary. 

We first interprete $\Gamma^{(a,b)}(\omega)$.
Under stationarity for all $t$ and $\tau$ and $(a,b)$ we have
\begin{eqnarray*}
 \cor\left[ X_{t}^{(a)} -
  P_{\mathcal{H}-(X_{t}^{(a)},X_{\tau}^{(b)})}(X_{t}^{(a)}), X_{\tau}^{(b)}
  -P_{\mathcal{H}-(X_{t}^{(a)},X_{\tau}^{(b)})}(X_{\tau}^{(b)})\right]
                         = \phi_{t-\tau}^{(a,b)}.
\end{eqnarray*}
Let  $\sigma_{a}^{2}= \var[X_{0}^{(a)} -
  P_{\mathcal{H}-(X_{0}^{(a)})}(X_{0}^{(a)})]$. 
By using Lemma \ref{lemma:D} for all $(t,a)\neq (\tau,b)$
\begin{eqnarray*}
[{\bf D}_{a,b}]_{t,\tau} = 
\left\{
\begin{array}{cc}
\frac{1}{\sigma_{a}^{2}} & t=\tau \textrm{ and }a=b \\
-\frac{1}{\sigma_{a}\sigma_{b}}\phi_{t-\tau}^{(a,b)} & \textrm{ otherwise }
\end{array}
\right.
\end{eqnarray*}
By using the above, we have
\begin{eqnarray*}
\Gamma^{(a,b)}(\omega)=
\left\{
\begin{array}{cc}
-\frac{1}{\sigma_{a}\sigma_{b}}\sum_{r\in
  \mathbb{Z}}\phi_{r}^{(a,b)}\exp(ir\omega) & a\neq b\\
\frac{1}{\sigma_{a}^{2}}\left(1-\sum_{r\neq
  0}\phi_{r}^{(a,a)}\exp(ir\omega)\right) &  a=b\\
\end{array}
\right..
\end{eqnarray*}
Thus the entries of $\Gamma(\omega)$ are the Fourier transforms
of the partial correlations. Let 
\begin{eqnarray*}
X_{t}^{(a)|\shortminus \{a\}} &=& X_{t}^{(a)} - P_{\mathcal{H} -
  (X^{(a)})}(X_{t}^{(a)})  \textrm{ for }t\in \mathbb{Z}. 
\end{eqnarray*}
By stationarity we have 
\begin{eqnarray*}
\rho_{t-\tau}^{(a,a)|\shortminus\{a\}} = \cov[X_{t}^{(a)|\shortminus\{a\}}, X_{\tau}^{(a)|\shortminus\{a\}}].
\end{eqnarray*}
By using Lemma \ref{lemma:spectral3} we have 
\begin{eqnarray*}
[\Gamma^{(a,a)}(\omega)]^{-1} = \sum_{r\in \mathbb{Z}}\rho_{r}^{(a,a)|\shortminus\{a\}}\exp(ir\omega).
\end{eqnarray*}

We now use the methods laid out in this paper to derive the partial
spectral coherence. For $a\neq b$ we define 
\begin{eqnarray*}
X_{t}^{(a)|\shortminus \{a,b\}} &=& X_{t}^{(a)} - P_{\mathcal{H} -
  (X^{(a)},X^{(b)})}(X_{t}^{(a)}) \textrm{ and }X_{\tau}^{(b)|\shortminus \{a,b\}} = X_{\tau}^{(b)} - P_{\mathcal{H} -
  (X^{(a)},X^{(b)})}(X_{\tau}^{(b)}).
\end{eqnarray*}
Since the time series is  stationary we  define the time series partial covariance  as
\begin{eqnarray*}
\left(
\begin{array}{cc}
\rho_{t-\tau}^{(a,a)|\shortminus\{a,b\}} & \rho_{t-\tau}^{(a,b)|\shortminus\{a,b\}} \\
\rho_{t-\tau}^{(b,a)|\shortminus\{a,b\}} &\rho_{t-\tau}^{(b,b)|\shortminus\{a,b\}}\\
\end{array}
\right) =
\cov\left[ \left(
\begin{array}{c}
X_{t}^{(a)|\shortminus\{a,b\}} \\
X_{t}^{(b)|\shortminus\{a,b\}} \\
\end{array}
 \right),
\left(
\begin{array}{c}
X_{\tau}^{(a)|\shortminus\{a,b\}} \\
X_{\tau}^{(b)|\shortminus\{a,b\}} \\
\end{array}
 \right)
\right]
\end{eqnarray*}
Thus by using Lemma \ref{lemma:spectral3} we have 
\begin{eqnarray*}
&& \sum_{r\in \mathbb{Z}}
\left(
\begin{array}{cc}
\rho_{r}^{(a,a)|\shortminus\{a,b\}} & \rho_{r}^{(a,b)|\shortminus\{a,b\}} \\
\rho_{r}^{(b,a)|\shortminus\{a,b\}} &\rho_{r}^{(b,b)|\shortminus\{a,b\}}\\
\end{array}
\right)
\exp(ir\omega) \\
&=& 
\frac{1}{\Gamma^{(a,b)}(\omega)\Gamma^{(b,b)}(\omega)  - |\Gamma^{(a,b)}(\omega)|^{2}}
\left(
\begin{array}{cc}
\Gamma^{(b,b)}(\omega) & -\Gamma^{(a,b)}(\omega) \\
-\Gamma^{(a,b)}(\omega)^{*} & \Gamma^{(a,a)}(\omega) \\
\end{array}
\right).
\end{eqnarray*}
Therefore for $c\in
\{a,b\}$ we have 
\begin{eqnarray*}
\sum_{r\in
  \mathbb{Z}}\rho_{r}^{(c,c)|\shortminus\{a,b\}}\exp(ir\omega) 
=
  \frac{\Gamma^{(c,c)}(\omega)}{\Gamma^{(a,a)}(\omega)\Gamma^{(b,b)}(\omega)
  - |\Gamma^{(a,b)}(\omega)|^{2}} =
H^{(c,c)}(\omega)
\end{eqnarray*}
and 
\begin{eqnarray*}
\sum_{r\in
  \mathbb{Z}}\rho_{r}^{(a,b)|\shortminus\{a,b\}}\exp(ir\omega) 
=  -\frac{\Gamma^{(a,b)}(\omega)}{\Gamma^{(a,a)}(\omega)\Gamma^{(b,b)}(\omega)
  - |\Gamma^{(a,b)}(\omega)|^{2}} = H^{(a,b)}(\omega). 
\end{eqnarray*}
Thus the partial spectral coherence between edge $(a,b)$ is
\begin{eqnarray}
\label{eq:Rab1}
R_{ab}(\omega) = \frac{H^{(a,b)}(\omega) }{\sqrt{H^{(a,a)}(\omega) H^{(b,b)}(\omega) }}=
-\frac{\Gamma^{(a,b)}(\omega)}{\sqrt{\Gamma^{(a,a)}(\omega)\Gamma^{(b,b)}(\omega)}}.
\end{eqnarray}

This coincides with with the partial spectral coherence given in
\cite{p:dah-00b}, equation (2.2) who shows that the 
partial spectral coherence is 
\begin{eqnarray}
\label{eq:Rab2}
R_{ab}(\omega) = \frac{g_{a,b}(\omega)}{\sqrt{g_{a,a}(\omega)g_{b,b}(\omega)}}
\end{eqnarray}
where 
\begin{eqnarray*}
g_{c,d}(\omega) = \Sigma_{c,d}(\omega)  -
  \Sigma_{c,-(a,b)}(\omega)\Sigma_{-(a,b)}(\omega)^{-1}\Sigma_{d,-(a,b)}(\omega)^{*}
  \quad c,d\in \{a,b\}
\end{eqnarray*}
and  $\Sigma_{a,-(a,b)}$ denotes the spectral cross
correlation between $\{X_{t}^{(a)}\}_{t}$ and $\{X_{t}^{(c)};c\notin
\{a,b\}\}$, $\Sigma_{b,-(a,b)}$ denotes the spectral cross
correlation  between $\{X_{t}^{(b)}\}_{t}$ and $\{X_{t}^{(c)};c\neq
\{a,b\}\}$ and $\Sigma_{-(a,b)}$ denotes the spectral cross correlation
of $\{X_{t}^{(c)};c\notin \{a,b\}\}_{t}$ i.e.
\begin{eqnarray*}
{\boldsymbol \Sigma}(\omega) = 
\left(
\begin{array}{ccc}
\Sigma_{a,a}(\omega) & \Sigma_{a,b}(\omega) & \Sigma_{a,-(a,b)}(\omega) \\
\Sigma_{b,a}(\omega) & \Sigma_{b,b}(\omega) & \Sigma_{b,-(a,b)}(\omega) \\
\Sigma_{a,-(a,b)}(\omega)^{*} & \Sigma_{b,-(a,b)}(\omega)^{*} & \Sigma_{-(a,b)}(\omega) \\
\end{array}
\right).
\end{eqnarray*}
\cite{p:dah-00b}, Theorem 2.4 shows that (\ref{eq:Rab1}) and
(\ref{eq:Rab2}) are equivalent. For completeness we give the proof 
using the block matrix inversion identity. 
The Schur complement of the $(p-2)\times(p-2)$ matrix
$\Sigma_{-(a,b)}(\omega)$ in ${\boldsymbol \Sigma}(\omega)$ is 
\begin{eqnarray*}
&& P^{(a,b)}(\omega) \\
&=&  \left(
\begin{array}{cc}
\Sigma_{aa}(\omega)  -
  \Sigma_{a,-(a,b)}(\omega)\Sigma_{-(a,b)}(\omega)^{-1}\Sigma_{-a,(a,b)}(\omega)^{*} & \Sigma_{ba}(\omega)  -
  \Sigma_{a,-(a,b)}(\omega)\Sigma_{-(a,b)}(\omega)^{-1}\Sigma_{-b,(a,b)}(\omega)^{*}  \\
\Sigma_{ab}(\omega)  -
  \Sigma_{b,-(a,b)}(\omega)\Sigma_{-(a,b)}(\omega)^{-1}\Sigma_{a,-(a,b)}(\omega)^{*} & \Sigma_{bb}(\omega)  -
  \Sigma_{b,-(a,b)}(\omega)\Sigma_{-(a,b)}(\omega)^{-1}\Sigma_{-b,(a,b)}(\omega)^{*} \\
\end{array}
\right).
\end{eqnarray*}
Using the block inverse identity we recall that 
$P^{(a,b)}(\omega)^{-1}$ is the top left hand matrix in ${\boldsymbol
  \Sigma}(\omega)^{-1} = {\boldsymbol \Gamma}(\omega)$. Thus 
\begin{eqnarray*}
P^{(a,b)}(\omega)^{-1} =  \left(
\begin{array}{cc}
\Gamma^{(a,a)}(\omega) & \Gamma^{(a,b)}(\omega)\\
{} \Gamma^{(b,a)}(\omega) & \Gamma^{(b,b)}(\omega) \\
\end{array}
\right).
\end{eqnarray*}
Therefore from the above we have
\begin{eqnarray*}
&&\left(
\begin{array}{cc}
\Sigma_{aa}(\omega)  -
  \Sigma_{a,-(a,b)}(\omega)\Sigma_{-(a,b)}(\omega)^{-1}\Sigma_{-a,(a,b)}(\omega)^{*} & \Sigma_{ba}(\omega)  -
  \Sigma_{a,-(a,b)}(\omega)\Sigma_{-(a,b)}(\omega)^{-1}\Sigma_{-b,(a,b)}(\omega)^{*}  \\
\Sigma_{ab}(\omega)  -
  \Sigma_{b,-(a,b)}(\omega)\Sigma_{-(a,b)}(\omega)^{-1}\Sigma_{a,-(a,b)}(\omega)^{*} & \Sigma_{bb}(\omega)  -
  \Sigma_{b,-(a,b)}(\omega)\Sigma_{-(a,b)}(\omega)^{-1}\Sigma_{-b,(a,b)}(\omega)^{*} \\
\end{array}
\right)\\
&=& \frac{1}{\Gamma^{(a,b)}(\omega)\Gamma^{(b,b)}(\omega)  - |\Gamma^{(a,b)}(\omega)|^{2}}
\left(
\begin{array}{cc}
\Gamma^{(b,b)}(\omega) & -\Gamma^{(a,b)}(\omega) \\
-\Gamma^{(b,a)}(\omega) & \Gamma^{(a,a)}(\omega) \\
\end{array}
\right).
\end{eqnarray*}
Comparing the entries in the above it is immediately clear that 
\begin{eqnarray*}
\frac{g_{a,b}(\omega)}{\sqrt{g_{a,a}(\omega)g_{b,b}(\omega)}} = 
-\frac{\Gamma^{(a,b)}(\omega)}{\sqrt{\Gamma^{(a,a)}(\omega)\Gamma^{(b,b)}(\omega)}}.
\end{eqnarray*}
Thus proving that (\ref{eq:Rab1}) and (\ref{eq:Rab2}) are equivalent
expression for multivariate stationary time series.